\documentclass[11pt]{amsart}

\usepackage[T1]{fontenc}
\usepackage[utf8]{inputenc}

\usepackage{hyperref}

\allowdisplaybreaks[4]
\usepackage{algorithm}
\usepackage{algpseudocode}

\usepackage{amsmath,mathtools}
\usepackage{amssymb,bbm}
\usepackage{mathrsfs}
\usepackage{amsthm}
\usepackage{bm}
\usepackage{graphicx}
\usepackage{tikz}
\usepackage{epstopdf}
\usepackage{booktabs}
\usepackage{epsfig}

\usepackage[caption=false]{subfig}
\usepackage{caption}

\usepackage{float}

\newfloat{longfigure}{htbp}{lof}
\floatname{longfigure}{Figure}
\newsubfloat{longfigure}

\usepackage{geometry}
\usepackage{color}

\usepackage{color}
\theoremstyle{plain}
\newtheorem{lemma}{Lemma}[section]
\newtheorem{theorem}[lemma]{Theorem}

\newtheorem{remark}[lemma]{Remark}
\newtheorem{proposition}[lemma]{Proposition}

 \makeatletter
\@namedef{subjclassname@2020}{%
  \textup{2020} Mathematics Subject Classification}
\makeatother
\numberwithin{equation}{section}

\begin{document}
\title[Weak Convergence analysis for SACE]
{Optimal Order Weak Error Estimates of a Full Discretization to Stochastic Allen-Cahn Equation with Nemytskii Type Diffusions}
\author{Minxing Zhang, Yongkui Zou}
\address{School of Mathematics, Jilin University, Changchun 130012, China}
\email{zhangminxingjlu@163.com, zouyk@jlu.edu.cn}

\author{Yanzhao Cao}
\thanks{Yanzhao Cao is the corresponding author.}
\address{Department of Mathematics and Statistics, Auburn University, Auburn, AL 36840, USA}
\email{yzc0009@auburn.edu}
\thanks{The work is supported by Jilin Provincial Department of Science and Technology (20240301017GX) and National Natural Science Foundation of China (12571428).}
\begin{abstract}
 	In this paper, we investigate a drift-implicit Euler scheme and spectral Galerkin approximation to a class of stochastic Allen-Cahn equations with Nemytskii type diffusion coefficients. A key difficulty arises from the fact that the second-order Fr\'echet derivative of the diffusion operator is unbounded on the square-integrable function space. We first show that this derivative becomes bounded when restricted to higher-order Sobolev product spaces. Building on this, we introduce a new strategy for estimating the solution to an associated Kolmogorov equation by analyzing a composed operator involving a positive fractional Laplacian, the second-order Fr\'echet derivative of the Kolmogorov solution, and a negative fractional Laplacian.  This allows us to prove that the spatial weak convergence order is twice the spatial regularity of the exact mild solution, yielding an optimal spatial weak rate. The temporal weak error analysis presents additional challenges due to the need to estimate the trace of an operator involving infinite-dimensional stochastic integrals, as well as to treat an anticipating stochastic integral arising from the implicit time discretization. We show that the temporal weak convergence rate is also twice the temporal H\"older continuity order of the exact mild solution.
\end{abstract}

\maketitle

\section{Introduction}
\label{Sec1}

In recent years, numerical approximations to  nonlinear stochastic partial differential equations (SPDEs) with non-globally Lipschitz continuity, such as stochastic Cahn–Hilliard equations \cite{CuiJianbo2021_2}, stochastic Burgers equations \cite{Jentzen2013} and stochastic Allen–Cahn equations \cite{ZhaoWeidong2024}, have been extensively studied to explore the statistical properties of the underlying mathematical models, where both strong and weak error estimates have been used to investigate the numerical convergence. The strong error characterizes the average deviation of a numerical trajectory from the exact one \cite{Kruse2014,Lord2013,CuiJianbo2021_2}, while the weak error measures the discrepancy in sampling the probability law of the exact solution \cite{Andersson2016,Brehier2018,WangXiaojie2016,CaiMeng2023}. These studies consistently show that weak convergence orders are higher than their strong counterparts.

The stochastic Allen-Cahn equation models thermal perturbations and  plays a key role in phase-field theory and in the simulation of rare events in infinite-dimensional stochastic systems \cite{Cerrai2003}. A great amount of efforts has been devoted to studying the strong errors \cite{Brehier2019_2,Brehier2019,Larsson2015,Kovacs2018,QiRuisheng2019,WangXiaojie2020} and weak errors \cite{CuiJianbo2019,CuiJianbo2021,WangXiaojie2021,CuiJianbo2024,JiangYingsong2026,ZhangZouCao2027JCAM} of numerical approximations to stochastic Allen-Cahn equation driven by additive noise. Meanwhile, a few works are also proposed to studying strong  errors \cite{LiuZhihui2021,Shenjie2023,ZhaoWeidong2024,ChenChuchu2025,Bauzet2025,QiXiao2026} and weak errors \cite{Breit2024,ZhangZou2026JSC} of numerical solutions to stochastic Allen-Cahn equation driven by multiplicative noise, where the key difficulties in analysis emanate from the interaction between stochastic forcing and solution. 

In this paper, we are concerned with a class of one-dimensional stochastic Allen-Cahn equations with Nemytskii operator-valued diffusion terms. We aim to study the optimal order weak error estimates of numerical solutions based on drift-implicit Euler and spectral Galerkin methods. Existing works on weak error analysis for SPDEs driven by multiplicative noise typically assume that the diffusion coefficient and its higher-order Fr\'echet derivatives are continuous and uniformly bounded \cite{Debussche2011,Jentzen2019,Naurois2021,Breit2024}. In contrast, we prove that the second order Fr\'echet derivative of the Nemytskii operator is unbounded and only becomes uniformly bounded on a more regular function space, thereby relaxing the standard assumptions on the diffusion coefficient. In addition, we show that the spatial weak convergence rate is twice the spatial regularity of the exact mild solution, yielding an optimal weak convergence order.

A main challenge of this work is to establish uniform a priori estimates for the Kolmogorov equation associated with a spatial spectral Galerkin semi-discretization. 
The key point lies in dealing with Nemytskii operator-valued diffusion terms.
To overcome the difficulties arising from the unboundedness of Nemytskii operator, we first establish uniform bounds for a composed operator involving a positive fractional power of the Laplacian, the second-order Fr\'echet derivative of the Kolmogorov equation solution, and a negative fractional power of the Laplacian. Second, we derive higher-order regularity estimates for the Malliavin derivatives of the numerical solution. These results form the basis for proving the optimal weak convergence rate.

The difficulties involved in the temporal weak error analysis arise from two sources. The first challenge is due to the implicit structure of the drift-implicit Euler scheme, which leads a non-adapted stochastic process. To overcome this challenge, we employ an It\^o's formula for Skorohod integrals \cite{CuiJianbo2021} to estimate such term. The second challenge concerns estimating the trace of an operator that contains a stochastic integral with respect to a Wiener process on a Hilbert space. A key step is to remove the stochastic integral appearing inside the trace. We utilize the spectral representations of both the $Q$-Wiener process and the directional Fr\'echet derivative of the solution, together with the expression of trace via a complete orthonormal basis to transform the trace into a multiple series of real-valued random variables involving stochastic integrals with respect to one-dimensional Wiener processes. We then apply the classical Malliavin integration-by-parts formula in the one-dimensional setting \cite{Brehier2018,Debussche2011,Nualart2006} to convert the stochastic integrals into deterministic time integrals. This finally enables us to establish the weak convergence order in time.

The rest of this paper is organized as follows. 
In Section~\ref{Sect2}, the stochastic Allen–Cahn equation, the underlying assumptions, and relevant Malliavin calculus are introduced.  
Section~\ref{Sect3} presents the spatially semi-discretized and fully discretized approximations based on the drift-implicit Euler scheme and the spectral Galerkin method, together with a priori estimates.  
In Section~\ref{Sect4}, we investigate a priori estimates for the Kolmogorov equation associated with the spatial semi-discretization.  
Section~\ref{Sect5} establishes the spatial optimal weak convergence order for the semi-discretization, and Section~\ref{Sect6} proves the weak convergence order for the fully discretized numerical solution.

\section{Preliminaries}
\label{Sect2}

\subsection{Some working spaces and properties}
Let $(U,(\cdot, \cdot)_{U})$ and $(V,(\cdot, \cdot)_{V})$ be two real separable Hilbert spaces and $(\Omega,\mathcal{F},\mathbb{P})$ be a complete probability space. By $\mathcal{L}(U;V)$ we denote the Banach space of all bounded linear operators from $U$ to $V$ endowed with operator norm. For simplicity, we write $\mathcal{L}(U) = \mathcal{L}(U;U)$.

For any $T\in\mathcal{L}(U;V)$, denote by $\{\sigma_k\}_{k\in\mathbb{N}}$ the eigenvalues of $(TT^*)^{\frac{1}{2}}$. For any $1\le p<\infty$, define a Banach space of linear operators
\begin{equation*}
	\mathcal{L}_p(U;V) = \big\{T\in\mathcal{L}(U;V):\ \|T\|_{\mathcal{L}_p(U;V)} = \big(\sum_{k=1}^{\infty}\sigma_k^p\big)^{\frac{1}{p}}<\infty\big\}.
\end{equation*}
Then, $\mathcal{L}_1(U;V)$ and $\mathcal{L}_2(U;V)$ are the spaces of trace-class operator and Hilbert-Schmidt operator, respectively.
Let $\{e_k\}_{k=1}^{\infty}$ be an orthonormal basis of $U$, then we have
\begin{equation*}
	\mathrm{Tr}(T) = \sum_{k=1}^{\infty}(Te_{k}, e_{k})_{U} \le \|T\|_{\mathcal{L}_{1}(U)},\quad \|T\|_{\mathcal{L}_{2}(U;V)} = \sqrt{\mathrm{Tr}(TT^*)}.
\end{equation*}
For simplifying notations, we write $\mathcal{L}_p(U)=\mathcal{L}_p(U;U)$. Furthermore, for any $S\in \mathcal{L}(U;V)$, $T_1,T_2 \in \mathcal{L}_2(U)$, there hold \cite{DaPrato2014,Kruse2014}
\begin{equation}\label{eq7.2.3}
	\begin{aligned}
		&\|ST_1\|_{\mathcal{L}(U;V)}\le\|ST_1\|_{\mathcal{L}_2(U;V)}\le\|S\|_{\mathcal{L}(U;V)}\|T_1\|_{\mathcal{L}_2(U)},\\ &\|T_1T_2\|_{\mathcal{L}_2(U)}\le\|T_1T_2\|_{\mathcal{L}_1(U)}\le\|T_1\|_{\mathcal{L}_2(U)}\|T_2\|_{\mathcal{L}_2(U)}.
	\end{aligned}
\end{equation}

For any integer $m\ge 1$, by $\mathcal{C}_b^m(U)$ denote the space of not necessarily bounded functions $\phi: U \rightarrow \mathbb{R}$ having continuous and bounded Fr\'echet derivatives $D^i\phi(u)$ $(1\le i\le m)$ equipped with a maximum seminorm $|\phi|_{\mathcal{C}_b^i(U)}$.
Notice that the first order Fr\'echet derivative $D\phi(u)$ of $\phi$ at $u\in U$ is a bounded linear functional on $U$. According to Riesz representation theorem, there exists a unique element in $U$, still denoted by $D\phi(u)$, such that
\begin{equation}\label{eq7.2.4}
	D\phi(u)h = ( D\phi(u),h)_U,\quad |D\phi(u)h| \le |\phi|_{\mathcal{C}_b^1(U)}\|h\|_U, \quad \forall u,h \in U.
\end{equation}

Let $\mathcal{D}=[0,1]$. Denote by $W^{k,p}$ $(k\ge 0,p\ge 1)$ the standard Sobolev spaces defined on $\mathcal{D}$, and simply write $H^k=W^{k,2}$. Especially, $H = W^{0,2}(\mathcal{D}) = L^2(\mathcal{D})$. The Sobolev embedding $W^{s,p}\hookrightarrow W^{r,q}$ holds for $s\ge r$ and $s-\frac{1}{p}\ge r-\frac{1}{q}$, cf. \cite{Brezis2001}. Let $V = L^{\infty}(\mathcal{D})$ be the space of essentially bounded functions with norm $\|\cdot\|_V$.

Let $A =-\Delta:\mathrm{Dom}(A) = H_0^1 \cap H^2 \subset H \to H$ be a Laplacian operator with homogeneous Dirichlet boundary condition. 
Then, the operator $A$ possesses a family of eigenvalues $\beta_i = \pi^2i^2$ with complete orthonormal eigenbases $\varphi_i(x)=\sqrt{2}\sin(i\pi x)$. 
For any $\alpha \in \mathbb{R}$, define a fractional power of $A$ by $\displaystyle A^{\alpha}v = \sum_{j=1}^{\infty}\beta_j^\alpha(v,\varphi_j)\varphi_j$ and a space $\displaystyle\dot{H}^{\alpha} = \{v \in H: \| v \|_{\alpha} =\| A^{\frac{\alpha}{2}}v\| < \infty\}$. Then $(\dot{H}^{\alpha},\| \cdot \|_{\alpha})$ is a Hilbert space. Especially, we have $\dot{H}^0 = H$ and $\dot{H}^1 = H_0^1$, cf.\ \cite{Thomee2006}.

The operator $-A$ generates an analytic semigroup $S(t) = e^{-tA}$ on $H$. It is pointed out in \cite{Thomee2006} that for any $\alpha \ge 0$, there exists a constant $C>0$ such that
\begin{equation}\label{eq3.2.2}
	A^{\alpha}S(t) = S(t)A^{\alpha} \ \text{ on } \dot{H}^{2\alpha}  \text{ and } 	\|A^{\alpha}S(t)\|_{\mathcal{L}(H)} \le Ct^{-\alpha}, \quad t\in(0,T].
\end{equation}

Let $Q$ be a trace-class, self-adjoint and strictly positive definite operator on $H$ with eigenpairs $(q_j,\varphi_j)$ $(j=1,2,\cdots)$. Then, $Q$ commutes with $A$ and $Q^{\frac{1}{2}}$ is a Hilbert-Schmidt operator due to
\begin{equation}\label{eq3.2.00}
	\|Q^{\frac{1}{2}}\|^2_{\mathcal{L}_2(H)}=\mathrm{Tr}(Q)=\sum_{j=1}^{\infty}q_j<\infty.
\end{equation}
Define a Hilbert space $U_0 = Q^{\frac{1}{2}}(H) \subset H$ with scalar product $(u,v)_{U_0}=(Q^{-\frac{1}{2}}u,Q^{-\frac{1}{2}}v)$ and norm $\|u\|_{U_0}=\|Q^{-\frac{1}{2}}u\|$. Then, $\displaystyle\{Q^{\frac{1}{2}}\varphi_j\}_{j=1}^{\infty}$ is a complete orthogonal basis of $U_0$, cf.\ \cite[Page 96]{DaPrato2014}. For any $\alpha \ge 0$, denote by $\mathcal{L}_2^\alpha = \mathcal{L}_2(U_0,\dot{H}^{\alpha})$ the space of Hilbert-Schmidt operators from $U_0$ to $\dot{H}^{\alpha}$. We emphasize that $\mathcal{L}_2^{\alpha}$ is a separable Hilbert space equipped with norm \cite{Shenjie2023}
\begin{equation}\label{eq3.2.norm}
	\|T\|^2_{\mathcal{L}_2^{\alpha}} = \|A^{\frac{\alpha}{2}}T\|^2_{\mathcal{L}_2^0} =  \|A^{\frac{\alpha}{2}}TQ^{\frac{1}{2}}\|^2_{\mathcal{L}_2(H)} = \sum_{k=1}^{\infty}\|TQ^{\frac{1}{2}}\varphi_k\|_{\alpha}^2.
\end{equation}

To end this section, we present a Burkholder-Davis-Gundy inequality for a sequence of $H$-valued discrete martingales.

\begin{lemma}\label{l3.2.00}\cite[Lemma 2.1]{Shenjie2023}
	Let $p\ge2$, $N\ge 2$ and $\{Z_m\}$ be a sequence of $H$-valued random variables with bounded $p$-moments such that $\mathbb{E}[Z_{m+1}|Z_0,\cdots,Z_m]=0$ for all $1\le m\le N-1$. Then there exists a constant $C=C(p)>0$ such that
	\begin{equation*}
		\big(\mathbb{E}\big\|\sum_{i=0}^mZ_i\big\|^p\big)^{\frac{1}{p}} \le C\big(\sum_{i=0}^m(\mathbb{E}\|Z_i\|^p)^{\frac{2}{p}}\big)^{\frac{1}{2}}.
	\end{equation*}
\end{lemma}

\subsection{Stochastic Allen-Cahn equation and assumptions}
In this section, we first describe a one-dimensional stochastic Allen--Cahn equation driven by multiplicative trace-class noise. Then we introduce some assumptions and collect some properties of its solution.

Fix $T>0$ and consider a stochastic Allen--Cahn equation driven by multiplicative noise
\begin{equation}\label{eq3.1.1}
	\begin{aligned}
		dX(t)+A X(t)dt =&\ F(X(t))dt+G(X(t))dW(t),\ t\in (0,T],\\
		\quad X(0) =&\ \xi,
	\end{aligned}
\end{equation}
where $F: L^{6}(\mathcal{D}) \rightarrow H$ and $G:H\to\mathcal{L}(H)$ are Nemytskii operators defined by, respectively, 
\begin{align*}
	&F(u)(x)=f(u(x))=u(x)-u^{3}(x),\  \forall x \in \mathcal{D}, u \in L^{6}(\mathcal{D}),\\
	&(G(u)v)(x) = g(u(x))v(x),\ \forall x \in \mathcal{D},\ \forall u,v \in H,\ g\in \mathcal{C}_b^3(\mathbb{R}),
\end{align*} 
$W(t)$ is an $H$-valued $Q$-Wiener process, and $\xi$ is an initial value function. 
By $X(t,\xi)$ we denote the solution of \eqref{eq3.1.1} satisfying $X(0,\xi)(x) = \xi(x)$. 

We assume

\begin{itemize}
	\item[(A1)] The initial value $\xi \in \dot{H}^{\gamma}$ is a deterministic function for a fixed $\gamma\in(\frac{1}{2},1)$.
	\item[(A2)] $g\in \mathcal{C}_b^3(\mathbb{R})$ with $g(0)=0$.
	\item[(A3)] There exists a $\rho>2\gamma+1$ such that the eigenvalues of $Q$ satisfy $q_j\le C j^{-\rho}$ for all $j$, where $\gamma\in(\frac{1}{2},1)$, and $C>0$ is a constant.
\end{itemize}

According to \cite{WangXiaojie2021,Larsson2015}, there exists a constant $C_F>0$ such that for all $u,v\in V$ and $w_1,w_2 \in L^{6}(\mathcal{D})$,
\begin{equation}\label{eq3.2.5}
	\begin{aligned}
		& \|F(u)\|\le C_F(1+\|u\|_V^2)\|u\|,\quad \|F(u)-F(v)\| \leq C_F(1+\|u\|_{V}^{2}+\|v\|_{V}^{2})\|u-v\|, \\
		& (F(u)-F(v),u-v) \le C_F\|u-v\|^2,\\
		& ( DF(w_2)w_1,w_1) \le C_F\|w_1\|^2,\quad \|DF(u)w_1\| \leq C_F(1+\|u\|_{V}^{2})\|w_1\|, \\
		& \|D^2F(u)(w_1,w_2)\|_{-1} \leq C_F\|u\|_V\|w_1\|\|w_2\|.
	\end{aligned}
\end{equation}

The next lemma studies properties of $G$.

\begin{lemma}\label{l3.2.1}
	Assume (A1)-(A3). Then there exists a constant $C>0$ such that for any $u\in H$
	\begin{align*}
		&\|G(u)\|_{\mathcal{L}_2^0} \le C(1+\|u\|),\\
		&\|DG(u)v_1\|_{\mathcal{L}_2^0} \le C\|v_1\|,\quad \forall v_1\in H,\\ 
		&\|D^2G(u)(v_1,v_2)\|_{\mathcal{L}_2^0} \le C\|v_1\|_{L^4(\mathcal{D})}\|v_2\|_{L^4(\mathcal{D})},\quad \forall v_1,v_2\in L^4(\mathcal{D}),\\
		&\|D^2G(u)(v_1,v_2)\|_{\mathcal{L}_2^0} \le C\|v_1\|\|v_2\|_V,\quad \forall v_1\in H,v_2\in V.
	\end{align*}
	Moreover, $G(\dot{H}^{\gamma})\subset \mathcal{L}_2^\gamma$ and
	\begin{equation*}
		\|G(u)\|_{\mathcal{L}_2^\gamma} \le C(1+\|u\|_\gamma),\quad \forall u\in \dot{H}^{\gamma}.
	\end{equation*}
\end{lemma}

\begin{proof}
	Noticing that $Q$ is trace-class, for any $u\in H$ we get
	\begin{align*}
		\|G(u)\|_{\mathcal{L}_2^0}^2 =&\ \sum_{j=1}^{\infty}\|G(u)Q^{\frac{1}{2}}\varphi_j\|^2 = \sum_{j=1}^{\infty}\int_\mathcal{D}|g(u(x))\sqrt{q_j}\varphi_j(x)|^2dx\\ 
		\le&\ \sum_{j=1}^{\infty}q_j \int_\mathcal{D}|g|_{\mathcal{C}_b^1}^2(1+|u(x)|)^2\cdot 2dx \le C(1+\|u\|^2).
	\end{align*}
	
	For any $v_1\in H$, we have
	\begin{align*}
		\|DG(u)v_1\|_{\mathcal{L}_2^0}^2 =&\ \sum_{j=1}^{\infty}\|DG(u)v_1Q^{\frac{1}{2}}\varphi_j\|^2 = \sum_{j=1}^{\infty}\int_\mathcal{D}|g'(u(x))v_1(x)\sqrt{q_j}\varphi_j(x)|^2dx\\ 
		\le&\ |g|_{\mathcal{C}_b^1}^2\sum_{j=1}^{\infty}q_j\int_\mathcal{D}|v_1(x)|^2\cdot2dx \le  C\|v_1\|^2.
	\end{align*}
	
	Furthermore, for any $v_1,v_2\in L^4(\mathcal{D})$, by H\"older inequality we obtain
	\begin{align*}
		&\ \|D^2G(u)(v_1,v_2)\|_{\mathcal{L}_2^0}^2 = \sum_{j=1}^{\infty}\|D^2G(u)(v_1,v_2)Q^{\frac{1}{2}}\varphi_j\|^2\\ 
		=&\ \sum_{j=1}^{\infty}\int_\mathcal{D}|g''(u(x))v_1(x)v_2(x)\sqrt{q_j}\varphi_j(x)|^2dx\\ 
		\le&\  |g|_{\mathcal{C}_b^2}^2\sum_{j=1}^{\infty}q_j\int_\mathcal{D}|v_1(x)v_2(x)|^2\cdot2dx 
		\le  C\|v_1\|^2_{L^4(\mathcal{D})}\|v_2\|^2_{L^4(\mathcal{D})}.
	\end{align*}
	
	For any $v_1\in H$ and $v_2\in V$, we have
	\begin{equation*}
		\|D^2G(u)(v_1,v_2)\|_{\mathcal{L}_2^0}^2 \le |g|_{\mathcal{C}_b^2}^2\sum_{j=1}^{\infty}q_j\int_\mathcal{D}|v_1(x)v_2(x)|^2\cdot2dx \le C\|v_1\|^2\|v_2\|_{V}^2.
	\end{equation*}
	
	For any $\eta\in[\gamma,\min\{\frac{\rho-1}{2},1\})$, there holds $2\eta-\rho<-1$. Note that $a\le 2^{1-\eta}a^{\eta}$ for any $a\in[0,2]$, then
	\begin{align*}
		\sum_{j=1}^{\infty}q_j\|\varphi_j\|^2_{\mathcal{C}^\eta(0,1)} \le&\ 2\sum_{j=1}^{\infty}C j^{-\rho}\Big(1+\sup _{\substack{x,y\in(0,1) \\ x \neq y}}\frac{|\sin(j\pi x)-\sin(j\pi y)|}{|x-y|^{\eta}}\Big)^2\\
		\le&\ 2\sum_{j=1}^{\infty}C j^{-\rho}\Big(1+\sup _{\substack{x,y\in(0,1) \\ x \neq y}}\frac{2^{1-\eta}\cdot|\sin(j\pi x)-\sin(j\pi y)|^{\eta}}{|x-y|^{\eta}}\Big)^2\\
		\le&\ 2\sum_{j=1}^{\infty}C j^{-\rho}\Big(1+2^{1-\eta}(j\pi)^{\eta}\Big)^2 \le 2\sum_{j=1}^{\infty}\nu j^{-\rho}(1+\pi j^{\eta})^2\\
		\le&\ 8C \pi^2(\sum_{j=1}^{\infty}j^{2\eta-\rho})<\infty.
	\end{align*}
	Then, according to \cite[(27)]{Jentzen2012}, we have for any $\alpha\in(0,\eta)$, $u\in W^{\alpha,2}(0,1)$ and $v\in U_0$,
	\begin{equation*}
		\|G(u)v\|_{W^{\alpha,2}(0,1)} \le \frac{C}{(\eta-\alpha)^2}\Big(\sum_{j=1}^{\infty}q_j\|\varphi_j\|^2_{\mathcal{C}^\eta(0,1)}\Big)^{\frac{1}{2}}(1+\|u\|_{W^{\alpha,2}(0,1)})\|v\|_{U_0}.
	\end{equation*}
	The last assertion follows from taking $\eta=\gamma$.
\end{proof}

The next lemma collects some properties of solution to \eqref{eq3.1.1}.

\begin{lemma}\label{l3.2.3}\cite[Theorem 2.1, Lemma 5.4]{Shenjie2023}
	Assume (A1)-(A3). Then, \eqref{eq3.1.1} possesses a unique mild solution in $\dot{H}^{\gamma}$ adapted to $W(t)$ of the form
	\begin{equation*}
		X(t) = S(t)\xi + \int_0^tS(t-s)F(X(s))ds + \int_0^tS(t-s)G(X(s))dW(s).
	\end{equation*}
	Moreover, for any $p \ge 1$, there exists a constant $C>0$, such that
	\begin{align*}
		\mathbb{E}\| X(t)\|_{\gamma}^p \le&\ C(1+\|\xi\|_{\gamma}^{3p}),\quad \forall t\in[0,T],\\
		\|X(t)-X(s)\|_{L^p(\Omega;H)} \le&\ C(t-s)^{\frac{\gamma}{2}},\quad \forall 0\le s\le t\le T.
	\end{align*}
\end{lemma}

By Sobolev embedding property $\dot{H}^{\gamma} \hookrightarrow V$, for any $p\ge1$, there exists a constant $C>0$ such that
\begin{equation}\label{eq3.2.11.1}
	\mathbb{E}\| X(t)\|_{V}^p \le C\mathbb{E}\| X(t)\|_{\gamma}^p \le C(1+\|\xi\|_{\gamma}^{3p}),\quad \forall t\in[0,T].
\end{equation}
Since $\dot{H}^{\gamma}$ is a Banach algebra \cite[Appendix (54c)]{Zeidler1990}, there exists a constant $C>0$ such that
\begin{equation}\label{eq3.2.14}
	\mathbb{E}\|F(X(t))\|_{V}^p \le C\mathbb{E}\|X(t)-X^3(t)\|_{\gamma}^p \le C\big(\mathbb{E}\|X(t)\|_{\gamma}^p+\mathbb{E}\|X(t)\|_{\gamma}^{3p}\big) \le C.
\end{equation}

\subsection{Malliavin calculus}
In this part, we introduce some properties of Malliavin calculus, including the definition of Malliavin derivatives, integration by parts formula and Skorohod integral, cf. \cite{Andersson2016,Nualart1998}. 

For any $h \in$ $L^{2}((0,T);U_0)$, define $W(h)$ as a centered isonormal Gaussian process with covariance structure
\begin{equation*}
	\mathbb{E}[W(h)W(g)]=(h, g)_{L^{2}((0,T);U_0)}, \quad \forall h,g \in L^{2}((0,T);U_0).
\end{equation*}
For any integer $n \ge 1$, let $\mathcal{C}_{p}^{\infty}(\mathbb{R}^{n})$ denote the space of all real-valued $\mathcal{C}^{\infty}$-functions defined on $\mathbb{R}^{n}$ with polynomial growth. We define two families of smooth cylindrical random variables with values in $\mathbb{R}$ and in $H$, respectively, by
\begin{equation*}
	\begin{aligned}
		\mathcal{S}=&\ \{X=f(W(h_{1}), \cdots, W(h_{n})):  f \in \mathcal{C}_{p}^{\infty}(\mathbb{R}^{n}),\\
		&\qquad h_{1}, \cdots, h_{n} \in L^{2}((0,T);U_0),\ n=1,2,\cdots\},\\
		\mathcal{S}(H)=&\ \{F=\sum_{i=1}^{m} X_{i} \varphi_{i}: X_{i} \in \mathcal{S}, 1\le i \le m, m=1,2,\cdots\}.
	\end{aligned}
\end{equation*}
For any $X \in \mathcal{S}$, define its Malliavin derivative by 
\begin{equation*}
	\mathcal{D} X=\sum_{i=1}^{n} \partial_{i} f(W(h_{1}),\cdots, W(h_{n})) h_{i},
\end{equation*}
which is an $L^{2}((0,T);U_0)$-valued random variable. Equivalently, $(\mathcal{D}_t X)_{t\in(0,T)}$ is an $H$-valued stochastic process. Then, for any $F \in \mathcal{S}(H)$, define its Malliavin derivative by
\begin{equation*}
	\mathcal{D} F=\sum_{i=1}^{m} \sum_{j=1}^{n} \partial_{j} f_{i}(W(h_{1}), \cdots, W(h_{n}))(\varphi_{i} \otimes h_{j}).
\end{equation*}
Note that
\begin{equation*}
	H \otimes L^{2}((0,T);U_0) \simeq \mathcal{L}_2(L^{2}((0,T);U_0);H) \simeq L^2((0,T);\mathcal{L}_2^0),
\end{equation*}
then $(\mathcal{D}_{t}F)_{t \in(0,T)}$ can be regarded as an $\mathcal{L}_2^0$-valued stochastic process.
By $\mathcal{D}_{t}^{u} F$ we denote the Malliavin derivative of $F$ in direction $u \in U_0$ at time $t$, i.e.
\begin{equation*}
	\mathcal{D}_{t}^{u}F = \mathcal{D}_{t}Fu=\sum_{i=1}^{m} \sum_{j=1}^{n} \partial_{j} f_{i}(W(h_{1}), \cdots, W(h_{n})) ( u, h_{j}(t))_{U_0} \varphi_{i}.
\end{equation*}
Define a Watanabe-Sobolev space $\mathbb{D}^{1,2}$ as the closure of $\mathcal{S}(H)$ with respect to a norm
\begin{equation*}
	\|F\|^2_{\mathbb{D}^{1,2}}=\mathbb{E}\|F\|^{2}+\mathbb{E}\int_{0}^{T}\|\mathcal{D}_{t} F\|_{\mathcal{L}_2^0}^{2} dt<\infty.
\end{equation*}

Denote by $\mathrm{dom}(\delta)$ the set of elements $\Phi \in L^{2}((0,T) \times \Omega; \mathcal{L}_2^0)$ satisfying
\begin{equation*}
	\mathbb{E}[(\mathcal{D}F,\Phi)_{L^2((0,T);\mathcal{L}_2^0)}] \le C\|F\|_{L^2(\Omega;H)},\quad \forall F\in L^2(\Omega;H),
\end{equation*}
where $C=C(\Phi)>0$ is a constant.
The adjoint of Malliavin derivative $\mathcal{D}$ is an operator $\delta:\mathrm{dom}(\delta)\subset L^{2}((0,T) \times \Omega; \mathcal{L}_2^0)\to L^2(\Omega;H)$ satisfying
\begin{equation*}
	\mathbb{E}[(\mathcal{D}F,\Phi)_{L^2((0,T);\mathcal{L}_2^0)}] = \mathbb{E}[(F,\delta(\Phi))],\quad \forall F\in \mathbb{D}^{1,2},\ \Phi\in \mathrm{dom}(\delta).
\end{equation*}
For any predictable process $\Phi\in L^2((0,T);\mathcal{L}_2^0)$, the action of $\delta$ on $\Phi$ coincides with that of It\^o integral, i.e.\
\begin{equation*}
	\delta(\Phi) = \int_0^T\Phi(t)dW(t).
\end{equation*}
Then $\delta$ is an extension of It\^o integral and is called the Skorohod integral.

Let $\mathbb{L}^{i,2}(H)$ $(i=1,2)$ be the space of processes $u\in L^2((0,T)\times \Omega;H)$ such that $u(t)\in \mathbb{D}^{i,2}(H)$ $a.s.$ Then $\mathbb{L}^{i,2}(H)\simeq L^2((0,T);\mathbb{D}^{i,2}(H))$. Let $X\in \mathbb{L}^{1,2}(H)$ be of the form
\begin{equation*}
	X_t = X_0 + \int_0^tv_sds + \int_0^tu_sdW(s),
\end{equation*}
where $X_0\in \mathbb{D}^{1,2}(H)$, $v\in\mathbb{L}^{1,2}(H)$ and $u\in \mathbb{L}^{2,2}(H)$.
Denote by $(\mathcal{D}^-X)_t$ the element in $L^2((0,T)\times\Omega;\mathcal{L}_2^0)$ satisfying
\begin{equation*}
	\lim_{n\to\infty}\int_0^T\sup_{(s-\frac{1}{n})\vee 0\le t<s}\mathbb{E}\|\mathcal{D}_sX_t-(\mathcal{D}^-X)_s\|^2_{\mathcal{L}_2^0}ds = 0.
\end{equation*} 
According to \cite[Proposition 3.1.1]{Nualart2006}, we have
\begin{equation}\label{eq6.2.13}
	(\mathcal{D}^-X)_t = \mathcal{D}_tX_0 + \int_0^t\mathcal{D}_tv_rdr + \int_0^t\mathcal{D}_tu_rdW(r).
\end{equation}
Let $U:H\to \mathbb{R}$ be a twice continuously differentiable function. Then, there holds the It\^o's formula for Skorohod integral \cite[Theorem 3.2.2]{Nualart2006},
\begin{equation}\label{eq6.2.14}
	\begin{aligned}
		U(X_t) =&\ U(X_0) + \int_0^t(DU(X_s),dX_s) + \frac{1}{2}\int_0^t\mathrm{Tr}\{D^2U(X_s)u_sQ^{\frac{1}{2}}(u_sQ^{\frac{1}{2}})^*\}ds\\
		&\ + \int_0^t\mathrm{Tr}\{D^2U(X_s)(\mathcal{D}^-X)_tQ^{\frac{1}{2}}(u_sQ^{\frac{1}{2}})^*\}ds.
	\end{aligned}
\end{equation}

To end this section, we introduce the integration by parts formula in Malliavin sense, which plays an important role in weak error analysis.

\begin{lemma}\label{l3.2.4}\cite[Lemma 2.2]{Andersson2016}
	For any given random variable $F \in \mathbb{D}^{1,2}(H)$ and any $\Phi \in L^{2}((0,T) \times \Omega, \mathcal{L}_2^0)$ adapted to $W(t)$, there holds
	\begin{equation*}
		\mathbb{E}\big(F,\int_0^T\Phi(t)dW(t)\big) = \mathbb{E}\int_0^T\big(\mathcal{D}_tF,\Phi(t)\big)_{\mathcal{L}_2^0}dt.
	\end{equation*}
\end{lemma}

\section{Full-discretization and regularity estimates}
\label{Sect3}

In this section, we contruct a full-discretized approximation with drift-implicit Euler scheme and spectral Galerkin method, and then derive a priori estimate to the numerical solution.

\subsection{Spectral Galerkin semi-discretization}
For any $N\ge 1$, define an $N$-dimensional subspace $H_N = \mathrm{Span}\{\varphi_1,\cdots,\varphi_N\}$. Let $P_N:H \to H_N$ be an orthogonal projection and then there holds $P_NA=AP_N$, cf.\ \cite{CaiMeng2023}.
Define $A_N = AP_N:H_N\to H_N$, which generates an analytic semigroup $S_N(t)=e^{-tA_N}$ on $H_N$. According to \cite[(3.1)]{WangXiaojie2021}, for any $\alpha\ge 0$ and $p\ge 1$, there exists a constant $C>0$ independent of $N$ such that 
\begin{equation}\label{eq3.3.1}
	\|(P_N-I)A^{-\alpha}\|_{\mathcal{L}(H)}\le CN^{-2\alpha}.
\end{equation}
Moreover, for $\alpha\in[0,1]$, there holds, cf.\ \cite[Lemma 3.1]{Shenjie2023},
\begin{equation}\label{eq3.3.3}
	\|A_N^{\alpha}S_N(t)\|_{\mathcal{L}(H)} \le Ct^{-\alpha}.
\end{equation}

The semi-discretized spectral Galerkin approximation to \eqref{eq3.1.1} is defined as: seeking $X^{N}(t)\in H_N$ such that $X^{N}(0) = P_N \xi$ and
\begin{equation}\label{eq3.3.5}
	dX^{N}(t) + A_{N}X^{N}(t)dt = P_NF(X^{N}(t))dt + P_NG(X^{N}(t))dW(t),\ t\in (0,T].
\end{equation}
Its mild solution $X^{N}(t)=X^{N}(t,P_N\xi)$ with $X^{N}(0,P_N\xi) = P_N\xi$ is given by
\begin{equation}\label{eq3.3.6}
	X^{N}(t) = S_N(t)P_N\xi +\! \int_0^t\!S_N(t-s)P_NF(X^{N}(s))ds +\! \int_0^t\!S_N(t-s)P_NG(X^{N}(s))dW(s).
\end{equation}
Moreover, for any $p \ge 1$, there exists a constant $C>0$ independent of $N$ such that for all $t \in [0,T]$, cf.\ \cite[Lemma 3.2]{Shenjie2023},
\begin{equation}\label{eq3.3.2}
	\mathbb{E}\| X^{N}(t)\|_{\gamma}^p +
	\mathbb{E}\| X^{N}(t)\|_{V}^p \le C(1+\|\xi\|_{\gamma}^{3p}).
\end{equation}

The next lemma presents a strong convergence result of the semi-discretized spectral Galerkin approximation.

\begin{lemma}\label{l3.3.2}\cite{Shenjie2023}
	Assume (A1)-(A3). Then there exists a constant $C>0$ independent of $N$, such that
	\begin{equation*}
		\big(\mathbb{E}\|X(t)-X^N(t)\|^2\big)^{\frac{1}{2}} \le CN^{-\gamma},\quad \forall t\in[0,T].
	\end{equation*}
\end{lemma}

\subsection{Full-discretization and a priori estimates}
For any integer $K> 0$, let $\tau=\frac{T}{K}$ be a time stepsize and $t_k=k\tau$ $(k=0,1,\cdots,K)$ be the partition nodes in $[0,T]$. We apply a drift-implicit Euler scheme to \eqref{eq3.3.5} and set up a fully discretized scheme: seeking $X^N_{k}\in H_N$ $(0\le k\le K)$ such that $X_0^N = P_N\xi$ and
\begin{equation}\label{eq3.3.7}
	X_{k+1}^N = X_k^N-\tau A_NX_{k+1}^N+\tau P_NF(X_{k+1}^N)+P_NG(X_k^N)\Delta W_k,
\end{equation}
where $\Delta W_k=W(t_{k+1})-W(t_k)$. 

Define $S_{N,\tau} = (I+\tau A_N)^{-1}$, then
\begin{equation}\label{eq3.3.3.1}
	S_{N,\tau}=I-\tau A_NS_{N,\tau}.
\end{equation}
According to \cite[Lemma 4.2]{QiRuisheng2019}, for any $\alpha\in[0,1]$, there exists a constant $C>0$ independent of both $\tau$ and $N$, such that
\begin{equation}\label{eq3.3.4}
	\|S_{N,\tau}P_N\|_{\mathcal{L}(H)} \le C,\quad \|A_N^{\alpha}S_{N,\tau}^kP_N\|_{\mathcal{L}(H)} \le Ct_{k}^{-\alpha}.
\end{equation}
By \eqref{eq3.3.3.1}, \eqref{eq3.3.7} can be rewritten as
\begin{equation}\label{eq3.3.8}
	\begin{aligned}
		X_{k+1}^N =&\ S_{N,\tau}X_k^N +  \tau S_{N,\tau}P_NF(X_{k+1}^N) + S_{N,\tau}P_NG(X_k^N)\Delta W_k\\
		=&\ X_k^N-\tau A_NS_{N,\tau}X_k^N +  \tau S_{N,\tau}P_NF(X_{k+1}^N) + S_{N,\tau}P_NG(X_k^N)\Delta W_k.
	\end{aligned}
\end{equation}
Define a continuous version $\widetilde{X}(t)$ of $X_{k}^N$ for any $0\le k\le K-1$ and $t \in [t_{k},t_{k+1}]$,
\begin{equation}\label{eq3.3.9}
	\widetilde{X}(t) = X_k^N+ \int_{t_{k}}^{t}-A_NS_{N,\tau}X_{k}^N+ S_{N,\tau}P_NF(X_{k+1}^N)dr + \int_{t_{k}}^{t}S_{N,\tau}P_NG(X_k^N)dW(r).
\end{equation}
Obviously, $\widetilde{X}(t_k)$ coincides with $X_k^N$ for all $k$. 
\eqref{eq3.3.8} minus \eqref{eq3.3.9} yields
\begin{equation}\label{eq3.3.9.1}
	\begin{aligned}
		X_{k+1}^N - \widetilde{X}(t) =&\ -(t_{k+1}-t)A_NS_{N,\tau}X_{k}^N + (t_{k+1}-t)S_{N,\tau}P_NF(X_{k+1}^N)\\ 
		&\ + \int_{t}^{t_{k+1}}S_{N,\tau}P_NG(X_k^N)dW(r).
	\end{aligned}
\end{equation}
By induction on $k\ge1$ in \eqref{eq3.3.8}, we obtain
\begin{equation}\label{eq3.3.10}
	X_{k+1}^N = S_{N,\tau}^{k+1}P_N\xi + \tau\sum_{i=0}^{k}S_{N,\tau}^{k+1-i}P_NF(X_{i+1}^N) + \sum_{i=0}^{k}S_{N,\tau}^{k+1-i}P_NG(X_i^N)\Delta W_i.
\end{equation}

The next lemma studies a priori estimates for the numerical solution $X_k^N$.

\begin{lemma}\label{l3.4.2}
	Assume (A1)-(A3). Then, for any integer $p\ge1$, there exists a constant $C>0$ independent of both $\tau$ and $N$, such that
	\begin{equation*}
		\|X_k^N\|_{L^{2p}(\Omega;H)} + \|X_k^N\|_{L^{2p}(\Omega;\dot{H}^{\gamma})} + \|X_k^N\|_{L^{2p}(\Omega;V)} \le C,\quad \forall \ 0\le k\le K.
	\end{equation*}
\end{lemma}

\begin{proof}
	For any integer $p\ge1$, we first estimate $X_k^N$ in $L^{2p}(\Omega;H)$-norm. Testing \eqref{eq3.3.7} by $X_{k+1}^N$, we obtain
	\begin{equation}\label{eq3.3.13}
		(X_{k+1}^N-X_k^N,X_{k+1}^N)+\tau(A_NX_{k+1}^N,X_{k+1}^N) 
		= \tau(P_NF(X_{k+1}^N),X_{k+1}^N) + (P_NG(X_k^N)\Delta W_k,X_{k+1}^N).
	\end{equation}
	Noticing the fact  $(u-v,u) = \frac{1}{2}(\|u\|^2-\|v\|^2)+\frac{1}{2}\|u-v\|^2$, we get 
	\begin{equation*}
		\begin{aligned}
			&\ (X_{k+1}^N-X_k^N,X_{k+1}^N) = (X_{k+1}^N-X_k^N,X_{k+1}^N)\\
			=&\ \frac{1}{2}(\|X_{k+1}^N\|^2-\|X_k^N\|^2) + \frac{1}{2}\|X_{k+1}^N-X_k^N\|^2.
		\end{aligned}
	\end{equation*}
	We notice that $(A_NX_{k+1}^N,X_{k+1}^N) = (A_N^{\frac{1}{2}}X_{k+1}^N,A_N^{\frac{1}{2}}X_{k+1}^N) = \|A_N^{\frac{1}{2}}X_{k+1}^N\|^2 \ge 0$.
	From \eqref{eq3.2.5}, it follows that
	\begin{equation*}
		(P_NF(X_{k+1}^N),X_{k+1}^N) = (F(X_{k+1}^N),X_{k+1}^N) \le C_F\|X_{k+1}^N\|^2.
	\end{equation*}
	Moreover
	\begin{equation*}
		\begin{aligned}
			&\ (P_NG(X_k^N)\Delta W_k,X_{k+1}^N) = (G(X_k^N)\Delta W_k,X_{k+1}^N)\\
			=&\ (G(X_k^N)\Delta W_k,X_{k+1}^N-X_k^N) + (G(X_k^N)\Delta W_k,X_k^N)\\
			\le&\ \frac{1}{2}\|G(X_k^N)\Delta W_k\|^2 + \frac{1}{2}\|X_{k+1}^N-X_k^N\|^2 + (G(X_k^N)\Delta W_k,X_k^N).
		\end{aligned}
	\end{equation*}
	Substituting these into \eqref{eq3.3.13} yields
	\begin{equation*}
		\frac{1}{2}(\|X_{k+1}^N\|^2-\|X_k^N\|^2) \le C_F\tau\|X_{k+1}^N\|^2 +  \frac{1}{2}\|G(X_k^N)\Delta W_k\|^2 + (G(X_k^N)\Delta W_k,X_k^N).
	\end{equation*}
	By induction on $k\ge1$ in above inequality, we have
	\begin{equation*}
		\|X_{k+1}^N\|^2 \le \|X_{0}^N\|^2
		+ 2C_F\tau\sum_{i=0}^{k}\|X_{i+1}^N\|^2 +  \sum_{i=0}^{k}\|G(X_i^N)\Delta W_i\|^2 + \sum_{i=0}^{k}2(G(X_i^N)\Delta W_i,X_i^N).
	\end{equation*}
	By continuous version Burkholder-Davis-Gundy inequality \cite[Theorem 6.1.2]{LiuWei2015}, we obtain
	\begin{align*}
		&\ \mathbb{E}\|X_{k+1}^N\|^{2p} 
		\le 4^{p-1}\mathbb{E}\|X_{0}^N\|^{2p}
		+ 4^{p-1}\mathbb{E}\big(2C_F\tau\sum_{i=0}^{k}\|X_{i+1}^N\|^2\big)^{p} \\ 
		&\ +  4^{p-1}\mathbb{E}\big(\sum_{i=0}^{k}\|G(X_i^N)\Delta W_i\|^2\big)^{p} + 4^{p-1}\mathbb{E}\big(\sum_{i=0}^{k}2(G(X_i^N)\Delta W_i,X_i^N)\big)^{p}\\
		\le&\ C(p)\|\xi\|^{2p}
		+ C(p)\tau^p(k+1)^{p-1}\sum_{i=0}^{k}\mathbb{E}\|X_{i+1}^N\|^{2p} \\ 
		&\ +  C(p)(k+1)^{p-1}\sum_{i=0}^k\mathbb{E}\|\int_{t_i}^{t_{i+1}}G(X_i^N)dW(s)\|^{2p}\\ 
		&\ + C(p)\mathbb{E}\big(\sum_{i=0}^{k}\int_{t_i}^{t_{i+1}}(G(X_i^N)dW(s),X_i^N)\big)^{p}\\
		\le&\ C	+ C(T,p)\tau\sum_{i=0}^{k}\mathbb{E}\|X_{i+1}^N\|^{2p} + C(p)(k+1)^{p-1}\tau^{p-1}\sum_{i=0}^{k}\mathbb{E}\int_{t_{i}}^{t_{i+1}}\|G(X_i^N)\|^{2p}_{\mathcal{L}_2^0}ds\\ 
		&\ + C(p)\mathbb{E}\big(\sum_{i=0}^k\tau\|G(X_i^N)\|^2_{\mathcal{L}_2^0}\|X_i^N\|^2\big)^{\frac{p}{2}}\\
		\le&\ C	+ C\tau \sum_{i=0}^{k}\mathbb{E}\|X_{i+1}^N\|^{2p} + C(T,p)\tau \sum_{i=0}^{k}\mathbb{E}[C_G(1+\|X_i^N\|)]^{2p}\\ 
		&\ + C(k+1)^{\frac{p}{2}-1}\tau^{\frac{p}{2}}\sum_{i=0}^{k}C(1+\mathbb{E}\|X_i^N\|^{2p})\\
		\le&\  C + C\tau \sum_{i=0}^{k}\mathbb{E}\|X_{i+1}^N\|^{2p} + 2\Big(C\tau \sum_{i=0}^{k}1 + C\tau\sum_{i=0}^{k}\mathbb{E}\|X_i^N\|^{2p}\Big)\\
		\le&\ C + C\tau\sum_{i=0}^{k+1}\mathbb{E}\|X_{i}^N\|^{2p},
	\end{align*}
	where $C>0$ is a constant independent of both $\tau$ and $N$.
	From a discrete Gronwall's inequality \cite[Lemma A.3]{Kruse2014}, it follows that
	\begin{equation*}
		\mathbb{E}\|X_{k+1}^N\|^{2p} \le C e^{C(k+1)\tau}\le C e^{CT} \le C,\quad \forall 1\le k\le K.
	\end{equation*}
	
	Note that for any $0<\alpha<1$ and $0\le k\le K-1$,
	\begin{equation}\label{eq3.3.15}
		\tau\sum_{i=0}^{k}t_{k+1-i}^{-\alpha} \le \tau\sum_{i=0}^{K-1}t_{K-i}^{-\alpha} \le \sum_{i=0}^{K-1}\int_{t_i}^{t_{i+1}}s^{-\alpha}ds = \int_0^Ts^{-\alpha}ds = \frac{T^{1-\alpha}}{1-\alpha} := C(\alpha).
	\end{equation}
	
	Now we turn to estimate $X_k^N$ in $L^{2p}(\Omega;\dot{H}^{\gamma})$-norm. By Sobolev embedding $H\hookrightarrow W^{-1,6}(\mathcal{D})$, \eqref{eq3.3.4}, \eqref{eq3.2.5}, Burkholder-Davis-Gundy inequality in Lemma \ref{l3.2.00}, \eqref{eq3.3.15} and noticing that $\frac{\gamma+1}{2}<1$, we have
	\begin{align*}
		&\ \|X_{k+1}^N\|_{L^{2p}(\Omega;\dot{H}^{\gamma})} = \|A^{\frac{\gamma}{2}}X_{k+1}^N\|_{L^{2p}(\Omega;H)} \\ 
		\le&\ \|S_{N,\tau}^{k+1}P_NA^{\frac{\gamma}{2}}\xi\|_{L^{2p}(\Omega;H)} + \tau\sum_{i=0}^{k}\|A^{\frac{\gamma}{2}+\frac{1}{2}}S_{N,\tau}^{k+1-i}P_NA^{-\frac{1}{2}}F(X_{i+1}^N)\|_{L^{2p}(\Omega;H)}\\ 
		&\ + \Big(\mathbb{E}\Big\|\sum_{i=0}^{k}A^{\frac{\gamma}{2}}S_{N,\tau}^{k+1-i}P_NG(X_i^N)\Delta W_i\Big\|^{2p}\Big)^{\frac{1}{2p}}\\
		\le&\ \|S_{N,\tau}^{k+1}P_N\|_{\mathcal{L}(H)}\|A^{\frac{\gamma}{2}}\xi\| + \tau\sum_{i=0}^{k}\|A^{\frac{\gamma+1}{2}}S_{N,\tau}^{k+1-i}P_N\|_{\mathcal{L}(H)}\big(\mathbb{E}\|A^{-\frac{1}{2}}[X_{i+1}^N-(X_{i+1}^N)^3]\|^{2p}\big)^{\frac{1}{2p}}\\ 
		&\ + C(p)\Big(\sum_{i=0}^{k}\big(\mathbb{E}\|A^{\frac{\gamma}{2}}S_{N,\tau}^{k+1-i}P_NG(X_i^N)\Delta W_i\|^{2p}\big)^{\frac{1}{p}}\Big)^{\frac{1}{2}}\\
		\le&\ C\|\xi\|_{\gamma} + \tau\sum_{i=0}^{k}Ct_{k+1-i}^{-\frac{\gamma+1}{2}} \cdot \big(C\mathbb{E}\|X_{i+1}^N\|_{W^{-1,6}(\mathcal{D)}}^{6p}\big)^{\frac{1}{2p}}\\ 
		&\ + C(p)\Big(\sum_{i=0}^{k}Ct_{k+1-i}^{-\gamma}\big(\mathbb{E}\|G(X_i^N)\Delta W_i\|^{2p}\big)^{\frac{1}{p}}\Big)^{\frac{1}{2}}\\
		\le&\ C + \tau\sum_{i=0}^{k}Ct_{k+1-i}^{-\frac{\gamma+1}{2}} \cdot \big(\mathbb{E}(C\|X_{i+1}^N\|)^{6p}\big)^{\frac{1}{2p}} + C\Big(\sum_{i=0}^{k}t_{k+1-i}^{-\gamma}\tau\big(\mathbb{E}\|G(X_i^N)\|_{\mathcal{L}_2^0}^{2p}\big)^{\frac{1}{p}}\Big)^{\frac{1}{2}}\\
		\le&\ C + \tau\sum_{i=0}^{k}Ct_{k+1-i}^{-\frac{\gamma+1}{2}}\cdot C +  C\Big(C\sum_{i=0}^{k}t_{k+1-i}^{-\gamma}\tau\Big)^{\frac{1}{2}} \le C,
	\end{align*}
	where $C>0$ is a constant independent of both $\tau$ and $N$. This together with Sobolev embedding property $\dot{H}^{\gamma} \hookrightarrow V$ implies $\|X_{k+1}^N\|_{L^{2p}(\Omega;V)} \le C$.
\end{proof}

It follows directly from this lemma and \eqref{eq3.3.9} that $\widetilde{X}(t)$ satisfies
\begin{equation}\label{eq3.3.14'r}
	\mathbb{E}\|\widetilde{X}(t)\|_{\gamma}^{2p} + \mathbb{E}\|\widetilde{X}(t)\|_V^{2p} \le C,\quad \forall p\ge1,\ t\in[0,T].
\end{equation}

\section{Regularity of Kolmogorov equation}
\label{Sect4}
In this part, we investigate a priori estimates of a Kolmogorov equation associated with the semi-discretization equation \eqref{eq3.3.5}. 

Let $X^N(t,\nu)$ be the solution to \eqref{eq3.3.5} with initial value $\nu\in H_N$. 
According to \cite[(12)]{CuiJianbo2021}, for any given $\phi \in \mathcal{C}_b^3(H)$, $U^N(t,\nu) := \mathbb{E}[\phi(X^N(t,\nu))]$ satisfies a Kolmogorov equation with initial value $\phi(\nu)$
\begin{equation}\label{eq3.4.1}
	\begin{aligned}
		\frac{\partial U^N(t,\nu)}{\partial t}=&\ (DU^N(t,\nu),-A\nu+P_NF(\nu))\\
		&\ +\frac{1}{2}\mathrm{Tr}\{D^2U^N(t,\nu)(P_NG(\nu)Q^{\frac{1}{2}})(P_NG(\nu)Q^{\frac{1}{2}})^*\}.
	\end{aligned}
\end{equation}
The main task of this section is to derive a priori estimates for $D^iU^N(t,\nu)$ $(i=1,2,3)$.

We begin with estimating $DU^N(t,\nu)$. According to \cite[Theorem 9.8]{DaPrato2014}, differentiating $U^N(t,\nu)$ with respect to $\nu$ along direction $h \in H_N$, we obtain
\begin{equation}\label{eq3.4.13.1}
	DU^N(t,\nu)h = \mathbb{E}[D\phi(X^N(t,\nu))\eta^{h,\nu}(t)],
\end{equation}
where $\eta^{h,\nu}(t) = DX^N(t,\nu)h\in H_N$ satisfies $\eta^{h,\nu}(0) =h$ and an SPDE
\begin{equation}\label{eq3.4.13.6}
	d\eta^{h,\nu}(t) = -A\eta^{h,\nu}(t)dt+P_NDF(X^N(t))\eta^{h,\nu}(t)dt + P_NDG(X^N(t))\eta^{h,\nu}(t)dW(t).
\end{equation}
Then the regularity estimate for $DU^N(t,\nu)$ is converted to studying $\eta^{h,\nu}(t)$.

\begin{lemma}\label{l3.4.1}
	Assume (A1)-(A3). Then for any $p \ge 1$, $0<\beta < \frac{1}{2}$ and $\alpha<\frac{1}{2}$, there exists a constant $C > 0$ independent of $N$ such that, for all $\nu, h \in H_N$ and $t \in [0,T]$,
	\begin{align*}
		\|\eta^{h,\nu}(t)\|_{L^{2p}(\Omega;H)} &\le C \|h\|,\\
		\|A^{\beta}\eta^{h,\nu}(t)\|_{L^{2p}(\Omega;H)} &\le C \bigl(1 + t^{-\beta} + \|\nu\|_{\gamma}^6\bigr) \|h\|,\\
		\|\eta^{A^{\alpha}h,\nu}(t)\|_{L^{2p}(\Omega;H)} &\le C \bigl(1 + t^{-\max\{\alpha,0\}} + \|\nu\|_{\gamma}^6\bigr) \|h\|.
	\end{align*}
	Moreover, if $\beta+\alpha\ge0$, there exists a constant $C > 0$ independent of $N$ such that, for all $\nu, h \in H_N$ and $t \in [0,T]$,
	\begin{equation*}
		\|A^{\beta}\eta^{A^{\alpha}h,\nu}(t)\|_{L^{2p}(\Omega;H)} \le C \bigl(1 + t^{-(\beta+\alpha)} + \|\nu\|_{\gamma}^{12}\bigr) \|h\|.
	\end{equation*}
\end{lemma}

\begin{proof}
	Let $p \ge 1$ and define a functional $\rho(u) = \|u\|^{2p}$ on $H$.  
	For any $u,v \in H$, cf.\ \cite[(3.22)]{Gawarecki2011}, we have
	\begin{equation}\label{eq3.4.10}
		(D\rho(u),v) = 2p \|u\|^{2p-2} (u,v), \quad
		\|D^2 \rho(u)\|_{\mathcal{L}(H)} \le 2p (2p - 1) \|u\|^{2p-2}.
	\end{equation}
	Applying It\^o's formula \cite[Theorem 1]{DaPrato2019} to $\rho\bigl(\eta^{h,\nu}(t)\bigr)$ yields
	\begin{align*}
		&\ \rho(\eta^{h,\nu}(t))
		= \rho(h)
		+ \int_0^t \bigl(D \rho(\eta^{h,\nu}(s)),
		-A \eta^{h,\nu}(s) + P_N D F(X^N(s)) \eta^{h,\nu}(s)\bigr) ds \\
		& + \int_0^t \bigl(D \rho(\eta^{h,\nu}(s)),
		P_N D G(X^N(s)) \eta^{h,\nu}(s) dW(s)\bigr) \\
		& + \frac{1}{2}\int_0^t \operatorname{Tr}\!\bigl\{
		D^2 \rho(\eta^{h,\nu}(s))
		(P_N D G(X^N(s)) \eta^{h,\nu}(s) Q^{1/2})
		(P_N D G(X^N(s)) \eta^{h,\nu}(s) Q^{1/2})^*
		\bigr\} ds.
	\end{align*}
	Clearly,
	\[
	(D \rho(\eta^{h,\nu}(s)), -A \eta^{h,\nu}(s))
	= -2p \|\eta^{h,\nu}(s)\|^{2p-2}
	\|A^{1/2} \eta^{h,\nu}(s)\|^2 \le 0.
	\]
	This together with \eqref{eq3.4.10}, \eqref{eq3.2.5} and Lemma \ref{l3.2.1} implies
	\begin{align*}
		&\ \mathbb{E}\|\eta^{h,\nu}(t)\|^{2p} = \mathbb{E}\rho(\eta^{h,\nu}(t))\\ 
		\le&\ \mathbb{E}\rho(h) + \mathbb{E}\int_0^t2p\|\eta^{h,\nu}(s)\|^{2p-2}(\eta^{h,\nu}(s),P_NDF(X^N(s))\eta^{h,\nu}(s))ds\\
		&+\frac{1}{2}\mathbb{E}\int_0^t\|D^2\rho(\eta^{h,\nu}(s))\|_{\mathcal{L}(H)}\|P_NDG(X^N(s))\eta^{h,\nu}(s)\|^2_{\mathcal{L}_2^0}ds\\
		\le&\ \|h\|^{2p} + \mathbb{E}\int_0^t2p\|\eta^{h,\nu}(s)\|^{2p-2}\cdot C_F\|\eta^{h,\nu}(s)\|^2ds\\
		&+ \frac{1}{2}\mathbb{E}\int_0^t2p(2p-1)\|\eta^{h,\nu}(s)\|^{2p-2}\cdot C_G\|\eta^{h,\nu}(s)\|^2ds\\
		\le&\ \|h\|^{2p} + C\mathbb{E}\int_0^t\|\eta^{h,\nu}(s)\|^{2p}ds,
	\end{align*}
	where $C>0$ is a constant independent of $N$.
	By Gronwall's inequality, we get
	\begin{equation}\label{eq3.4.13.2}
		\mathbb{E}\|\eta^{h,\nu}(t)\|^{2p} \le C\|h\|^{2p},\quad\forall\ t\in[0,T].
	\end{equation}
	This implies the first assertion.
	
	For any $\beta\in(0,\frac{1}{2})$, applying $A^{\beta}$ to mild solution of \eqref{eq3.4.13.6}, we have
	\begin{equation}\label{eq3.4.13.3}
		\begin{aligned}
			A^{\beta}\eta^{h,\nu}(t) =&\ A^{\beta}S(t)h + \int_0^tA^{\beta}S(t-s)P_NDF(X^N(s))\eta^{h,\nu}(s)ds\\ 
			&\ + \int_0^tA^{\beta}S(t-s)P_NDG(X^N(s))\eta^{h,\nu}(s)dW(s).
		\end{aligned}
	\end{equation}
	By H\"older inequality and Burkholder-Davis-Gundy inequality, \eqref{eq3.3.3}, \eqref{eq3.2.5}, Lemma \ref{l3.2.1}, \eqref{eq3.3.2} and \eqref{eq3.4.13.2}, we obtain
	\begin{align*}
		&\ \mathbb{E}\|A^{\beta}\eta^{h,\nu}(t)\|^{2p}
		\le 3^{2p-1}\|A^{\beta}S(t)h\|^{2p}\\ 
		&\ + 3^{2p-1}\mathbb{E}\Big(t\int_0^t\|A^{\beta}S(t-s)P_NDF(X^N(s))\eta^{h,\nu}(s)\|^2ds\Big)^{p}\\
		&\ + 3^{2p-1}\mathbb{E}\Big(\int_0^t\|A^{\beta}S(t-s)P_NDG(X^N(s))\eta^{h,\nu}(s)\|_{\mathcal{L}_2^0}^2ds\Big)^{p}\\
		\le&\ 3^{2p-1}\|A^{\beta}S(t)\|_{\mathcal{L}(H)}^{2p}\|h\|^{2p} \!+\! 3^{2p-1}T^{p}\mathbb{E}\Big(\int_0^t\|A^{\beta}S(t-s)P_N\|_{\mathcal{L}(H)}^2\|DF(X^N(s))\eta^{h,\nu}(s)\|^2ds\Big)^{p}\\
		&\ + 3^{2p-1}\mathbb{E}\Big(\int_0^t\|A^{\beta}S(t-s)P_N\|_{\mathcal{L}(H)}^2\|DG(X^N(s))\eta^{h,\nu}(s)\|_{\mathcal{L}_2^0}^2ds\Big)^{p}\\
		\le&\ Ct^{-2p\beta}\|h\|^{2p} + C\mathbb{E}\Big(\int_0^tC(t-s)^{-2\beta}C_F^2(1+\|X^N(s)\|_V^2)^2\|\eta^{h,\nu}(s)\|^2ds\Big)^{p}\\
		&\ + C\mathbb{E}\Big(\int_0^tC(t-s)^{-2\beta}C_G^2\|\eta^{h,\nu}(s)\|^2ds\Big)^{p}\\
		\le&\ Ct^{-2p\beta}\|h\|^{2p} + C\mathbb{E}[\sup_{r\in[0,T]}C_F(1+\|X^N(r)\|_V^2)^{2p}\|\eta^{h,\nu}(r)\|^{2p}] + C\mathbb{E}[\sup_{r\in[0,T]}C_G\|\eta^{h,\nu}(r)\|^{2p}]\\ 
		\le&\ Ct^{-2p\beta}\|h\|^{2p} + C\big(\mathbb{E}\sup_{r\in[0,T]}(1+\|X^N(r)\|_V^2)^{4p}\big)^{\frac{1}{2}}\big(\mathbb{E}\sup_{r\in[0,T]}\|\eta^{h,\nu}(r)\|^{4p}\big)^{\frac{1}{2}} + C\|h\|^{2p}\\
		\le&\ C(t^{-2p\beta}+\|\nu\|_{\gamma}^{12p}+1)\|h\|^{2p}.
	\end{align*}
	This implies the second assertion.
	
	For any $\alpha<\frac{1}{2}$, define $\tilde{\eta}^{A^{\alpha}h,\nu}(t) = \eta^{A^{\alpha}h,\nu}(t)-S(t)A^{\alpha}h$, then $\tilde{\eta}^{A^{\alpha}h,\nu}(0) = 0$ and
	\begin{equation*}
		d\tilde{\eta}^{A^{\alpha}h,\nu}(t) = d\eta^{A^{\alpha}h,\nu}(t)+ AS(t)A^{\alpha}hdt.
	\end{equation*}
	According to \eqref{eq3.4.13.6} and above equation, we get an equation for $\tilde{\eta}^{A^{\alpha}h,\nu}$ as
	\begin{equation}\label{eq3.4.6}
		\begin{aligned}
			&\ d\tilde{\eta}^{A^{\alpha}h,\nu}(t)\\ =&\ -A\tilde{\eta}^{A^{\alpha}h,\nu}(t)dt+P_NDF(X^N(t))\tilde{\eta}^{A^{\alpha}h,\nu}(t)dt + P_NDG(X^N(t))\tilde{\eta}^{A^{\alpha}h,\nu}(t)dW(t)\\
			&\ + P_NDF(X^N(t))S(t)A^{\alpha}hdt + P_NDG(X^N(t))S(t)A^{\alpha}hdW(t).
		\end{aligned}
	\end{equation}
	Applying It\^o's formula to $\rho\big(\eta^{h,\nu}(t)\big) = \|\tilde{\eta}^{h,\nu}(t)\|^{2p}$, we have
	\begin{align*}
		& \|\tilde{\eta}^{A^{\alpha}h,\nu}(t)\|^{2p}\\ 
		=& \int_0^t2p\|\tilde{\eta}^{A^{\alpha}h,\nu}(s)\|^{2p-2}\big(\tilde{\eta}^{A^{\alpha}h,\nu}(s),-A\tilde{\eta}^{A^{\alpha}h,\nu}(s)+P_NDF(X^N(s))\tilde{\eta}^{A^{\alpha}h,\nu}(s)\big)ds\\
		+& \int_0^t2p\|\tilde{\eta}^{A^{\alpha}h,\nu}(s)\|^{2p-2}\big(\tilde{\eta}^{A^{\alpha}h,\nu}(s),P_NDF(X^N(s))S(s)A^{\alpha}h\big)ds\\ 
		+& \int_0^t2p\|\tilde{\eta}^{A^{\alpha}h,\nu}(s)\|^{2p-2}\big(\tilde{\eta}^{A^{\alpha}h,\nu}(s),P_NDG(X^N(s))\tilde{\eta}^{A^{\alpha}h,\nu}(s)dW(s)\big)\\
		+& \int_0^t2p\|\tilde{\eta}^{A^{\alpha}h,\nu}(s)\|^{2p-2}\big(\tilde{\eta}^{A^{\alpha}h,\nu}(s),P_NDG(X^N(s))S(s)A^{\alpha}hdW(s)\big)\\
		+& \frac{1}{2}\int_0^t\mathrm{Tr}\big\{D^2\rho(\tilde{\eta}^{A^{\alpha}h,\nu}(s))(P_NDG(X^N(s))\tilde{\eta}^{A^{\alpha}h,\nu}(s)Q^{\frac{1}{2}})(P_NDG(X^N(s))\tilde{\eta}^{A^{\alpha}h,\nu}(s)Q^{\frac{1}{2}})^*\big\}ds\\
		+& \frac{1}{2}\int_0^t\mathrm{Tr}\big\{D^2\rho(\tilde{\eta}^{A^{\alpha}h,\nu}(s))(P_NDG(X^N(s))\tilde{\eta}^{A^{\alpha}h,\nu}(s)Q^{\frac{1}{2}})(P_NDG(X^N(s))S(s)A^{\alpha}hQ^{\frac{1}{2}})^*\big\}ds\\
		+& \frac{1}{2}\int_0^t\mathrm{Tr}\big\{D^2\rho(\tilde{\eta}^{A^{\alpha}h,\nu}(s))(P_NDG(X^N(s))S(s)A^{\alpha}hQ^{\frac{1}{2}})(P_NDG(X^N(s))\tilde{\eta}^{A^{\alpha}h,\nu}(s)Q^{\frac{1}{2}})^*\big\}ds\\
		+& \frac{1}{2}\int_0^t\mathrm{Tr}\big\{D^2\rho(\tilde{\eta}^{A^{\alpha}h,\nu}(s))(P_NDG(X^N(s))S(s)A^{\alpha}hQ^{\frac{1}{2}})(P_NDG(X^N(s))S(s)A^{\alpha}hQ^{\frac{1}{2}})^*\big\}ds.
	\end{align*}
	As $0\le\alpha<\frac{1}{2}$, noticing that $(\tilde{\eta}^{A^{\alpha}h,\nu}(s),-A\tilde{\eta}^{A^{\alpha}h,\nu}(s)) = -(A^{\frac{1}{2}}\tilde{\eta}^{A^{\alpha}h,\nu}(s),A^{\frac{1}{2}}\tilde{\eta}^{A^{\alpha}h,\nu}(s)) = -\|\tilde{\eta}^{A^{\alpha}h,\nu}(s)\|_1^2\le 0$ and by \eqref{eq3.3.3}, \eqref{eq3.4.10}, \eqref{eq3.2.5}, Lemma \ref{l3.2.1} and H\"older inequality with $1= \frac{2p-1}{2p}+\frac{1}{2p}$, we obtain
	\begin{align*}
		& \mathbb{E}\|\tilde{\eta}^{A^{\alpha}h,\nu}(t)\|^{2p}
		\le 2p\mathbb{E}\int_0^t\|\tilde{\eta}^{A^{\alpha}h,\nu}(s)\|^{2p-2}(\tilde{\eta}^{A^{\alpha}h,\nu}(s),P_NDF(X^N(s))\tilde{\eta}^{A^{\alpha}h,\nu}(s))ds\\ 
		& + 2p\mathbb{E}\int_0^t\|\tilde{\eta}^{A^{\alpha}h,\nu}(s)\|^{2p-2}\big(\frac{1}{2}\|\tilde{\eta}^{A^{\alpha}h,\nu}(s)\|^2+\frac{1}{2}\|P_NDF(X^N(s))S(s)A^{\alpha}h\|^2\big)ds\\
		& +\frac{1}{2}\mathbb{E}\int_0^t2p(2p-1)\|\tilde{\eta}^{A^{\alpha}h,\nu}(s)\|^{2p-2}\|P_NDG(X^N(s))\tilde{\eta}^{A^{\alpha}h,\nu}(s)\|^2_{\mathcal{L}_2^0}ds\\ 
		& + \mathbb{E}\int_0^t2p(2p-1)\|\tilde{\eta}^{A^{\alpha}h,\nu}(s)\|^{2p-2}\\ 
		&\ \cdot \|P_NDG(X^N(s))\tilde{\eta}^{A^{\alpha}h,\nu}(s)\|_{\mathcal{L}_2^0}\|P_NDG(X^N(s))S(s)A^{\alpha}h\|_{\mathcal{L}_2^0}ds\\ 
		& + \frac{1}{2}\mathbb{E}\int_0^t2p(2p-1)\|\tilde{\eta}^{A^{\alpha}h,\nu}(s)\|^{2p-2}\|P_NDG(X^N(s))S(s)A^{\alpha}h\|^2_{\mathcal{L}_2^0}ds\\
		\le& C\mathbb{E}\int_0^t\|\tilde{\eta}^{A^{\alpha}h,\nu}(s)\|^{2p-2}\cdot C_F\|\tilde{\eta}^{A^{\alpha}h,\nu}(s)\|^2ds\\ 
		& + C\mathbb{E}\int_0^t\|\tilde{\eta}^{A^{\alpha}h,\nu}(s)\|^{2p}ds + C\mathbb{E}\int_0^t\|\tilde{\eta}^{A^{\alpha}h,\nu}(s)\|^{2p-2}C_F^2(1+\|X^N(s)\|_V^2)^2\|S(s)A^{\alpha}h\|^2ds\\
		& +C\mathbb{E}\int_0^tC_G^2\|\tilde{\eta}^{A^{\alpha}h,\nu}(s)\|^{2p}ds + C\mathbb{E}\int_0^tC_G\|\tilde{\eta}^{A^{\alpha}h,\nu}(s)\|^{2p-2}C_G\|\tilde{\eta}^{A^{\alpha}h,\nu}(s)\|C_G\|S(s)A^{\alpha}h\|ds\\
		& + C\mathbb{E}\int_0^t\|\tilde{\eta}^{A^{\alpha}h,\nu}(s)\|^{2p-2}C_G^2\|S(s)A^{\alpha}h\|^2ds\\
		\le& C\mathbb{E}\int_0^t\|\tilde{\eta}^{A^{\alpha}h,\nu}(s)\|^{2p}ds +  C\mathbb{E}\int_0^t\|\tilde{\eta}^{A^{\alpha}h,\nu}(s)\|^{2p}ds\\
		&+ C\mathbb{E}\int_0^t\|\tilde{\eta}^{A^{\alpha}h,\nu}(s)\|^{2p-2}(1+\|X^N(s)\|_V^2)^2\cdot Cs^{-2\alpha}\|h\|^2ds\\
		&+ 	C\mathbb{E}\int_0^t\|\tilde{\eta}^{A^{\alpha}h,\nu}(s)\|^{2p}ds + C\mathbb{E}\int_0^t\|\tilde{\eta}^{A^{\alpha}h,\nu}(s)\|^{2p-1}\cdot Cs^{-\alpha}\|h\|ds\\
		&+ C\mathbb{E}\int_0^t\|\tilde{\eta}^{A^{\alpha}h,\nu}(s)\|^{2p-2}\cdot Cs^{-2\alpha}\|h\|^2ds\\
		\le& C\mathbb{E}\int_0^t\|\tilde{\eta}^{A^{\alpha}h,\nu}(s)\|^{2p}ds + C\int_0^t\big(\mathbb{E}\|\tilde{\eta}^{A^{\alpha}h,\nu}(s)\|^{2p}\big)^{\frac{2p-2}{2p}}\big(\mathbb{E}(1+\|X^N(s)\|_V^2)^{2p}\big)^{\frac{1}{p}}s^{-2\alpha}\|h\|^2ds\\
		& + C\mathbb{E}\int_0^t\|\tilde{\eta}^{A^{\alpha}h,\nu}(s)\|^{2p-1}\|h\|s^{-\alpha}ds + C\mathbb{E}\int_0^t\|\tilde{\eta}^{A^{\alpha}h,\nu}(s)\|^{2p-2}\|h\|^2s^{-2\alpha}ds\\
		\le& C\mathbb{E}\int_0^t\|\tilde{\eta}^{A^{\alpha}h,\nu}(s)\|^{2p}ds + C\int_0^t\big(\mathbb{E}\|\tilde{\eta}^{A^{\alpha}h,\nu}(s)\|^{2p}\big)^{\frac{2p-2}{2p}}(1+\|\nu\|_{\gamma}^6)^2\|h\|^2s^{-2\alpha}ds\\
		& + C\mathbb{E}\int_0^t\|\tilde{\eta}^{A^{\alpha}h,\nu}(s)\|^{2p-1}\|h\|s^{-\alpha}ds + C\mathbb{E}\int_0^t\|\tilde{\eta}^{A^{\alpha}h,\nu}(s)\|^{2p-2}\|h\|^2s^{-2\alpha}ds.
	\end{align*}
	By Young's inequalities $abs^{-2\alpha} \le \frac{2p-2}{2p}a^{\frac{2p}{2p-2}}s^{-2\alpha} + \frac{1}{p}b^{p}s^{-2\alpha}$ and $abs^{-\alpha} \le \frac{2p-1}{2p}a^{\frac{2p}{2p-1}}s^{-\alpha} + \frac{1}{2p}b^{2p}s^{-\alpha}$ for $a,b\ge 0$, we obtain
	\begin{align*}
		&\ \mathbb{E}\|\tilde{\eta}^{A^{\alpha}h,\nu}(t)\|^{2p}\\
		\le&\ C\int_0^t\mathbb{E}\|\tilde{\eta}^{A^{\alpha}h,\nu}(s)\|^{2p}ds + C\int_0^t\mathbb{E}\|\tilde{\eta}^{A^{\alpha}h,\nu}(s)\|^{2p}s^{-2\alpha}ds + C\int_0^t(1+\|\nu\|_{\gamma}^6)^{2p}\|h\|^{2p}s^{-2\alpha}ds\\
		&\ + C\int_0^t\mathbb{E}\|\tilde{\eta}^{A^{\alpha}h,\nu}(s)\|^{2p}s^{-\alpha}ds + C\int_0^t\|h\|^{2p}s^{-\alpha}ds\\
		&\ + C\int_0^t\mathbb{E}\|\tilde{\eta}^{A^{\alpha}h,\nu}(s)\|^{2p}s^{-2\alpha}ds + C\int_0^t\|h\|^{2p}s^{-2\alpha}ds\\
		\le&\ C\int_0^t(1+s^{-\alpha}+s^{-2\alpha})\mathbb{E}\|\tilde{\eta}^{A^{\alpha}h,\nu}(s)\|^{2p}ds + C(1+\|\nu\|_{\gamma}^6)^{2p}\|h\|^{2p}.
	\end{align*}
	Applying Gronwall's inequality, we get
	\begin{equation*}
		\mathbb{E}\|\tilde{\eta}^{A^{\alpha}h,\nu}(t)\|^{2p} \le C(1+\|\nu\|_{\gamma}^6)^{2p}\|h\|^{2p}\cdot\exp\Big(C\int_0^t1+s^{-\alpha}+s^{-2\alpha}ds\Big) \le C(1+\|\nu\|_{\gamma}^6)^{2p}\|h\|^{2p}.
	\end{equation*}
	As $\alpha<0$, we have $\|A^{\alpha}\|_{\mathcal{L}(H)} \le C$. Similar argument as above yields
	\begin{align*}
		& \mathbb{E}\|\tilde{\eta}^{A^{\alpha}h,\nu}(t)\|^{2p}
		\le C\mathbb{E}\int_0^t\|\tilde{\eta}^{A^{\alpha}h,\nu}(s)\|^{2p-2}\cdot C_F\|\tilde{\eta}^{A^{\alpha}h,\nu}(s)\|^2ds\\ 
		& + C\mathbb{E}\int_0^t\|\tilde{\eta}^{A^{\alpha}h,\nu}(s)\|^{2p}ds + C\mathbb{E}\int_0^t\|\tilde{\eta}^{A^{\alpha}h,\nu}(s)\|^{2p-2}C_F^2(1+\|X^N(s)\|_V^2)^2\|S(s)A^{\alpha}h\|^2ds\\
		& +C\mathbb{E}\int_0^tC_G^2\|\tilde{\eta}^{A^{\alpha}h,\nu}(s)\|^{2p}ds + C\mathbb{E}\int_0^tC_G\|\tilde{\eta}^{A^{\alpha}h,\nu}(s)\|^{2p-2}C_G\|\tilde{\eta}^{A^{\alpha}h,\nu}(s)\|C_G\|S(s)A^{\alpha}h\|ds\\
		& + C\mathbb{E}\int_0^t\|\tilde{\eta}^{A^{\alpha}h,\nu}(s)\|^{2p-2}C_G^2\|S(s)A^{\alpha}h\|^2ds\\
		\le& C\mathbb{E}\int_0^t\|\tilde{\eta}^{A^{\alpha}h,\nu}(s)\|^{2p}ds + C\mathbb{E}\int_0^t\|\tilde{\eta}^{A^{\alpha}h,\nu}(s)\|^{2p-2}(1+\|X^N(s)\|_V^2)^2\cdot C\|h\|^2ds\\
		&+  C\mathbb{E}\int_0^t\|\tilde{\eta}^{A^{\alpha}h,\nu}(s)\|^{2p-1}\cdot C\|h\|ds + C\mathbb{E}\int_0^t\|\tilde{\eta}^{A^{\alpha}h,\nu}(s)\|^{2p-2}\cdot C\|h\|^2ds\\
		\le& C\mathbb{E}\int_0^t\|\tilde{\eta}^{A^{\alpha}h,\nu}(s)\|^{2p}ds + C\int_0^t\big(\mathbb{E}\|\tilde{\eta}^{A^{\alpha}h,\nu}(s)\|^{2p}\big)^{\frac{2p-2}{2p}}\big(\mathbb{E}(1+\|X^N(s)\|_V^2)^{2p}\big)^{\frac{1}{p}}\|h\|^2ds\\
		& + C\mathbb{E}\int_0^t\|\tilde{\eta}^{A^{\alpha}h,\nu}(s)\|^{2p-1}\|h\|ds + C\mathbb{E}\int_0^t\|\tilde{\eta}^{A^{\alpha}h,\nu}(s)\|^{2p-2}\|h\|^2ds\\
		\le& C\mathbb{E}\int_0^t\|\tilde{\eta}^{A^{\alpha}h,\nu}(s)\|^{2p}ds + C\int_0^t\big(\mathbb{E}\|\tilde{\eta}^{A^{\alpha}h,\nu}(s)\|^{2p}\big)^{\frac{2p-2}{2p}}(1+\|\nu\|_{\gamma}^6)^2\|h\|^2ds\\
		&\ + C\int_0^t\mathbb{E}\|\tilde{\eta}^{A^{\alpha}h,\nu}(s)\|^{2p}ds + C\int_0^t\|h\|^{2p}ds\\
		&\ + C\int_0^t\mathbb{E}\|\tilde{\eta}^{A^{\alpha}h,\nu}(s)\|^{2p}ds + C\int_0^t\|h\|^{2p}ds\\
		\le&\ C\mathbb{E}\int_0^t\|\tilde{\eta}^{A^{\alpha}h,\nu}(s)\|^{2p}ds + C\int_0^t\mathbb{E}\|\tilde{\eta}^{A^{\alpha}h,\nu}(s)\|^{2p}ds + C\int_0^t(1+\|\nu\|_{\gamma}^6)^{2p}\|h\|^{2p}ds\\
		&\ + C\mathbb{E}\int_0^t\|\tilde{\eta}^{A^{\alpha}h,\nu}(s)\|^{2p}ds + C\int_0^t\|h\|^{2p}ds\\
		\le&\ C\mathbb{E}\int_0^t\|\tilde{\eta}^{A^{\alpha}h,\nu}(s)\|^{2p}ds + C(1+\|\nu\|_{\gamma}^6)^{2p}\|h\|^{2p}.
	\end{align*}
	By Gronwall's inequality, we have
	\begin{equation*}
		\mathbb{E}\|\tilde{\eta}^{A^{\alpha}h,\nu}(t)\|^{2p} \le  C(1+\|\nu\|_{\gamma}^6)^{2p}\|h\|^{2p}.
	\end{equation*}
	Note that $\|S(t)A^{\alpha}\|_{\mathcal{L}(H)} \le Ct^{-\alpha}$ for $0\le\alpha<\frac{1}{2}$ and $\|S(t)A^{\alpha}\|_{\mathcal{L}(H)} \le C$ for $\alpha<0$.
	Thus for $\alpha<\frac{1}{2}$ we obtain
	\begin{align*}
		&\ \mathbb{E}\|\eta^{A^{\alpha}h,\nu}(t)\|^{2p} \le C(p)\mathbb{E}\|\tilde{\eta}^{A^{\alpha}h,\nu}(t)\|^{2p} + C(p)\mathbb{E}\|S(t)A^{\alpha}h\|^{2p}\\ 
		\le&\ C(1+\|\nu\|_{\gamma}^{6})^{2p}\|h\|^{2p} + \|S(t)A^{\alpha}\|_{\mathcal{L}(H)}^{2p}\|h\|^{2p} \le C(1+\|\nu\|_{\gamma}^{12p}+t^{-2p\max\{\alpha,0\}})\|h\|^{2p}.
	\end{align*}
	This implies the third assertion in lemma.
	
	Applying $A^{\beta}$ to the mild solution of \eqref{eq3.4.6}, we have for any $\alpha<\frac{1}{2}$ satisfying $\beta+\alpha\ge0$,
	\begin{equation}\label{eq3.4.6.1}
		\begin{aligned}
			A^{\beta}\eta^{A^{\alpha}h,\nu}(t) =&\ A^{\beta}S(t)A^{\alpha}h + \int_0^tA^{\beta}S(t-s)P_NDF(X^N(s))\eta^{A^{\alpha}h,\nu}(s)ds\\ 
			&\ + \int_0^tA^{\beta}S(t-s)P_NDG(X^N(s))\eta^{A^{\alpha}h,\nu}(s)dW(s).
		\end{aligned}
	\end{equation}
	By Burkholder-Davis-Gundy inequality, \eqref{eq3.3.3}, \eqref{eq3.2.5}, Lemma \ref{l3.2.1}, \eqref{eq3.3.2} and \eqref{eq3.4.13.2}, we obtain
	\begin{align*}
		&\ \|A^{\beta}\eta^{A^{\alpha}h,\nu}(t)\|_{L^{2p}(\Omega;H)}\\
		\le&\ \|A^{\beta}S(t)A^{\alpha}h\|_{L^{2p}(\Omega;H)} + \int_0^t\|A^{\beta}S(t-s)P_NDF(X^N(s))\eta^{A^{\alpha}h,\nu}(s)\|_{L^{2p}(\Omega;H)}ds\\
		&\ + \Big[\int_0^t\big(\mathbb{E}\|A^{\beta}S(t-s)P_NDG(X^N(s))\eta^{A^{\alpha}h,\nu}(s)\|_{\mathcal{L}_2^0}^{2p}\big)^{\frac{1}{p}}ds\Big]^{\frac{1}{2}}\\
		\le&\ Ct^{-(\beta+\alpha)}\|h\| + \int_0^tC(t-s)^{-\beta}\cdot [\mathbb{E}\big(C_F(1+\|X^N(s)\|_V^2)^{2p}\|\eta^{A^{\alpha}h,\nu}(s)\|^{2p}\big)]^{\frac{1}{2p}}ds\\
		&\ + \Big[\int_0^tC(t-s)^{-2\beta}\cdot C_G\big(\mathbb{E}\|\eta^{A^{\alpha}h,\nu}(s)\|^{2p}\big)^{\frac{1}{p}}ds\Big]^{\frac{1}{2}}\\ 
		\le&\ Ct^{-(\beta+\alpha)}\|h\| + C\int_0^t(t-s)^{-\beta}\cdot(1+\|\nu\|_{\gamma}^6)\cdot(1+s^{-\max\{\alpha,0\}}+\|\nu\|_{\gamma}^6)\|h\|ds\\ 
		&\ + C\Big[\int_0^t(t-s)^{-2\beta}\cdot(1+s^{-2\max\{\alpha,0\}}+\|\nu\|_{\gamma}^{12})\|h\|^2ds\Big]^{\frac{1}{2}}\\
		\le&\ C(1+ t^{-(\beta+\alpha)} + \|\nu\|_{\gamma}^{12})\|h\|.
	\end{align*}
	This implies the last assertion in lemma.
\end{proof}

Now we are ready to estimate the regularity of $DU^N(t,\nu)$.

\begin{theorem}\label{t3.4.3}
	Assume (A1)-(A3). Then for any given $\phi \in \mathcal{C}_b^2(H)$ and $0<\alpha<\frac{1}{2}$, there exists a constant $C>0$ independent of $N$, such that for all $\nu\in H_N$ and $t\in[0,T]$
	\begin{equation*}
		\begin{aligned}
			\|DU^N(t,\nu)\| \le&\ C,\\
			\|A^{\alpha}DU^N(t,\nu)\| \le&\ C(1+t^{-\alpha}+\|\nu\|_{\gamma}^6).
		\end{aligned}
	\end{equation*}
\end{theorem}

\begin{proof}
	For any $h\in H$, by \eqref{eq3.4.13.1}, H\"older inequality and Lemma \ref{l3.4.1} with $p=1$, we obtain
	\begin{equation*}
		|DU^N(t,\nu)h| = |\mathbb{E}[D\phi(X^N(t,\nu))\eta^{h,\nu}(t)]|
		\le |\phi|_{\mathcal{C}_b^1(H)}(\mathbb{E}\|\eta^{h,\nu}(t)\|^2)^{\frac{1}{2}}
		\le C\|h\|.
	\end{equation*}
	Similarly, by \eqref{eq3.4.13.1}, H\"older inequality and Lemma \ref{l3.4.1} with $p=1$, we obtain
	\begin{equation*}
		\begin{aligned}
			&\ |(A^{\alpha}DU^N(t,\nu),h)| = |\mathbb{E}(D\phi(X^N(t,\nu))DX^N(t,\nu),A^{\alpha}h)|\\
			=&\ |\mathbb{E}[D\phi(X^N(t,\nu))\eta^{A^{\alpha}h,\nu}(t)]|\le |\phi|_{\mathcal{C}_b^1(H)}(\mathbb{E}\|\eta^{A^{\alpha}h,\nu}(t)\|^2)^{\frac{1}{2}}\\
			\le&\ C(1+t^{-\alpha}+\|\nu\|_{\gamma}^6)\|h\|.
		\end{aligned}
	\end{equation*}
	Thus the proof is completed.
\end{proof}

Next, we estimate $D^2U^N(t,\nu)$. According to \cite[Theorem 9.9]{DaPrato2014}, differentiating $DU^N(t,\nu)h$ with respect to $\nu$ along direction $g \in H_N$ yields
\begin{equation}\label{eq3.4.13.4}
	D^2U^N(t,\nu)(h,g) = \mathbb{E}[D^2\phi(X^N(t,\nu))(\eta^{h,\nu}(t),\eta^{g,\nu}(t))]+ \mathbb{E}[D\phi(X^N(t,\nu))\zeta^{h,g,\nu}(t)],
\end{equation}
where $\zeta^{h,g,\nu}(t) = D^2X^N(t,\nu)(h,g) = D\eta^{h,\nu}(t)g\in H_N$ satisfies
\begin{equation}\label{eq3.4.13.7}
	\begin{aligned}
		&\ d\zeta^{h,g,\nu}(t) = -A\zeta^{h,g,\nu}(t)dt\\ &\ + \Big(P_ND^2F(X^N(t))(\eta^{h,\nu}(t),\eta^{g,\nu}(t)) + P_NDF(X^N(t))\zeta^{h,g,\nu}(t)\Big) dt\\
		&\ + \Big(P_ND^2G(X^N(t))(\eta^{h,\nu}(t),\eta^{g,\nu}(t)) + P_NDG(X^N(t))\zeta^{h,g,\nu}(t)\Big)dW(t)\\
		&\ := \big(a_{11}(t)+a_{12}(t)+a_{13}(t)\big)dt + \big(a_{21}(t)+a_{22}(t)\big)dW(t)
	\end{aligned}
\end{equation}
with $\zeta^{h,g,\nu}(0) = 0$.
By Lemma \ref{l3.2.1} and Sobolev embedding $\dot{H}^{\frac{1}{4}}\hookrightarrow L^4(\mathcal{D})$, we have
\begin{equation}\label{eq3.4.new}
	\begin{aligned}
		\|a_{21}(t)\|_{\mathcal{L}_2^0}^2 \le&\ C\|\eta^{h,\nu}(t)\|_{L^4(\mathcal{D})}^2\|\eta^{g,\nu}(t)\|_{L^4(\mathcal{D})}^2 \le C\|A^{\frac{1}{8}}\eta^{h,\nu}(t)\|^2\|A^{\frac{1}{8}}\eta^{g,\nu}(t)\|^2,\\
		\|a_{22}(t)\|_{\mathcal{L}_2^0}^2 \le&\ C\|\zeta^{h,g,\nu}(t)\|^2.
	\end{aligned}
\end{equation}
We emphasize that $\delta\in(0,\frac{\gamma}{2}-\frac{1}{4})$ represent a small parameter throughout the following analysis.

\begin{lemma}\label{l3.4.3}
	Assume (A1)-(A3). Then for any $0<\alpha_1<\frac{1}{2}$ and $\alpha_2$ satisfying $-\frac{1}{4}-\delta\le\alpha_2< \min\{\frac{3}{4}-2\alpha_1-\delta,0\}$, there exists a constant $C>0$ independent of $N$, such that for all $\nu,h,g\in H_N$ and $t \in [0,T]$,
	\begin{equation*}
		\begin{aligned}
			\|\zeta^{h,g,\nu}(t)\|_{L^{2}(\Omega;H)} \le&\ C(1+\|\nu\|^{12}_{\gamma})\|h\|\|g\|,\\
			\|\zeta^{A^{\alpha_1}h,A^{\alpha_2}g,\nu}(t)\|_{L^{2}(\Omega;H)} \le&\ C(1+\|\nu\|_{\gamma}^{18})\|h\|\|g\|,
		\end{aligned}  
	\end{equation*}
	where $\delta\in(0,\frac{\gamma}{2}-\frac{1}{4})$ is a small parameter.
\end{lemma}

\begin{proof}
	Applying It\^o's formula to $\rho(\zeta^{h,g,\nu}(t)) = \|\zeta^{h,g,\nu}(t)\|^{2}$, we obtain 
	\begin{equation}\label{eq3.4.12}
		\begin{aligned}
			&\ \rho(\zeta^{h,g,\nu}(t)) = \int_0^t (D\rho(\zeta^{h,g,\nu}(s)),a_{11}(s)+a_{12}(s)+a_{13}(s)) ds\\
			&\ + \int_0^t (D\rho(\zeta^{h,g,\nu}(s)),(a_{21}(s)+a_{22}(s))dW(s))\\
			&\ + \underset{\psi(t)}{\underbrace{\frac{1}{2}\int_0^t\mathrm{Tr}\{D^2\rho(\zeta^{h,g,\nu}(s))(a_{21}(s)+a_{22}(s))Q^{\frac{1}{2}}((a_{21}(s)+a_{22}(s))Q^{\frac{1}{2}})^*\}ds}}.
		\end{aligned}
	\end{equation}
	By \eqref{eq3.4.10}, we have
	\begin{equation}\label{eq3.4.13}
		\mathbb{E}(D\rho(\zeta^{h,g,\nu}(s)),a_{11}(s)) = -2\mathbb{E}(A^{\frac{1}{2}}\zeta^{h,g,\nu}(s),A^{\frac{1}{2}}\zeta^{h,g,\nu}(s)) = -2\mathbb{E}\|\zeta^{h,g,\nu}(s)\|_1^2 \le 0.
	\end{equation}
	According to \eqref{eq3.2.5}, H\"older inequality with $1 = \frac{1}{2} + \frac{1}{4} + \frac{1}{4}$ and Lemma \ref{l3.4.1}, we obtain
	\begin{equation}\label{eq3.4.14}
		\begin{aligned}
			&\ \mathbb{E}(D\rho(\zeta^{h,g,\nu}(s)),a_{12}(s)) = 2\mathbb{E}(\zeta^{h,g,\nu}(s),D^2F(X^N(t))(\eta^{h,\nu}(t),\eta^{g,\nu}(t)))\\ 
			=&\ 2\mathbb{E}[\big(A^{\frac{1}{2}}\zeta^{h,g,\nu}(s),A^{-\frac{1}{2}}D^2F(X^N(t))(\eta^{h,\nu}(t),\eta^{g,\nu}(t))\big)]\\
			\le&\ \mathbb{E}\|\zeta^{h,g,\nu}(s)\|^2_1 + \mathbb{E}\|D^2F(X^N(t))(\eta^{h,\nu}(t),\eta^{g,\nu}(t))\|_{-1}^2\\
			\le&\ \mathbb{E}\|\zeta^{h,g,\nu}(s)\|^2_1 + C_F^2\mathbb{E}\|X^N(t)\|^2_V\|\eta^{h,\nu}(t)\|^2\|\eta^{g,\nu}(t)\|^2\\
			\le&\ \mathbb{E}\|\zeta^{h,g,\nu}(s)\|^2_1 + (\mathbb{E}\|X^N(t)\|^{4}_V)^{\frac{1}{2}}(\mathbb{E}\|\eta^{h,\nu}(t)\|^{8})^{\frac{1}{4}}(\mathbb{E}\|\eta^{g,\nu}(t)\|^{8})^{\frac{1}{4}}\\
			\le&\ \mathbb{E}\|\zeta^{h,g,\nu}(s)\|^2_1 + C(1+\|\nu\|_{\gamma}^3)^2\|h\|^2\|g\|^2.
		\end{aligned}
	\end{equation}
	By \eqref{eq3.2.5}, we get
	\begin{equation}\label{eq3.4.15}
		\begin{aligned}
			&\ \mathbb{E}(D\rho(\zeta^{h,g,\nu}(s)),a_{13}(s)) = 2\mathbb{E}(\zeta^{h,g,\nu}(s),P_NDF(X^N(s))\zeta^{h,g,\nu}(s))\\
			=&\ 2\mathbb{E}(\zeta^{h,g,\nu}(s),DF(X^N(s))\zeta^{h,g,\nu}(s)) \le C_F\mathbb{E}\|\zeta^{h,g,\nu}(s)\|^{2}.
		\end{aligned}
	\end{equation}
	Applying \eqref{eq3.4.10}, \eqref{eq3.4.new} and Lemma \ref{l3.4.1}, we get
	\begin{equation}\label{eq3.4.16}
		\begin{aligned}
			&\ \mathbb{E}\psi(t)\\ \le&\ \frac{1}{2}\mathbb{E}\int_0^t\|D^2\rho(\zeta^{h,g,\nu}(s))\|_{\mathcal{L}(H)}\|a_{21}(s)+a_{22}(s)\|_{\mathcal{L}_2^0}^2ds\\
			\le&\ 2\mathbb{E}\int_0^t \|a_{21}(s)\|_{\mathcal{L}_2^0}^2+\|a_{22}(s)\|_{\mathcal{L}_2^0}^2ds\\
			\le&\ 2\mathbb{E}\int_0^tC\|A^{\frac{1}{8}}\eta^{h,\nu}(s)\|^2\|A^{\frac{1}{8}}\eta^{g,\nu}(s)\|^2 + C\|\zeta^{h,g,\nu}(s)\|^{2}ds\\
			\le&\ C\int_0^tC(1+s^{-\frac{1}{8}}+\|\nu\|_{\gamma}^6)^2\|h\|^{2} C(1+s^{-\frac{1}{8}}+\|\nu\|_{\gamma}^6)^2\|g\|^{2}ds + C\int_0^t\mathbb{E}\|\zeta^{h,g,\nu}(s)\|^{2}ds\\
			\le&\ C(1+\|\nu\|_{\gamma}^{12})^2\|h\|^2\|g\|^2 + C\int_0^t\mathbb{E}\|\zeta^{h,g,\nu}(s)\|^{2}ds.
		\end{aligned}
	\end{equation}
	Here we have used the fact that $\int_0^t(s^{-\frac{1}{8}})^2(s^{-\frac{1}{8}})^2ds=\int_0^ts^{-\frac{1}{2}}ds<\infty$.
	Substituting all estimates in \eqref{eq3.4.13}, \eqref{eq3.4.14}, \eqref{eq3.4.15} and \eqref{eq3.4.16} into the right hand side of \eqref{eq3.4.12} yields
	\begin{equation*}
		\mathbb{E}\|\zeta^{h,g,\nu}(t)\|^{2}
		\le C(1+\|\nu\|_{\gamma}^{12})^{2}\|h\|^{2}\|g\|^{2} + \int_0^tC\mathbb{E}\|\zeta^{h,g,\nu}(s)\|^{2}ds,
	\end{equation*}
	By Gronwall's inequality, we get
	\begin{equation*}
		\mathbb{E}\|\zeta^{h,g,\nu}(t)\|^{2}\le C(1+\|\nu\|_{\gamma}^{12})^{2}\|h\|^{2}\|g\|^{2}.
	\end{equation*}
	This implies the first assertion.
	
	For any $0<\alpha_1<\frac{1}{2}$ and $-\frac{1}{4}-\delta\le\alpha_2< \min\{\frac{3}{4}-2\alpha_1-\delta,0\}$, replacing  $h$ and $g$ by $A^{\alpha_1}h$ and $A^{\alpha_2}g$ in \eqref{eq3.4.13.7}, we have
	\begin{align*}
		&\ d\zeta^{A^{\alpha_1}h,A^{\alpha_2}g,\nu}(t) = -A\zeta^{A^{\alpha_1}h,A^{\alpha_2}g,\nu}(t)dt\\ &\ + \Big(P_ND^2F(X^N(t))(\eta^{A^{\alpha_1}h,\nu}(t),\eta^{A^{\alpha_2}g,\nu}(t)) + P_NDF(X^N(t))\zeta^{A^{\alpha_1}h,A^{\alpha_2}g,\nu}(t)\Big) dt\\
		&\ + \Big(P_ND^2G(X^N(t))(\eta^{A^{\alpha_1}h,\nu}(t),\eta^{A^{\alpha_2}g,\nu}(t)) + P_NDG(X^N(t))\zeta^{A^{\alpha_1}h,A^{\alpha_2}g,\nu}(t)\Big)dW(t)\\
		&\ := b_{11}(t)dt+\big(b_{12}(t)+b_{13}(t)\big)dt + \big(b_{21}(t)+b_{22}(t)\big)dW(t)
	\end{align*}
	with  $\zeta^{A^{\alpha_1}h,A^{\alpha_2}g,\nu}(0) = 0$.
	Applying It\^o's formula to $\rho(\zeta^{A^{\alpha_1}h,A^{\alpha_2}g,\nu}(t)) = \|\zeta^{A^{\alpha_1}h,A^{\alpha_2}g,\nu}(t)\|^2$, we get
	\begin{equation}\label{eq3.4.7}
		\begin{aligned}
			&\ \|\zeta^{A^{\alpha_1}h,A^{\alpha_2}g,\nu}(t)\|^2 = 2\int_0^t(\zeta^{A^{\alpha_1}h,A^{\alpha_2}g,\nu}(s),b_{11}(s)+b_{12}(s)+b_{13}(s))ds\\
			&\ + 2\int_0^t(\zeta^{A^{\alpha_1}h,A^{\alpha_2}g,\nu}(s),(b_{21}(s)+b_{22}(s))dW(s))\\
			&\ + \frac{1}{2}\int_0^t\mathrm{Tr}\{D^2\rho(\zeta^{A^{\alpha_1}h,A^{\alpha_2}g,\nu}(s))(b_{21}(s)+b_{22}(s))Q^{\frac{1}{2}}((b_{21}(s)+b_{22}(s))Q^{\frac{1}{2}})^*\}ds.
		\end{aligned}
	\end{equation}
	Note that
	\begin{equation*}
		\mathbb{E}(\zeta^{A^{\alpha_1}h,A^{\alpha_2}g,\nu}(s),b_{11}(s))
		= -\mathbb{E}\|\zeta^{A^{\alpha_1}h,A^{\alpha_2}g,\nu}(s)\|_1^2 \le 0.
	\end{equation*}
	According to \eqref{eq3.2.5}, H\"older inequality with $1=\frac{1}{2}+\frac{1}{4}+\frac{1}{4}$, \eqref{eq3.3.2} and Lemma \ref{l3.4.1}, we have
	\begin{align*}
		&\ \mathbb{E}(\zeta^{A^{\alpha_1}h,A^{\alpha_2}g,\nu}(s),b_{12}(s))\\ 
		=&\ \mathbb{E}\big(A^{\frac{1}{2}}\zeta^{A^{\alpha_1}h,A^{\alpha_2}g,\nu}(s),A^{-\frac{1}{2}}D^2F(X^N(s))(\eta^{A^{\alpha_1}h,\nu}(s),\eta^{A^{\alpha_2}g,\nu}(s))\big)\\ 
		\le&\ \frac{1}{2}\mathbb{E}\big[\|\zeta^{A^{\alpha_1}h,A^{\alpha_2}g,\nu}(s)\|_1^2+\|D^2F(X^N(s))(\eta^{A^{\alpha_1}h,\nu}(s),\eta^{A^{\alpha_2}g,\nu}(s))\|_{-1}^2\big]\\
		\le&\ \frac{1}{2}\mathbb{E}\|\zeta^{A^{\alpha_1}h,A^{\alpha_2}g,\nu}(s)\|_1^2 + \frac{1}{2}C_F^2\mathbb{E}\big[\|X^N(s)\|^2_V\|\eta^{A^{\alpha_1}h,\nu}(s)\|^2\|\eta^{A^{\alpha_2}g,\nu}(s)\|^2\big]\\
		\le&\ \frac{1}{2}\mathbb{E}\|\zeta^{A^{\alpha_1}h,A^{\alpha_2}g,\nu}(s)\|_1^2 + C(\mathbb{E}\|X^N(s)\|^4_V)^{\frac{1}{2}}(\mathbb{E}\|\eta^{A^{\alpha_1}h,\nu}(s)\|^8)^{\frac{1}{4}}(\mathbb{E}\|\eta^{A^{\alpha_2}g,\nu}(s)\|^8)^{\frac{1}{4}}\\
		\le&\ \frac{1}{2}\mathbb{E}\|\zeta^{A^{\alpha_1}h,A^{\alpha_2}g,\nu}(s)\|_1^2 + C(1+\|\nu\|_{\gamma}^3)^2(1+s^{-\alpha_1}+\|\nu\|^6_{\gamma})^2\|h\|^2(1+s^{-\max\{\alpha_2,0\}}+\|\nu\|^6_{\gamma})^2\|g\|^2\\
		\le&\ \frac{1}{2}\mathbb{E}\|\zeta^{A^{\alpha_1}h,A^{\alpha_2}g,\nu}(s)\|_1^2 + C(1+\|\nu\|_{\gamma}^9)^2(1+s^{-\alpha_1}+\|\nu\|_{\gamma}^{6})^2\|h\|^2\|g\|^2.
	\end{align*}
	From \eqref{eq3.2.5}, it follows that
	\begin{align*}
		&\ \mathbb{E}(\zeta^{A^{\alpha_1}h,A^{\alpha_2}g,\nu}(s),b_{13}(s))\\
		=&\ 2\mathbb{E}(\zeta^{A^{\alpha_1}h,A^{\alpha_2}g,\nu}(s),P_NDF(X^N(s))\zeta^{A^{\alpha_1}h,A^{\alpha_2}g,\nu}(s))\le C\mathbb{E}\|\zeta^{A^{\alpha_1}h,A^{\alpha_2}g,\nu}(s)\|^{2}.
	\end{align*}
	Noticing the Sobolev embedding $\dot{H}^{\frac{1}{2}+2\delta}\hookrightarrow V$ and by Lemma \ref{l3.2.1}, we have
	\begin{align*}
		&\|b_{21}(t)\|^2_{\mathcal{L}_2^0} \le C\|\eta^{A^{\alpha_1}h,\nu}(t)\|^2\|\eta^{A^{\alpha_2}g,\nu}(t)\|_{V}^2 \le C\|\eta^{A^{\alpha_1}h,\nu}(t)\|^2\|A^{\frac{1}{4}+\delta}\eta^{A^{\alpha_2}g,\nu}(t)\|^2,\\
		&\|b_{22}(t)\|^2_{\mathcal{L}_2^0} \le C\|\zeta^{A^{\alpha_1}h,A^{\alpha_2}g,\nu}(t)\|^2.
	\end{align*}
	According to the assumption on $\alpha_1$ and $\alpha_2$, we have $\frac{1}{4}+\delta+\alpha_2\ge 0$, $-2\alpha_1-\frac{1}{4}-\delta-\alpha_2>-1$.
	Thus by H\"older inequality and Lemma \ref{l3.4.1}, we obtain
	\begin{align*}
		&\ \mathbb{E}\mathrm{Tr}\{D^2\rho(\zeta^{A^{\alpha_1}h,A^{\alpha_2}g,\nu}(s))(b_{21}(s)+b_{22}(s))Q^{\frac{1}{2}}((b_{21}(s)+b_{22}(s))Q^{\frac{1}{2}})^*\}\\ 
		\le&\ C\mathbb{E}(\|b_{21}(s)\|_{\mathcal{L}_2^0}^2+\|b_{22}(s)\|_{\mathcal{L}_2^0}^2)ds\\ 
		\le&\  C\big(\mathbb{E}\|\eta^{A^{\alpha_1}h,\nu}(s)\|^4\big)^{\frac{1}{2}}\big(\mathbb{E}\|A^{\frac{1}{4}+\delta}\eta^{A^{\alpha_2}g,\nu}(s)\|^4\big)^{\frac{1}{2}} + C\mathbb{E}\|\zeta^{A^{\alpha_1}h,A^{\alpha_2}g,\nu}(s)\|^2\\
		\le&\ C(1+s^{-\alpha_1}+\|\nu\|_{\gamma}^6)^2\|h\|^2(1+s^{-\frac{1}{4}-\delta-\alpha_2}+\|\nu\|_{\gamma}^{12})^2\|g\|^2 + C\mathbb{E}\|\zeta^{A^{\alpha_1}h,A^{\alpha_2}g,\nu}(t)\|^2\\
		\le&\  \big((1+\|\nu\|_{\gamma}^{36})+(1+\|\nu\|_{\gamma}^{24})s^{-\alpha_1}+(1+\|\nu\|_{\gamma}^{12})s^{-\frac{1}{4}-\delta-\alpha_2}+s^{-2\alpha_1-\frac{1}{4}-\delta-\alpha_2}\big)\|h\|^2\|g\|^2\\ 
		&\ + C\mathbb{E}\|\zeta^{A^{\alpha_1}h,A^{\alpha_2}g,\nu}(t)\|^2.
	\end{align*}
	Finally, substituting all above estimates into the right hand side of \eqref{eq3.4.7} and noticing that $0\le2\max\{\alpha_1+\alpha_2,0\}<\frac{3}{4}<1$, we obtain
	\begin{equation}\label{eq3.4.8}
		\mathbb{E}\|\zeta^{A^{\alpha_1}h,A^{\alpha_2}g,\nu}(t)\|^2\\
		\le C(1+\|\nu\|^{36}_{\gamma})\|h\|^2\|g\|^2 + \int_0^tC\mathbb{E}\|\zeta^{A^{\alpha_1}h,A^{\alpha_2}g,\nu}(s)\|^{2}ds,
	\end{equation}
	By Gronwall's inequality, we get
	\begin{equation}\label{eq3.4.9}
		\mathbb{E}\|\zeta^{A^{\alpha_1}h,A^{\alpha_2}g,\nu}(t)\|^{2} 
		\le C(1+\|\nu\|^{36}_{\gamma})\|h\|^2\|g\|^2.
	\end{equation}
	Thus, the proof is completed.
\end{proof}

Now we are ready to estimate the regularity of $D^2U^N(t,\nu)$.

\begin{theorem}\label{t3.4.4}
	Assume (A1)-(A3). Then for any given $\phi \in \mathcal{C}_b^2(H)$, $0<\alpha_1<\frac{1}{2}$ and $-\frac{1}{4}-\delta\le\alpha_2< \min\{\frac{3}{4}-2\alpha_1-\delta,0\}$, there exists a constant $C>0$ independent of $N$, such that for all $\nu\in H_N$ and $t\in[0,T]$
	\begin{equation*}
		\begin{aligned}
			\|D^2U^N(t,\nu)\|_{\mathcal{L}(H)} \le&\ C(1+\|\nu\|_{\gamma}^{12}),\\
			\|A^{\alpha_1}D^2U^N(t,\nu)A^{\alpha_2}\|_{\mathcal{L}(H)} \le&\ C(1+t^{-\alpha_1}+\|\nu\|_{\gamma}^{18}),
		\end{aligned}
	\end{equation*}
	where $\delta\in(0,\frac{\gamma}{2}-\frac{1}{4})$ is a small parameter.
\end{theorem}

\begin{proof}
	For any $h,g\in H$, by Lemmas \ref{l3.4.1}, \ref{l3.4.3} and H\"older inequality, we have
	\begin{align*}
		&\ |D^2U^N(t,\nu)(h,g)| \\
		\le&\ \mathbb{E}|D^2\phi(X^N(t,\nu))\big(\eta^{h,\nu}(t),\eta^{g,\nu}(t)\big)| + \mathbb{E}|D\phi(X^N(t,\nu))\zeta^{h,g,\nu}(t)|\\
		\le&\ |\phi|_{\mathcal{C}_b^2(H)}\mathbb{E}\big[\|\eta^{h,\nu}(t)\|\cdot\|\eta^{g,\nu}(t)\|\big] + |\phi|_{\mathcal{C}_b^1(H)}C(\mathbb{E}\|\zeta^{h,g,\nu}(t)\|^2)^{\frac{1}{2}}\\
		\le&\ C(\mathbb{E}\|\eta^{h,\nu}(t)\|^2)^{\frac{1}{2}}(\mathbb{E}\|\eta^{g,\nu}(t)\|^2)^{\frac{1}{2}} + C(\mathbb{E}\|\zeta^{h,g,\nu}(t)\|^2)^{\frac{1}{2}}\\
		\le&\ C(1+\|\nu\|_{\gamma}^{12})\|h\|\|g\|.
	\end{align*}
	By Lemmas \ref{l3.4.1}, \ref{l3.4.3} and H\"older inequality, we obtain
	\begin{align*}
		&\ |(A^{\alpha_1}D^2U^N(t,\nu)A^{\alpha_2}h,g)| = |(D^2U^N(t,\nu)A^{\alpha_2}h,A^{\alpha_1}g)|\\
		\le&\ |\mathbb{E}D^2\phi(X^N(t,\nu))\big(\eta^{A^{\alpha_1}g,\nu}(t),\eta^{A^{\alpha_2}h,\nu}(t)\big)| + |\mathbb{E}[D\phi(X^N(t,\nu))\zeta^{A^{\alpha_1}g,A^{\alpha_2}h,\nu}(t)]|\\
		\le&\ |\phi|_{\mathcal{C}_b^2(H)}\mathbb{E}\big[\|\eta^{A^{\alpha_1}g,\nu}(t)\|\cdot\|\eta^{A^{\alpha_2}h,\nu}(t)\|\big] + |\phi|_{\mathcal{C}_b^1(H)}(\mathbb{E}\|\zeta^{A^{\alpha_1}g,A^{\alpha_2}h,\nu}(t)\|^2)^{\frac{1}{2}}\\
		\le&\ C(\mathbb{E}\|\eta^{A^{\alpha_1}g,\nu}(t)\|^2)^{\frac{1}{2}}(\mathbb{E}\|\eta^{A^{\alpha_2}h,\nu}(t)\|^2)^{\frac{1}{2}} + C(\mathbb{E}\|\zeta^{A^{\alpha_1}g,A^{\alpha_2}h,\nu}(t)\|^2)^{\frac{1}{2}}\\
		\le&\ C(1+t^{-\alpha_1}+\|\nu\|_{\gamma}^6)\|g\|\cdot C(1+\|\nu\|_{\gamma}^6)\|h\| + C(1+\|\nu\|_{\gamma}^{18})\|g\|\|h\|\\
		\le&\ C(1+t^{-\alpha_1}+\|\nu\|_{\gamma}^{18})\|g\|\|h\|.
	\end{align*}
	Thus the proof is completed.
\end{proof}

In a similar way, we can derive a priori estimate for $D^3U^N(t,\nu)$, and we omit the proof.

\begin{proposition}\label{p3.4.5}
	Assume (A1)-(A3). For any given $\phi \in \mathcal{C}_b^3(H)$, there exists a constant $C>0$ independent of $N$, such that for all $\nu,h,g,l\in H_N$ and $t \in [0,T]$
	\begin{equation*}
		|D^3U^N(t,\nu)(h,g,A^{-\frac{1}{4}-\delta}l)| \le C(1+\|\nu\|_{\gamma}^{24})\|h\|\|g\|\|l\|,
	\end{equation*}
	where $\delta\in(0,\frac{\gamma}{2}-\frac{1}{4})$ is a small parameter.
\end{proposition}

\section{Weak convergence analysis for spatial semi-discretization}
\label{Sect5}
In this section, we analyze the weak convergence order of the semi-discretized spectral Galerkin approximation. 

\begin{theorem}\label{th2.5.1}
	Assume (A1)-(A3). Then for any given $\phi \in \mathcal{C}_b^2(H)$, there exists a constant $C>0$ independent of $N$, such that
	\begin{equation*}
		|\mathbb{E}[\phi(X(T))]-\mathbb{E}[\phi(X^{N}(T))]| \le CN^{-2\gamma}.
	\end{equation*}
\end{theorem}

\begin{proof}
	For any $M>N>0$, by \eqref{eq3.3.1} we have
	\begin{equation}\label{eq3.5.3}
		\begin{aligned}
			\|(P_M-P_{N})A^{-\gamma}\|_{\mathcal{L}(H)} 
			\le&\ \|(P_M-I)A^{-\gamma}\|_{\mathcal{L}(H)} + \|(I-P^{N})A^{-\gamma}\|_{\mathcal{L}(H)}\\ 
			\le&\ C(M^{-2\gamma}+N^{-2\gamma}).
		\end{aligned}
	\end{equation}
	For any $\phi \in \mathcal{C}_b^2(H)$, according to Lemma \ref{l3.3.2}, there holds
	\begin{equation*}
		\begin{aligned}
			&\ |\mathbb{E}[\phi(X(T))]-\mathbb{E}[\phi(X^{N}(T))]|\\
			\le&\  |\mathbb{E}[\phi(X(T))]-\mathbb{E}[\phi(X^{M}(T))]| + |\mathbb{E}[\phi(X^M(T))]-\mathbb{E}[\phi(X^{N}(T))]|\\
			\le&\ |\phi|_{\mathcal{C}_b^1(H)}(\mathbb{E}\|X(T)-X^M(T)\|^2)^{\frac{1}{2}} + |\mathbb{E}[\phi(X^M(T))]-\mathbb{E}[\phi(X^{N}(T))]|\\
			\le&\ CM^{-\gamma} + |\mathbb{E}[\phi(X^M(T))]-\mathbb{E}[\phi(X^{N}(T))]|.
		\end{aligned}
	\end{equation*}
	For any $\mathcal{F}_{T}$-measurable, $H_M$-valued random variable $\eta$, there holds, cf.\ \cite[(5.1)]{Andersson2016},
	\begin{equation*}
		\mathbb{E}[\phi(\eta)] = \mathbb{E}[\mathbb{E}[\phi(\eta)|\mathcal{F}_T]] = \mathbb{E}[\mathbb{E}[\phi(X^M(0,\eta))|\mathcal{F}_T]] = \mathbb{E}[U^M(0,\eta)].
	\end{equation*}
	Noticing that $X^N(t)\in H_N\subset H_M$, we split the weak error as
	\begin{equation}\label{eq3.5.2}
		\begin{aligned}
			&\ |\mathbb{E}[\phi(X^M(T))] - \mathbb{E}[\phi(X^N(T))]| = |\mathbb{E}[U^M(T,P_M\xi)] - \mathbb{E}[U^M(0,X^N(T))]|\\
			\le&\ |\mathbb{E}[U^M(T,P_M\xi)] - \mathbb{E}[U^M(T,P_N\xi)]|
			+ |\mathbb{E}[U^M(T,P_N\xi)] - \mathbb{E}[U^M(0,X^N(T))]|.
		\end{aligned}
	\end{equation}
	
	By mean value theorem, \eqref{eq7.2.4}, Theorem \ref{t3.4.3} and \eqref{eq3.5.3}, we have
	\begin{align*}
		&\ |\mathbb{E}[U^M(T,P_M\xi)] - \mathbb{E}[U^M(T,P_N\xi)]|\\
		=&\ \big|\mathbb{E}\int_0^1 DU^M(T,\theta P_M\xi+(1-\theta)P_N\xi)(P_M-P_N)\xi d\theta\big|\\
		=&\ \big|\mathbb{E}\int_0^1 \big( DU^M(T,\theta P_M\xi+(1-\theta)P_N\xi),(P_M-P_N)\xi\big) d\theta\big|\\
		\le&\ \mathbb{E}\int_0^1 \big|\big(A^{\frac{\gamma}{2}}DU^M(T,\theta P_M\xi+(1-\theta)P_{N}\xi),(P_M-P_{N})A^{-\gamma}A^{\frac{\gamma}{2}}\xi\big)\big|d\theta\\
		\le&\ C(1+T^{-\frac{\gamma}{2}}+\|\xi\|_{\gamma}^6)\cdot \|(P_M-P_{N})A^{-\gamma}\|_{\mathcal{L}(H)}\|\xi\|_{\gamma}\\
		\le&\ C(M^{-2\gamma}+N^{-2\gamma}).
	\end{align*}
	
	Now we estimate the second term in \eqref{eq3.5.2}.
	Applying It\^o's formula to $U^M(T-t,X^{N}(t))$ for $t \in [0,T]$, we get
	\begin{equation}\label{eq3.5.1}
		\begin{aligned}
			&\ U^M(0,X^{N}(T)) - U^M(T,X^{N}(0))\\
			=&\ \int_0^T -\frac{\partial U^M}{\partial t}(T-t,X^N(t)) + (DU^M(T-t,X^N(t)),-AX^N(t)+P_NF(X^N(t)))dt\\
			&+ \int_0^T\big(DU^M(T-t, X^{N}(t)),P_NG(X^{N}(t))dW(t)\big)\\
			&+ \frac{1}{2}\int_0^T\mathrm{Tr}\{D^2U^M(T-t,X^{N}(t))(P_NG(X^{N}(t))Q^{\frac{1}{2}})(P_NG(X^{N}(t))Q^{\frac{1}{2}})^*\}dt.
		\end{aligned}
	\end{equation}
	According to Kolmogorov equation \eqref{eq3.4.1}, we have
	\begin{equation*}
		\begin{aligned}
			&\ \frac{\partial U^{M}}{\partial t}(T-t,X^{N}(t)) = (DU^M(T-t,X^{N}(t)), -AX^{N}(t)+P_MF(X^{N}(t)))\\
			&+ \frac{1}{2}\mathrm{Tr}\{D^2U^M(T-t,X^{N}(t))(P_MG(X^{N}(t))Q^{\frac{1}{2}})(P_MG(X^{N}(t))Q^{\frac{1}{2}})^*\}.
		\end{aligned}
	\end{equation*}
	Substituting this into \eqref{eq3.5.1}, we obtain
	\begin{align*}
		&\ |\mathbb{E}[U^M(T,P_N\xi)] - \mathbb{E}[U^M(0,X^N(T))]|\\ 
		\le&\ \big|\mathbb{E}\int_{0}^{T}(DU^M(T-t, X^{N}(t)), (P_N-P_M)F(X^N(t))) d t\big|\\
		&+ \frac{1}{2}\big| \mathbb{E}\int_0^T \mathrm{Tr}\big\{D^2U^M(T-t,X^{N}(t))(P_NG(X^{N}(t))Q^{\frac{1}{2}})(P_NG(X^{N}(t))Q^{\frac{1}{2}})^*\\
		&\quad\quad\quad-D^2U^M(T-t,X^{N}(t))(P_MG(X^{N}(t))Q^{\frac{1}{2}})(P_MG(X^{N}(t))Q^{\frac{1}{2}})^*\big\} dt\big|\\
		=:&\ e^1 + e^2.
	\end{align*}
	According to Banach algebra theory, there holds
	\begin{equation*}
		\|A^{\frac{\gamma}{2}}F(X^N(t))\| = \|X^N(t)-(X^N(t))^3\|_{\gamma} \le \|X^N(t)\|_{\gamma}+\|X^N(t)\|^3_{\gamma}.
	\end{equation*}
	By Theorem \ref{t3.4.3} and \eqref{eq3.5.3}, we have
	\begin{align*}
		e^1
		\le&\ \mathbb{E}\int_{0}^{T}\big|\big(A^{\frac{\gamma}{2}}DU^M(T-t, X^{N}(t)), (P_N-P_M)A^{-\gamma}A^{\frac{\gamma}{2}}F(X^N(t))\big)\big| dt\\
		\le&\ \mathbb{E}\int_{0}^{T} \|A^{\frac{\gamma}{2}}DU^M(T-t, X^{N}(t))\|\|(P_N-P_M)A^{-\gamma}\|_{\mathcal{L}(H)}\|A^{\frac{\gamma}{2}}F(X^N(t))\|dt\\
		\le&\ \mathbb{E}\int_{0}^{T} C(1+(T-t)^{-\frac{\gamma}{2}}+\|X^N(t)\|_{\gamma}^6)\cdot C(M^{-2\gamma}+N^{-2\gamma}) \cdot(\|X^N(t)\|_{\gamma}+\|X^N(t)\|^3_{\gamma})dt\\
		\le&\ C(M^{-2\gamma}+N^{-2\gamma}).
	\end{align*}
	Next we estimate $e^2$. Noticing that
	\begin{equation*}
		\begin{aligned}
			&\ P_NG(X^{N}(t))QG^*(X^{N}(t))P_N - P_MG(X^{N}(t))QG^*(X^{N}(t))P_M\\
			=&\ P_NG(X^{N}(t))QG^*(X^{N}(t))P_N - P_NG(X^{N}(t))QG^*(X^{N}(t))P_M\\ 
			&\ + P_NG(X^{N}(t))QG^*(X^{N}(t))P_M - P_MG(X^{N}(t))QG^*(X^{N}(t))P_M\\
			=&\ P_NG(X^{N}(t))QG^*(X^{N}(t))(P_N-P_M) + (P_N-P_M)G(X^{N}(t))QG^*(X^{N}(t))P_M.
		\end{aligned}
	\end{equation*}
	Then we get
	\begin{equation*}
		\begin{aligned}
			e^2
			\le&\ \frac{1}{2} \mathbb{E}\int_0^T\big| \mathrm{Tr}\big\{D^2U^M(T-t,X^{N}(t))P_NG(X^{N}(t))QG^*(X^{N}(t))(P_N-P_M)\big\}\big| dt\\
			&\ + \frac{1}{2} \mathbb{E}\int_0^T\big| \mathrm{Tr}\big\{D^2U^M(T-t,X^{N}(t))(P_N-P_M)G(X^{N}(t))QG^*(X^{N}(t))P_M\big\}\big| dt\\
			=:&\ e^{2,1}+e^{2,2}.
		\end{aligned}
	\end{equation*}
	For any $\delta\in(0,\frac{\gamma}{2}-\frac{1}{4})$, there holds $\dot{H}^{\gamma}\hookrightarrow\dot{H}^{\frac{1}{2}+2\delta}$. By \eqref{eq7.2.3}, Lemma \ref{l3.2.1}, Theorem \ref{t3.4.3} together with \eqref{eq3.5.3}, we have 
	\begin{align*}
		&\ e^{2,1} = \frac{1}{2} \mathbb{E}\int_0^T\big| \mathrm{Tr}\big\{A^{\frac{\gamma}{2}}D^2U^M(T-t,X^{N}(t))A^{-\frac{1}{4}-\delta}\\
		& \quad A^{\frac{1}{4}+\delta} P_NG(X^{N}(t))QG^*(X^{N}(t))A^{\frac{\gamma}{2}}(P_N-P_M)A^{-\gamma}\big\}\big|dt \\
		\le&\ \frac{1}{2} \mathbb{E}\int_0^T \|A^{\frac{\gamma}{2}}D^2U^M(T-t,X^{N}(t))A^{-\frac{1}{4}-\delta}\|_{\mathcal{L}(H)}\\
		& \quad \cdot\big|\mathrm{Tr}\big\{A^{\frac{1}{4}+\delta}P_NG(X^{N}(t))QG^*(X^{N}(t))A^{\frac{\gamma}{2}}\big\}\big|\cdot\|(P_N-P_M)A^{-\gamma}\|_{\mathcal{L}(H)}dt \\
		\le&\ \mathbb{E}\int_0^TC(1+(T-t)^{-\frac{\gamma}{2}}+\|X^N(t)\|_{\gamma}^{18})\\
		&\quad\cdot\big|\mathrm{Tr}\big\{A^{\frac{1}{4}+\delta}P_NG(X^{N}(t))Q^{\frac{1}{2}}Q^{\frac{1}{2}}G^*(X^{N}(t))A^{\frac{\gamma}{2}}\big\}\big|\cdot C(M^{-2\gamma}+N^{-2\gamma})dt\\
		\le&\ \mathbb{E}\int_0^TC(1+(T-t)^{-\frac{\gamma}{2}}+\|X^N(t)\|_{\gamma}^{18})\cdot \|A^{\frac{1}{4}+\delta}P_NG(X^{N}(t))Q^{\frac{1}{2}}\|_{\mathcal{L}_2(H)}\\
		&\quad \cdot\|(A^{\frac{\gamma}{2}}G(X^{N}(t))Q^{\frac{1}{2}})^*\|_{\mathcal{L}_2(H)} \cdot(M^{-2\gamma}+N^{-2\gamma})dt\\
		\le&\ \mathbb{E}\int_0^TC(1+(T-t)^{-\frac{\gamma}{2}}+\|X^N(t)\|_{\gamma}^{18})\cdot\|A^{\frac{\gamma}{2}}G(X^{N}(t))\|_{\mathcal{L}_2^0}^2 \cdot(M^{-2\gamma}+N^{-2\gamma})dt\\
		\le&\ \mathbb{E}\int_0^TC(1+(T-t)^{-\frac{\gamma}{2}}+\|X^N(t)\|_{\gamma}^{18})\cdot C(1+\|X^N(t)\|_\gamma)^2 \cdot(M^{-2\gamma}+N^{-2\gamma})dt\\
		\le&\ C(M^{-2\gamma}+N^{-2\gamma}).
	\end{align*}
	Similarly, we obtain
	\begin{equation*}
		e^{2,2} \le C(M^{-2\gamma}+N^{-2\gamma}).
	\end{equation*}
	To sum up, we deduce that
	\begin{equation*}
		e^2 \le C(M^{-2\gamma}+N^{-2\gamma}).
	\end{equation*}
	Thus we have
	\begin{equation*}
		|\mathbb{E}[\phi(X(T))]-\mathbb{E}[\phi(X^{N}(T))]| \le CM^{-\gamma} + C(M^{-2\gamma}+N^{-2\gamma}).
	\end{equation*}
	The proof is finished by letting $M\to \infty$ in above equation.
\end{proof}

\begin{remark}
	We emphasize that the spatial weak convergence order is twice of the regularity order of the exact mild solution. This is the optimal weak convergence rate.
\end{remark}

\section{Weak convergence analysis for full-discretization}
\label{Sect6}

In this section, we study the optimal weak convergence order of the full discretization. We begin with establishing a priori estimates for the Malliavin derivative of the numerical solution, then investigate a weak error estimate between the spatial semi-discretization and the fully discretized solution. Finally we derive the weak convergence order between the exact and numerical solutions.

\subsection{Malliavin regularity of full-discretization}

In this section, we prove the boundedness of the Malliavin derivative of the numerical solution $X^N_k$ and its continuous version $\widetilde{X}(t)$. For any $t\in[0,T]$, denote by $\lfloor t \rfloor_{\tau} = \max\{k:k\tau \le t\}$. We first study a priori estimates for the Malliavin derivative of $X^N_k$.

\begin{lemma}\label{l3.6.1}
	Assume (A1)-(A3) and $0<\tau < \frac{1}{2C_F}$. Then for any integer $p\ge1$, there exists a constant $C>0$ independent of both $\tau$ and $N$, such that for all $0< s\le T$ and $0\le k\le K-1$,
	\begin{align*}
		\|\mathcal{D}_sX_{k+1}^N\|_{L^p(\Omega;\mathcal{L}_2^0)} \le&\ C,\\
		\|A^{\frac{\gamma}{2}}\mathcal{D}_sX_{k+1}^N\|_{L^p(\Omega;\mathcal{L}_2^{0})} \le&\ C.
	\end{align*}
\end{lemma}

\begin{proof}
	For any $u\in U_0$ and $0\le k\le K-1$, denote by $\mathcal{D}_s^uX_{k+1}^N = \mathcal{D}_sX_{k+1}^N u$ the Malliavin derivative of $X_{k+1}^N$ in the direction $u$. Let $1_{(a,b]}(s)$ denote the indicator function of interval $(a,b]$, then for any $s\in(0,T]$,
	\begin{equation*}
		\mathcal{D}_s^u\Delta W_k = \mathcal{D}_s^u\Big(\int_0^T1_{(0,t_{k+1}]}(r)dW(r)-\int_0^T1_{(0,t_{k}]}(r)dW(r)\Big) = (1_{(0,t_{k+1}]}(s) - 1_{(0,t_{k}]}(s))u.
	\end{equation*}
	Noticing that $X_{k}^N$ is $\mathcal{F}_{t_k}$-adapted, then taking Malliavin derivative in the direction $u$ on both sides of \eqref{eq3.3.7} yields  for $0< s\le t_k$,
	\begin{equation}\label{eq3.6.11.0}
		\begin{aligned}
			\mathcal{D}_s^uX_{k+1}^N =&\ \mathcal{D}_s^uX_{k}^N - \tau A_N\mathcal{D}_s^uX_{k+1}^N + \tau P_NDF(X_{k+1}^N)\mathcal{D}_s^uX_{k+1}^N\\ 
			&\ + P_NDG(X_{k}^N)\mathcal{D}_s^uX_{k}^N\Delta W_k + P_NG(X_k^N)(1_{(0,t_{k+1}]}(s) - 1_{(0,t_{k}]}(s))u\\
			=&\ \mathcal{D}_s^uX_{k}^N - \tau A_N\mathcal{D}_s^uX_{k+1}^N + \tau P_NDF(X_{k+1}^N)\mathcal{D}_s^uX_{k+1}^N + P_NDG(X_{k}^N)\mathcal{D}_s^uX_{k}^N\Delta W_k,
		\end{aligned}
	\end{equation}
	for $t_k< s\le t_{k+1}$, 
	\begin{equation}\label{eq3.6.11}
		\begin{aligned}
			\mathcal{D}_s^uX_{k+1}^N =&\ -\tau A_N\mathcal{D}_s^uX_{k+1}^N + \tau P_NDF(X_{k+1}^N)\mathcal{D}_s^uX_{k+1}^N\\ 
			&\ + P_NG(X_k^N)(1_{(0,t_{k+1}]}(s) - 1_{(0,t_{k}]}(s))u\\
			=&\ -\tau A_N\mathcal{D}_s^uX_{k+1}^N + \tau P_NDF(X_{k+1}^N)\mathcal{D}_s^uX_{k+1}^N + P_NG(X_k^N)u,
		\end{aligned}
	\end{equation}
	and for $s>t_{k+1}$, $\mathcal{D}_s^uX_{k+1}^N\equiv 0$.
	
	As $0< s\le t_k$, testing \eqref{eq3.6.11.0} with $\mathcal{D}_s^uX_{k+1}^N$, we obtain
	\begin{equation}\label{eq3.6.12}
		\begin{aligned}
			&\ (\mathcal{D}_s^uX_{k+1}^N-\mathcal{D}_s^uX_{k}^N,\mathcal{D}_s^uX_{k+1}^N) + \tau (A_N\mathcal{D}_s^uX_{k+1}^N,\mathcal{D}_s^uX_{k+1}^N)\\ 
			=&\  \tau (P_NDF(X_{k+1}^N)\mathcal{D}_s^uX_{k+1}^N,\mathcal{D}_s^uX_{k+1}^N) + (P_NDG(X_{k}^N)\mathcal{D}_s^uX_{k}^N\Delta W_k,\mathcal{D}_s^uX_{k+1}^N).
		\end{aligned}
	\end{equation}
	By $(u-v,u) = \frac{1}{2}(\|u\|^2-\|v\|^2) + \frac{1}{2}\|u-v\|^2$, we have
	\begin{equation*}
		(\mathcal{D}_s^uX_{k+1}^N-\mathcal{D}_s^uX_{k}^N,\mathcal{D}_s^uX_{k+1}^N) = \frac{1}{2}(\|\mathcal{D}_s^uX_{k+1}^N\|^2-\|\mathcal{D}_s^uX_{k}^N\|^2) + \frac{1}{2}\|\mathcal{D}_s^uX_{k+1}^N-\mathcal{D}_s^uX_{k}^N\|^2.
	\end{equation*}
	Note that $(A_N\mathcal{D}_s^uX_{k+1}^N,\mathcal{D}_s^uX_{k+1}^N) = (A_N^{\frac{1}{2}}\mathcal{D}_s^uX_{k+1}^N,A_N^{\frac{1}{2}}\mathcal{D}_s^uX_{k+1}^N) \ge 0$. By \eqref{eq3.2.5} , we get
	\begin{equation*}
		(P_NDF(X_{k+1}^N)\mathcal{D}_s^uX_{k+1}^N,\mathcal{D}_s^uX_{k+1}^N) \le (DF(X_{k+1}^N)\mathcal{D}_s^uX_{k+1}^N,\mathcal{D}_s^uX_{k+1}^N) \le C_F\|\mathcal{D}_s^uX_{k+1}^N\|^2.
	\end{equation*}
	Moreover
	\begin{align*}
		&\ (P_NDG(X_{k}^N)\mathcal{D}_s^uX_{k}^N\Delta W_k,\mathcal{D}_s^uX_{k+1}^N) = (DG(X_{k}^N)\mathcal{D}_s^uX_{k}^N\Delta W_k,\mathcal{D}_s^uX_{k+1}^N)\\
		=&\ (DG(X_{k}^N)\mathcal{D}_s^uX_{k}^N\Delta W_k,\mathcal{D}_s^uX_{k+1}^N-\mathcal{D}_s^uX_{k}^N) + (DG(X_{k}^N)\mathcal{D}_s^uX_{k}^N\Delta W_k,\mathcal{D}_s^uX_{k}^N)\\
		\le&\ \frac{1}{2}\|DG(X_{k}^N)\mathcal{D}_s^uX_{k}^N\Delta W_k\|^2 + \frac{1}{2}\|\mathcal{D}_s^uX_{k+1}^N-\mathcal{D}_s^uX_{k}^N\|^2 + (DG(X_{k}^N)\mathcal{D}_s^uX_{k}^N\Delta W_k,\mathcal{D}_s^uX_{k}^N).
	\end{align*}
	Substituting these estimates into \eqref{eq3.6.12} yields
	\begin{equation}\label{eq3.6.4}
		\begin{aligned}
			&\ \frac{1}{2}\|\mathcal{D}_s^uX_{k+1}^N\|^2-\frac{1}{2}\|\mathcal{D}_s^uX_{k}^N\|^2\\ 
			\le&\ C_F\tau\|\mathcal{D}_s^uX_{k+1}^N\|^2 + \frac{1}{2}\|DG(X_{k}^N)\mathcal{D}_s^uX_{k}^N\Delta W_k\|^2 + (DG(X_{k}^N)\mathcal{D}_s^uX_{k}^N\Delta W_k,\mathcal{D}_s^uX_{k}^N).
		\end{aligned}
	\end{equation}
	As $t_k< s\le t_{k+1}$, testing \eqref{eq3.6.11} with $\mathcal{D}_s^uX_{k+1}^N$ and by \eqref{eq3.2.5}, we have
	\begin{align*}
		\|\mathcal{D}_s^uX_{k+1}^N\|^2
		=&\ -\tau (A_N\mathcal{D}_s^uX_{k+1}^N,\mathcal{D}_s^uX_{k+1}^N) + \tau (P_NDF(X_{k+1}^N)\mathcal{D}_s^uX_{k+1}^N,\mathcal{D}_s^uX_{k+1}^N)\\ 
		&\ + (P_NG(X_k^N)u,\mathcal{D}_s^uX_{k+1}^N)\\
		\le&\ \tau (DF(X_{k+1}^N)\mathcal{D}_s^uX_{k+1}^N,\mathcal{D}_s^uX_{k+1}^N) + \frac{1}{2}\|P_NG(X_k^N)u\|^2+\frac{1}{2}\|\mathcal{D}_s^uX_{k+1}^N\|^2\\
		\le&\ C_F\tau\|\mathcal{D}_s^uX_{k+1}^N\|^2 + \frac{1}{2}\|G(X_k^N)u\|^2 +\frac{1}{2}\|\mathcal{D}_s^uX_{k+1}^N\|^2\\
		=&\ (\frac{1}{2}+C_F\tau)\|\mathcal{D}_s^uX_{k+1}^N\|^2 + \frac{1}{2}\|G(X_k^N)u\|^2,
	\end{align*}
	which implies for small $\tau>0$, 
	\begin{equation*}
		\|\mathcal{D}_s^uX_{k+1}^N\|^2 \le C\|G(X_k^N)u\|^2.
	\end{equation*}
	this together with Lemma \ref{l3.2.1} leads to
	\begin{equation}\label{eq3.6.5}
		\begin{aligned}
			&\ \mathbb{E}\|\mathcal{D}_sX_{k+1}^N\|^{2p}_{\mathcal{L}_2^0} = \mathbb{E}\Big(\sum_{j=1}^{\infty}\|\mathcal{D}_sX_{k+1}^NQ^{\frac{1}{2}}\varphi_j\|^{2}\Big)^p = \mathbb{E}\Big(\sum_{j=1}^{\infty}\|\mathcal{D}_s^{Q^{\frac{1}{2}}\varphi_j}X_{k+1}^N\|^{2}\Big)^p\\ 
			\le&\ \mathbb{E}\Big(\sum_{j=1}^{\infty}C\|G(X_k^N)Q^{\frac{1}{2}}\varphi_j\|^2\Big)^p= C\mathbb{E}\|G(X_k^N)\|^{2p}_{\mathcal{L}_2^0} \le C\mathbb{E}C_G(1+\|X_k^N\|)^{2p} \le C.
		\end{aligned}
	\end{equation}
	
	For any given $0\le k\le K-1$ and $0<s\le t_{k+1}$, we have $t_{\lfloor s\rfloor }<s\le t_{\lfloor s\rfloor +1}$ and 
	\begin{equation}\label{eq3.6.6}
		\mathcal{D}_s^uX_{i}^N = 0,\quad  \forall \ 0\le i\le \lfloor s\rfloor .
	\end{equation}
	By induction on $k\ge \lfloor s\rfloor +1$ in \eqref{eq3.6.4}, we have
	\begin{align*}
		\|\mathcal{D}_s^uX_{k+1}^N\|^2 \le&\ \|\mathcal{D}_s^uX_{\lfloor s\rfloor +1}^N\|^2 + \sum_{i=\lfloor s\rfloor +1}^kC\tau\|\mathcal{D}_s^uX_{i+1}^N\|^2 + \sum_{i=\lfloor s\rfloor +1}^k\|DG(X_{i}^N)\mathcal{D}_s^uX_{i}^N\Delta W_i\|^2\\ 
		&\ + \sum_{i=\lfloor s\rfloor +1}^k2(DG(X_{i}^N)\mathcal{D}_s^uX_{i}^N\Delta W_i,\mathcal{D}_s^uX_{i}^N).
	\end{align*}
	By a continuous version Burkholder-Davis-Gundy inequality, we have
	\begin{equation*}
		\mathbb{E}\|\int_{t_i}^{t_{i+1}}DG(X_{i}^N)\mathcal{D}_s^uX_{i}^NdW(r)\|^{2p} \le C(p)\tau^{p}\mathbb{E}\|DG(X_{i}^N)\mathcal{D}_s^uX_{i}^N\|_{\mathcal{L}_2^0}^{2p}.
	\end{equation*}
	According to \eqref{eq3.6.6} and Lemma \ref{l3.2.00}, we obtain
	\begin{align*}
		&\ \mathbb{E}\|\mathcal{D}_s^uX_{k+1}^N\|^{2p} \le 4^{p-1}\|\mathcal{D}_s^uX_{\lfloor s\rfloor +1}^N\|^{2p} + 4^{p-1}\tau^p(k+1)^{p-1}\sum_{i=0}^k\mathbb{E}\|\mathcal{D}_s^uX_{i+1}^N\|^{2p}\\ 
		&\ + 4^{p-1}(k+1)^{p-1}\sum_{i=0}^k\mathbb{E}\|\int_{t_i}^{t_{i+1}}DG(X_{i}^N)\mathcal{D}_s^uX_{i}^NdW(r)\|^{2p}\\ 
		&\ + 4^{p-1}\mathbb{E}\big(\sum_{i=0}^k(DG(X_{i}^N)\mathcal{D}_s^uX_{i}^N\Delta W_i,\mathcal{D}_s^uX_{i}^N)\big)^p\\
		\le&\ C(p)\|\mathcal{D}_s^uX_{\lfloor s\rfloor +1}^N\|^{2p} + C(p)\tau T^{p-1}\sum_{i=0}^k\mathbb{E}\|\mathcal{D}_s^uX_{i+1}^N\|^{2p}\\ 
		&\ + C(p)\tau T^{p-1}\sum_{i=0}^k\mathbb{E}\|DG(X_{i}^N)\mathcal{D}_s^uX_{i}^N\|^{2p}_{\mathcal{L}_2^0}\\ 
		&\ + C(p)\mathbb{E}\big(\sum_{i=0}^k\tau\|DG(X_{i}^N)\mathcal{D}_s^uX_{i}^N\|^2_{\mathcal{L}_2^0}\|\mathcal{D}_s^uX_{i}^N\|^2\big)^{\frac{p}{2}}\\
		\le&\ C\|\mathcal{D}_s^uX_{\lfloor s\rfloor +1}^N\|^{2p} + C\tau\sum_{i=0}^k\mathbb{E}\|\mathcal{D}_s^uX_{i+1}^N\|^{2p} + C\tau\sum_{i=0}^k\mathbb{E}\|\mathcal{D}_s^uX_{i}^N\|^{2p}\\ 
		&\ + C(p)(k+1)^{\frac{p}{2}-1}\tau^{\frac{p}{2}}\sum_{i=0}^kC_G\|\mathcal{D}_s^uX_{i}^N\|^{2p}\\
		\le&\ C\|\mathcal{D}_s^uX_{\lfloor s\rfloor +1}^N\|^{2p} + C\tau\sum_{i=0}^{k+1}\mathbb{E}\|\mathcal{D}_s^uX_{i}^N\|^{2p}.
	\end{align*}
	From a discrete Gronwall's lemma, it follows that
	\begin{equation}\label{eq6.5.1}
		\mathbb{E}\|\mathcal{D}_s^uX_{k+1}^N\|^{2p} \le C\mathbb{E}\|\mathcal{D}_s^uX_{\lfloor s\rfloor +1}^N\|^{2p}e^{C(k+1)\tau} \le C\mathbb{E}\|\mathcal{D}_s^uX_{\lfloor s\rfloor +1}^N\|^{2p}.
	\end{equation}
	This together with \eqref{eq3.6.5} implies the first assertion.
	
	For any given $0\le k\le K-1$ and $0<s\le t_{k+1}$, taking Malliavin derivative in direction $u$ on both sides of \eqref{eq3.3.10}, we get
	\begin{align*}
		\mathcal{D}_s^uX_{k+1}^N =&\ \tau\sum_{i=\lfloor s\rfloor}^{k}S_{N,\tau}^{k+1-i}P_NDF(X_{i+1}^N)\mathcal{D}_s^uX_{i+1}^N\\ 
		&\ + \sum_{i=\lfloor s\rfloor+1}^{k}S_{N,\tau}^{k+1-i}P_NDG(X_i^N)\mathcal{D}_s^uX_{i}^N\Delta W_i+ S_{N,\tau}^{k+1-\lfloor s\rfloor }P_NG(X_{\lfloor s\rfloor }^N)u.
	\end{align*}
	By \eqref{eq3.3.4}, \eqref{eq3.2.5}, \eqref{eq3.3.15}, H\"older inequality, Burkholder-Davis-Gundy inequality in Lemma \ref{l3.2.00} and Lemma \ref{l3.4.2}, we have
	\begin{align*}
		&\ \|\mathcal{D}_s^uX_{k+1}^N\|_{L^{2p}(\Omega;\dot{H}^{\gamma})}
		\le \tau\sum_{i=\lfloor s\rfloor}^{k}\|A^{\frac{\gamma}{2}}S_{N,\tau}^{k+1-i}P_NDF(X_{i+1}^N)\mathcal{D}_s^uX_{i+1}^N\|_{L^{2p}(\Omega;H)}\\ 
		&\  + \Big(\mathbb{E}\Big\|\sum_{i=\lfloor s\rfloor+1}^{k}A^{\frac{\gamma}{2}}S_{N,\tau}^{k+1-i}P_NDG(X_i^N)\mathcal{D}_s^uX_{i}^N\Delta W_i\Big\|^{2p}\Big)^{\frac{1}{2p}} + \|A^{\frac{\gamma}{2}}S_{N,\tau}^{k+1-\lfloor s\rfloor }P_NG(X_{\lfloor s\rfloor }^N)u\|_{L^{2p}(\Omega;H)}\\
		\le&\ \tau\sum_{i=\lfloor s\rfloor}^{k}\|A^{\frac{\gamma}{2}}S_{N,\tau}^{k+1-i}P_N\|_{\mathcal{L}(H)}\|DF(X_{i+1}^N)\mathcal{D}_s^uX_{i+1}^N\|_{L^{2p}(\Omega;H)}\\  
		&\ + C(p)\Big(\sum_{i=\lfloor s\rfloor+1}^{k}\big(\mathbb{E}\|A^{\frac{\gamma}{2}}S_{N,\tau}^{k+1-i}P_NDG(X_i^N)\mathcal{D}_s^uX_{i}^N\Delta W_i\|^{2p}\big)^{\frac{1}{p}}\Big)^{\frac{1}{2}}\\
		&\ + \|S_{N,\tau}^{k+1-\lfloor s\rfloor }P_N\|_{\mathcal{L}(H)}\|A^{\frac{\gamma}{2}}G(X_{\lfloor s\rfloor }^N)u\|_{L^{2p}(\Omega;H)}\\
		\le&\ \tau\sum_{i=\lfloor s\rfloor}^{k}Ct_{k+1-i}^{-\frac{\gamma}{2}}\big(\mathbb{E}[C_F(1+\|X_{i+1}^N\|_V^2)^{2p}\|\mathcal{D}_s^uX_{i+1}^N\|^{2p}]\big)^{\frac{1}{2p}}\\ 
		&\ + C\Big(\sum_{i=\lfloor s\rfloor+1}^{k}\tau t^{-\gamma}_{k+1-i}\big(\mathbb{E}C_G\|\mathcal{D}_s^uX_{i}^N\|^{2p}\big)^{\frac{1}{p}}\Big)^{\frac{1}{2}} + C\|A^{\frac{\gamma}{2}}G(X_{\lfloor s\rfloor }^N)u\|_{L^{2p}(\Omega;H)}\\
		\le&\ \tau\sum_{i=\lfloor s\rfloor}^{k}Ct_{k+1-i}^{-\frac{\gamma}{2}}\cdot C\|\mathcal{D}_s^uX_{\lfloor s\rfloor +1}^N\|_{L^{2p}(\Omega;H)} + C\Big(\sum_{i=\lfloor s\rfloor+1}^{k}\tau t^{-\gamma}_{k+1-i}\cdot C\|\mathcal{D}_s^uX_{\lfloor s\rfloor +1}^N\|_{L^{2p}(\Omega;H)}^2\Big)^{\frac{1}{2}}\\ 
		&\ + C\|A^{\frac{\gamma}{2}}G(X_{\lfloor s\rfloor }^N)u\|_{L^{2p}(\Omega;H)}.
	\end{align*}
	This together with \eqref{eq6.5.1} implies the second assertion.
\end{proof}

Now we study a priori estimates for the continuous version $\widetilde{X}(t)$ of $X_k^N$.

\begin{lemma}\label{t3.4.1}
	Assume (A1)-(A3) and $0<\tau < \frac{1}{2C_F}$. Then for any integer $p\ge 1$, there exists a constant $C > 0$ independent of both $\tau$ and $N$, such that 
	\begin{align*}
		\|\mathcal{D}_s \widetilde{X}(t)\|_{L^p(\Omega;\mathcal{L}_2^0)} + \|A^{\frac{\gamma}{2}}\mathcal{D}_s \widetilde{X}(t)\|_{L^p(\Omega;\mathcal{L}_2^0)} \le&\ C,\quad \forall \ 0 < s,t \le T,\\
		\|(\mathcal{D}^-\widetilde{X})_t\|_{L^p(\Omega;\mathcal{L}_2^0)} \le&\ C\tau,\quad \forall \ 0\le t\le T.
	\end{align*}
\end{lemma}

\begin{proof}
	For any $t \in (0,T]$, there exists a unique integer $k$ with $0 \le k \le K-1$ such that $t \in (t_k, t_{k+1}]$.  
	Taking the Malliavin derivative in direction $u\in U_0$ on both sides of \eqref{eq3.3.9}, we obtain, for $0 < s \le t_k \le t \le t_{k+1}$,
	\begin{align*}
		\mathcal{D}_s^u\widetilde{X}(t) =&\ \mathcal{D}_s^uX_k^N - (t-t_k)A_NS_{N,\tau}\mathcal{D}_s^uX_k^N +  (t-t_k)S_{N,\tau}P_NDF(X_{k+1}^N)\mathcal{D}_s^uX_{k+1}^N\\ 
		&\ + \int_{t_k}^tS_{N,\tau}P_NDG(X_{k}^N)\mathcal{D}_s^uX_{k}^NdW(r),
	\end{align*} 
	for $t_k< s\le t\le t_{k+1}$,
	\begin{equation*}
		\mathcal{D}_s^u\widetilde{X}(t) = (t-t_k)S_{N,\tau}P_NDF(X_{k+1}^N)\mathcal{D}_s^uX_{k+1}^N + S_{N,\tau}P_NG(X_k^N)u,
	\end{equation*}
	and for $t_k\le t< s\le t_{k+1}$,
	\begin{align*}
		\mathcal{D}_s^u\widetilde{X}(t) =&\ (t-t_k)S_{N,\tau}P_NDF(X_{k+1}^N)\mathcal{D}_s^u(X_{k+1}^N) + S_{N,\tau}P_NG(X_k^N)\mathcal{D}_s^u(W(t)-W(t_k))\\
		=&\ (t-t_k)S_{N,\tau}P_NDF(X_{k+1}^N)\mathcal{D}_s^u(X_{k+1}^N).
	\end{align*}
	By Lemma \ref{l3.6.1}, we can obtain
	\begin{equation*}
		\|\mathcal{D}_s \widetilde{X}(t)\|_{L^p(\Omega;\mathcal{L}_2^0)} + \|A^{\frac{\gamma}{2}}\mathcal{D}_s\widetilde{X}(t)\|_{L^p(\Omega;\mathcal{L}_2^0)} \le C.
	\end{equation*}
	
	Specifically, according to \eqref{eq6.2.13}, we have
	\begin{equation*}
		(\mathcal{D}^-\widetilde{X})_t = (t-t_k)S_{N,\tau}P_NDF(X_{k+1}^N)\mathcal{D}_sX_{k+1}^N,\quad\forall\ t_k< t\le t_{k+1},\ 0\le k\le K-1.
	\end{equation*}
	By \eqref{eq3.3.4}, \eqref{eq3.2.5}, H\"older inequality and Lemmas \ref{l3.4.2} and \ref{l3.6.1}, we have
	\begin{align*}
		&\ \mathbb{E}\|S_{N,\tau}P_NDF(X_{k+1}^N)\mathcal{D}_sX_{k+1}^N\|^p_{\mathcal{L}_2^0} \le \|S_{N,\tau}P_N\|^p_{\mathcal{L}(H)}\mathbb{E}\big(\|DF(X_{k+1}^N)\|_{\mathcal{L}(H)}^p\|\mathcal{D}_sX_{k+1}^N\|^p_{\mathcal{L}_2^0}\big)\\
		\le&\ C\mathbb{E}\big(C_F(1+\|X_{k+1}^N\|_V^2)^p\|\mathcal{D}_sX_{k+1}^N\|_{\mathcal{L}_2^0}^p\big) \le C\big(\mathbb{E}(1+\|X_{k+1}^N\|_V^2)^{2p}\big)^{\frac{1}{2}}\big(\mathbb{E}\|\mathcal{D}_sX_{k+1}^N\|_{\mathcal{L}_2^0}^{2p}\big)^{\frac{1}{2}} \le C.
	\end{align*}
	This implies the second assertion.
\end{proof}

\subsection{Weak convergence order for full discretization}

We first analyze the weak error between the spectral Galerkin semi-discretization and the fully discrete scheme,  
and then derive weak error estimates between the exact and numerical solutions.

\begin{theorem}\label{th2.6.1}
	Assume (A1)-(A3) and $0<\tau < \frac{1}{2C_F}$.  
	Then for any $\phi\in \mathcal{C}_{b}^{3}(H)$, there exists a constant $C>0$ independent of both $\tau$ and $N$, such that
	\begin{equation*}
		\bigl|\mathbb{E}\bigl[\phi\bigl(X^{N}(T)\bigr)\bigr]
		- \mathbb{E}\bigl[\phi\bigl(X_{K}^{N}\bigr)\bigr]\bigr|
		\le C\,\tau^{\gamma}.
	\end{equation*}
\end{theorem}

\begin{proof}
	Noting  that $X_0^N = X^N(0) = P_N\xi$, the global weak error can be  converted to the sum of local weak errors as follows
	\begin{equation*}
		\begin{aligned}
			&\ |\mathbb{E}[\phi(X^N(T))]-\mathbb{E}[\phi(X_K^N)]|\\
			=&\ |\mathbb{E}[U^N(T,X_0^N)] - \mathbb{E}[U^N(0,X_K^N)]|\\
			=&\ \big|\sum_{k=0}^{K-1}\big(\mathbb{E}[U^N(T-t_k,X_k^N)] - \mathbb{E}[U^N(T-t_{k+1},X_{k+1}^N)]\big)\big|\\
			\le&\ \sum_{k=0}^{K-1}|\mathbb{E}[U^N(T-t_k,\widetilde{X}(t_{k})) - U^N(T-t_{k+1},\widetilde{X}(t_{k+1}))]|.
		\end{aligned}
	\end{equation*}
	Applying It\^o's formula for Skorohod integral in \eqref{eq6.2.14} to $U^N(T-t,\widetilde{X}(t))$ for $t \in [t_{k},t_{k+1}]$ and then taking $t=t_{k+1}$, we obtain
	\begin{equation}\label{eq3.6.0}
		\begin{aligned}
			&U^N(T-t_{k+1},\widetilde{X}(t_{k+1})) = U^N(T-t_k,\widetilde{X}(t_k)) - \int_{t_k}^{t_{k+1}}\frac{\partial U^N}{\partial t}(T-t,\widetilde{X}(t))dt\\
			&+ \int_{t_k}^{t_{k+1}}\big(DU^N(T-t,\widetilde{X}(t)),-A_NS_{N,\tau}X_{k}^N + S_{N,\tau}P_NF(X_{k+1}^N)\big)dt\\
			&+ \int_{t_k}^{t_{k+1}}\big(DU^N(T-t,\widetilde{X}(t)), S_{N,\tau}P_NG(X_k^N)dW(t)\big)\\
			&+ \frac{1}{2}\int_{t_k}^{t_{k+1}}\!\!\!\!\mathrm{Tr}\big\{D^2U^N(T-t,\widetilde{X}(t))(S_{N,\tau}P_NG(X_k^N)Q^{\frac{1}{2}})(S_{N,\tau}P_NG(X_k^N)Q^{\frac{1}{2}})^*\big\}dt\\
			&+ \int_{t_k}^{t_{k+1}}\mathrm{Tr}\big\{D^2U^N(T-t,\widetilde{X}(t))((\mathcal{D}^-\widetilde{X})_tQ^{\frac{1}{2}})(S_{N,\tau}P_NG(X_k^N)Q^{\frac{1}{2}})^*\big\}dt.
		\end{aligned}
	\end{equation}
	The Kolmogorov equation \eqref{eq3.4.1} yields
	\begin{equation*}
		\begin{aligned}
			&\ \frac{\partial U^N}{\partial t}(T-t,\widetilde{X}(t)) = (DU^N(T-t,\widetilde{X}(t)), -A_N\widetilde{X}(t)+P_NF(\widetilde{X}(t)))\\
			&+ \frac{1}{2}\mathrm{Tr}\{D^2U^N(T-t,\widetilde{X}(t))(P_NG(\widetilde{X}(t))Q^{\frac{1}{2}})(P_NG(\widetilde{X}(t))Q^{\frac{1}{2}})^*\}.
		\end{aligned}
	\end{equation*}
	Substitute this into \eqref{eq3.6.0}, we get
	\begin{equation}\label{eq3.6.2}
		\begin{aligned}
			&\ |\mathbb{E}[U^N(T-t_k,\widetilde{X}(t_{k})) - U^N(T-t_{k+1},\widetilde{X}(t_{k+1}))]|\\
			\le&\ |\mathbb{E}\int_{t_k}^{t_{k+1}}\big(DU^N(T-t,\widetilde{X}(t)),A_N\widetilde{X}(t)-A_NS_{N,\tau}X_{k}^N\big)dt|\\
			&\ + |\mathbb{E}\int_{t_k}^{t_{k+1}}\big(DU^N(T-t,\widetilde{X}(t)),S_{N,\tau}P_NF(X_{k+1}^N)-P_NF(\widetilde{X}(t))\big)dt|\\
			&\ + \frac{1}{2}|\mathbb{E}\int_{t_k}^{t_{k+1}}\mathrm{Tr}\big\{D^2U^N(T-t,\widetilde{X}(t))\\
			&\ \quad\cdot\big(S_{N,\tau}P_NG(X_k^N)Q(S_{N,\tau}P_NG(X_k^N))^* - P_NG(\widetilde{X}(t))Q(P_NG(\widetilde{X}(t)))^*\big)\big\}dt|\\
			&\ + |\mathbb{E}\int_{t_k}^{t_{k+1}}\mathrm{Tr}\big\{D^2U^N(T-t,\widetilde{X}(t))((\mathcal{D}^-\widetilde{X})_tQ^{\frac{1}{2}})(S_{N,\tau}P_NG(X_k^N)Q^{\frac{1}{2}})^*\big\}dt|\\
			=:&\ a_k + b_k + c_k + d_k.
		\end{aligned}
	\end{equation}
	
	The detailed analyses of the estimates for $a_{k}$, $b_{k}$, $c_{k}$, and $d_{k}$ are deferred to  
	Appendices A–D, respectively; here we only state the resulting estimates.
	\begin{equation*}
		\sum_{k=0}^{K-1}a_k \le C\tau^{\gamma},\quad \sum_{k=0}^{K-1}b_k \le C\tau^{\gamma},\quad \sum_{k=0}^{K-1}c_k \le C\tau^{\gamma},\quad \sum_{k=0}^{K-1}d_k \le C\tau.
	\end{equation*}
	Thus we finally obtain
	\begin{equation*}
		|\mathbb{E}[\phi(X^{N}(T))]-\mathbb{E}[\phi(X_K^N)]| \le \sum_{k=0}^{K-1}(a_k+b_k+c_k+d_k) \le C\tau^{\gamma}.
	\end{equation*}
	This completes the proof.
\end{proof}

\begin{remark}
	We notice that the temporal weak convergence order of numerical solution is twice of the temporal H\"older continuity order of exact mild solution.
\end{remark}

Now we are ready to derive the weak convergence order of the full-discretization.

\begin{theorem}\label{t3.6.2}
	Assume (A1)-(A3) and $0<\tau < \frac{1}{2C_F}$.  
	Then for any $\phi\in \mathcal{C}_{b}^{3}(H)$, there exists a constant $C>0$ independent of both $\tau$ and $N$, such that
	\begin{equation*}
		|\mathbb{E}[\phi(X(T))]-\mathbb{E}[\phi(X_K^N)]| \le C(\tau^{\gamma}+N^{-2\gamma}).
	\end{equation*}
\end{theorem}

\begin{proof}
	According to Theorems \ref{th2.5.1} and \ref{th2.6.1}, we have
	\begin{equation*}
		\begin{aligned}
			&\ |\mathbb{E}[\phi(X(T))]-\mathbb{E}[\phi(X_K^N)]|\\
			\le&\ |\mathbb{E}[\phi(X(T))]-\mathbb{E}[\phi(X^N(T))]| + |\mathbb{E}[\phi(X^N(T))]-\mathbb{E}[\phi(X_K^N)]|\\
			\le&\ CN^{-2\gamma}+C\tau^{\gamma}.
		\end{aligned}
	\end{equation*}
	Thus the proof is completed.
\end{proof}

\section{Numerical experiment}
\label{Sect7}

Consider a stochastic Allen--Cahn equation with Nemytskii type diffusion for $(t,x)\in (0,1)\times(0,1)$
\begin{equation*}
	\begin{aligned}
		dX(t)-\Delta X(t)dt =&\ (X(t)-X^3(t))dt+\sin(X(t))dW(t),\\
		X(0,x) =&\ x^{0.2501}(1-x)^{0.2501} \in \dot{H}^{\frac{3}{4}},
	\end{aligned}
\end{equation*}
where the covariance operator $Q$ of $W(t)$ has eigenvalues $q_j=j^{-2.5001}$.  Then, (A3) is fulfilled with $\gamma=\frac{3}{4}$. The test function is chosen as $\phi(X)=\sin(\|X\|)$. 

\begin{figure}[H]
	\centering
	\includegraphics[width=0.48\textwidth]{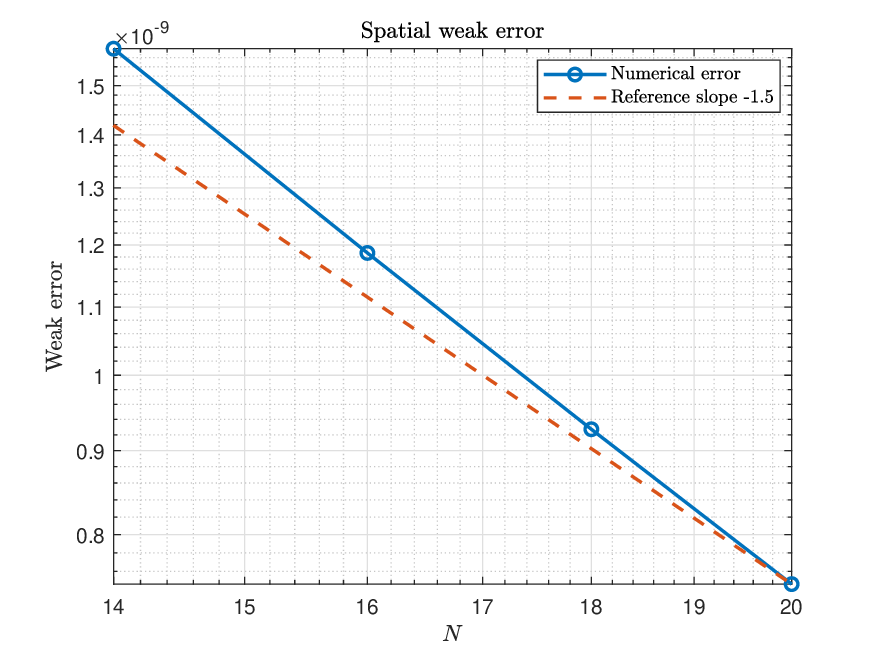}
	\includegraphics[width=0.48\textwidth]{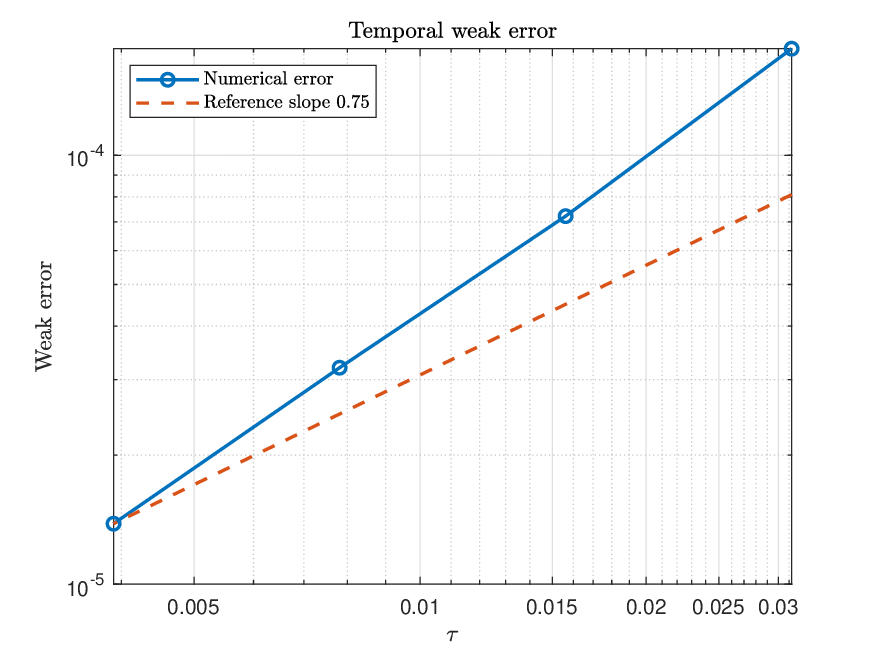}
	\caption{Temporal and spatial weak convergence orders}
	\label{figure6.1}
\end{figure}

We fix a time step size $\tau = \frac{1}{1000}$ and compute $M=1000$ independent sample trajectories for $N_i = 14,16,18,20$, respectively, and regard a numerical solution computed with $N_{\text{ref}} = 64$ as the reference 'exact' solution. We detect the spatial weak errors and plot them in log-log scale, see the left hand side in Figure \ref{figure6.1}. Then, we fix  $N=64$ and  compute $M=1000$ independent sample trajectories for each time step size $\tau_j = 2^{-j}$ $(j=5,6,7,8)$, respectively, and regard a numerical solution computed with $\tau_{\text{ref}} = 2^{-11}$ as the reference 'exact' solution. We compute the temporal weak errors and plot them in log-log scale, see the right hand side in Figure \ref{figure6.1}. We find that the weak convergence orders are consistent with predicted theoretical rates.

\section*{Funding}
The work is supported by National Natural Science Foundation of China (12171199, 12571428) and the Project for Enhancing the Research and Innovation Capacity of Doctoral Candidates at Jilin University (2025KC028).

\section*{Acknowledgement}  
The authors would like to thank the Editor and the anonymous reviewers for their careful reading and valuable comments, which have greatly improved the quality of this manuscript.

\bibliographystyle{plain}
\bibliography{references}
\appendix

\section{Estimation for $a_k$ in \eqref{eq3.6.2}}
\vskip1mm

From \eqref{eq3.3.3.1}, it follows that $ A_N(I-S_{N,\tau}) = \tau S_{N,\tau}A_N^2$. Then we have
\begin{equation*}
	\begin{aligned}
		&\ A_N\widetilde{X}(t)-A_NS_{N,\tau}X_{k}^N
		= A_N\widetilde{X}(t)-A_NX_{k}^N+A_NX_{k}^N-A_NS_{N,\tau}X_{k}^N\\
		=&\ A_N(\widetilde{X}(t)-X_{k}^N) + A_N(I-S_{N,\tau})X_{k}^N
		= A_N(\widetilde{X}(t)-X_{k}^N) + \tau S_{N,\tau}A_N^2X_{k}^N.
	\end{aligned}
\end{equation*}
Thus
\begin{equation}\label{app.1}
	\begin{aligned}
		a_k \le&\ |\mathbb{E}\int_{t_k}^{t_{k+1}}( DU^N(T-t,\widetilde{X}(t)), A_N(\widetilde{X}(t)-X_{k}^N))dt|\\
		&\ + |\mathbb{E}\int_{t_k}^{t_{k+1}}(DU^N(T-t,\widetilde{X}(t)), \tau S_{N,\tau}A_N^2X_{k}^N)dt|:= a_k^1 + a_k^2.
	\end{aligned}
\end{equation}

{\bf (a.1) Estimation of $a_k^{1}$}.
By \eqref{eq3.3.9.1}, we have
\begin{align*}
	a_k^1 \le&\ \tau|\mathbb{E}\int_{t_k}^{t_{k+1}}(DU^N(T-t,\widetilde{X}(t)), -A_N^2S_{N,\tau}X_{k}^N)dt|\\
	&\ + \tau|\mathbb{E}\int_{t_k}^{t_{k+1}}(DU^N(T-t,\widetilde{X}(t)), A_NS_{N,\tau}P_NF(X_{k+1}^N))dt|\\
	&\ + |\mathbb{E}\int_{t_k}^{t_{k+1}}(DU^N(T-t,\widetilde{X}(t)), \int_{t}^{t_{k+1}}A_NS_{N,\tau}P_NG(X_{k}^N)dW(s))dt|\\
	=:&\ a_k^{1,1} + a_k^{1,2} + a_k^{1,3}.
\end{align*}

{\bf (a.1.1) Estimation of $a_k^{1,1}$}.
In terms of \eqref{eq3.3.10}, we get
\begin{align*}
	a_k^{1,1}
	\le&\ \tau|\mathbb{E}\int_{t_k}^{t_{k+1}}(DU^N(T-t,\widetilde{X}(t)), -A_N^2S_{N,\tau}S_{N,\tau}^kP_N\xi)dt|\\
	&\ + \tau|\mathbb{E}\int_{t_k}^{t_{k+1}}(DU^N(T-t,\widetilde{X}(t)), -A_N^2S_{N,\tau}\tau\sum_{i=0}^{k-1}S_{N,\tau}^{k-i}P_NF(X_{i+1}^N))dt|\\
	&\ + \tau|\mathbb{E}\int_{t_k}^{t_{k+1}}\!\!(DU^N(T-t,\widetilde{X}(t)), -A_N^2S_{N,\tau}\!\!\sum_{i=0}^{k-1}\int_{t_i}^{t_{i+1}}\!\!\!\!S_{N,\tau}^{k-i}P_NG(X_i^N)dW(s))dt|\\
	=&\ a_k^{1,1,1} + a_k^{1,1,2} + a_k^{1,1,3}.
\end{align*}
According to Theorem \ref{t3.4.3}, \eqref{eq3.3.4} and \eqref{eq3.3.14'r}, we have for any $0<\varepsilon<\frac{1-\gamma}{2}$,
\begin{align*}
	&\ a_k^{1,1,1}\\ 
	=&\ \tau|\mathbb{E}\int_{t_k}^{t_{k+1}}(DU^N(T-t,\widetilde{X}(t)), A_N^2S_{N,\tau}^{k+1}P_N\xi)dt|\\
	=&\ \tau|\mathbb{E}\int_{t_k}^{t_{k+1}}(A_N^{\frac{\gamma}{2}+\varepsilon}DU^N(T-t,\widetilde{X}(t)), A_N^{1-\gamma}S_{N,\tau}A_N^{1-\varepsilon}S_{N,\tau}^{k}P_NA_N^{\frac{\gamma}{2}}\xi)dt|\\
	\le&\ \tau\mathbb{E}\int_{t_k}^{t_{k+1}}C(1+(T-t)^{-\frac{\gamma}{2}-\varepsilon}+\|\widetilde{X}(t)\|^6_{\gamma})\cdot \|A_N^{1-\gamma}S_{N,\tau}\|_{\mathcal{L}(H)} \|A_N^{1-\varepsilon}S_{N,\tau}^{k+1}P_N\|_{\mathcal{L}(H)}\|A_N^{\frac{\gamma}{2}}\xi\|dt\\
	\le&\ \tau\int_{t_k}^{t_{k+1}}C\mathbb{E}(1+(T-t)^{-\frac{\gamma}{2}-\varepsilon}+\|\widetilde{X}(t)\|^6_{\gamma})\cdot C\tau^{\gamma-1}\cdot Ct_{k+1}^{-1+\varepsilon}\cdot Cdt\\
	\le&\ C\tau^{\gamma} t_{k+1}^{-1+\varepsilon}\int_{t_k}^{t_{k+1}}1+(T-t)^{-\frac{\gamma}{2}-\varepsilon}dt.
\end{align*}
Due to Theorem \ref{t3.4.3}, \eqref{eq3.2.5}, \eqref{eq3.3.14'r}, \eqref{eq3.3.15}, Lemma \ref{l3.4.2} and Banach algebra property of $\dot{H}^{\gamma}$, we obtain
\begin{align*}
	&\ a_k^{1,1,2}\\
	=&\ \tau|\mathbb{E}\int_{t_k}^{t_{k+1}}\tau\sum_{i=0}^{k-1}(A_N^{\frac{\gamma}{2}+\varepsilon}DU^N(T-t,\widetilde{X}(t)), A_N^{1-\gamma}S_{N,\tau}A_N^{1-\varepsilon}S_{N,\tau}^{k-i}A_N^{\frac{\gamma}{2}}P_NF(X_{i+1}^N))dt|\\
	\le&\ \tau\mathbb{E}\int_{t_k}^{t_{k+1}}\tau\sum_{i=0}^{k-1}C(1+(T-t)^{-\frac{\gamma}{2}-\varepsilon}+\|\widetilde{X}(t)\|^6_{\gamma})\\
	&\quad\quad\quad\cdot\|A_N^{1-\gamma}S_{N,\tau}\|_{\mathcal{L}(H)}\|A_N^{1-\varepsilon}S_{N,\tau}^{k-i}\|_{\mathcal{L}(H)}\|A_N^{\frac{\gamma}{2}}P_NF(X_{i+1}^N)\|dt\\
	\le&\ \tau\mathbb{E}\int_{t_k}^{t_{k+1}}\tau\sum_{i=0}^{k-1} C(1+(T-t)^{-\frac{\gamma}{2}-\varepsilon}+\|\widetilde{X}(t)\|^6_{\gamma})\cdot C\tau^{\gamma-1}\cdot Ct_{k-i}^{-1+\varepsilon}\cdot \|F(X_{i+1}^N)\|_{\gamma}dt\\
	\le&\ C\tau^{\gamma}\int_{t_k}^{t_{k+1}}\tau\sum_{i=0}^{k-1}t_{k-i}^{-1+\varepsilon}\cdot (T-t)^{-\frac{\gamma}{2}-\varepsilon}\mathbb{E}\|F(X_{i+1}^N)\|_\gamma + \mathbb{E}\big((1+\|\widetilde{X}(t)\|^6_{\gamma})\|F(X_{i+1}^N)\|_\gamma\big)dt\\
	\le&\ C\tau^{\gamma}\int_{t_k}^{t_{k+1}}C\cdot(T-t)^{-\frac{\gamma}{2}-\varepsilon}\cdot C+Cdt\\ 
	\le&\  C\tau^{\gamma}\int_{t_k}^{t_{k+1}}1+(T-t)^{-\frac{\gamma}{2}-\varepsilon}dt.
\end{align*}
By Malliavin integration by parts formula in Lemma \ref{l3.2.4}, Theorems \ref{t3.4.4} and \ref{t3.4.1}, \eqref{eq3.3.15}, Lemma \ref{l3.2.1}, \eqref{eq3.3.14'r} and Lemma \ref{l3.4.2}, we get
\begin{align*}
	a_k^{1,1,3} 
	=&\ \tau|\mathbb{E}\int_{t_k}^{t_{k+1}}\sum_{i=0}^{k-1}(DU^N(T-t,\widetilde{X}(t)), \int_{t_i}^{t_{i+1}}A_N^2S_{N,\tau}^{k+1-i}P_NG(X_i^N)dW(s))dt|\\
	=&\ |\mathbb{E}\int_{t_k}^{t_{k+1}}\!\!\!\!\int_{t_i}^{t_{i+1}}\!\!\!\!\tau\sum_{i=0}^{k-1}\big(D^2U^N(T-t,\widetilde{X}(t))\mathcal{D}_s\widetilde{X}(t), A_N^2S_{N,\tau}^{k+1-i}P_NG(X_i^N)\big)_{\mathcal{L}_2^0}dsdt|\\
	=&\ |\mathbb{E}\int_{t_k}^{t_{k+1}}\!\!\!\!\int_{t_i}^{t_{i+1}}\!\!\!\!\tau\sum_{i=0}^{k-1}\big(A_N^{\frac{\gamma}{2}+\varepsilon}D^2U^N(T\!-\!t,\widetilde{X}(t))A_N^{-\frac{1}{4}-\delta}A_N^{\frac{1}{4}+\delta}\mathcal{D}_s\widetilde{X}(t),\\
	&\quad\quad\quad\quad\quad\quad A_N^{1-\gamma}S_{N,\tau}A_N^{1-\varepsilon}S_{N,\tau}^{k-i}A_N^{\frac{\gamma}{2}}P_NG(X_i^N)\big)_{\mathcal{L}_2^0}dsdt|\\
	\le&\
	\mathbb{E}\int_{t_k}^{t_{k+1}}\!\!\!\!\int_{t_i}^{t_{i+1}}\tau\sum_{i=0}^{k-1}\|A^{\frac{\gamma}{2}+\varepsilon}D^2U^N(T-t,\widetilde{X}(t))A^{-\frac{1}{4}-\delta}\|_{\mathcal{L}(H)}\|A^{\frac{1}{4}+\delta}\mathcal{D}_s\widetilde{X}(t)\|_{\mathcal{L}_2^0}\\
	&\quad\quad\quad\quad\quad\quad\cdot\|A^{1-\gamma}S_{N,\tau}\|_{\mathcal{L}(H)} \|A^{1-\varepsilon}S_{N,\tau}^{k-i}\|_{\mathcal{L}(H)}\|A^{\frac{\gamma}{2}}G(X_i^N)\|_{\mathcal{L}_2^0}dsdt\\
	\le&\
	\mathbb{E}\int_{t_k}^{t_{k+1}}\!\!\!\!\int_{t_i}^{t_{i+1}}\tau\sum_{i=0}^{k-1} C(1+(T-t)^{-\frac{\gamma}{2}-\varepsilon}+\|\widetilde{X}(t)\|^{18}_{\gamma})\cdot C\|A^{\frac{\gamma}{2}}\mathcal{D}_s\widetilde{X}(t)\|_{\mathcal{L}_2^0}\\
	&\quad\quad\quad\quad\quad\quad\cdot C\tau^{\gamma-1}\cdot Ct_{k-i}^{-1+\varepsilon}\cdot C\|G(X_i^N)\|_{\mathcal{L}_2^\gamma}dsdt\\
	\le&\ C\tau^{\gamma}\int_{t_k}^{t_{k+1}}1+(T-t)^{-\frac{\gamma}{2}-\varepsilon}dt.
\end{align*}
Thus we have
\begin{equation}\label{app.2}
	\sum_{k=0}^{K-1}a_k^{1,1} \le C\tau^{\gamma}\int_{0}^{T}1+(T-t)^{-\frac{\gamma}{2}-\varepsilon}dt \le C\tau^{\gamma},
\end{equation} 
where $C>0$ is a constant independent of both $\tau$ and $N$.

{\bf (a.1.2) Estimation of $a_k^{1,2}$}.
By virtue of Theorem \ref{t3.4.3}, \eqref{eq3.3.4}, \eqref{eq3.2.5}, H\"older inequality, \eqref{eq3.3.14'r} and Lemma \ref{l3.4.2}, we have
\begin{equation*}
	\begin{aligned}
		a_k^{1,2} =&\ \tau|\mathbb{E}\int_{t_k}^{t_{k+1}}\!\!\!(A_N^{\frac{\gamma}{2}}DU^N(T-t,\widetilde{X}(t)), A_N^{1-\gamma}S_{N,\tau}A_N^{\frac{\gamma}{2}}P_NF(X_{k+1}^N))dt|\\
		\le&\ \tau\mathbb{E}\int_{t_k}^{t_{k+1}}\!\!\!C(1+(T-t)^{-\frac{\gamma}{2}}+\|\widetilde{X}(t)\|^6_{\gamma})\cdot\|A_N^{1-\gamma}S_{N,\tau}\|_{\mathcal{L}(H)}\|A_N^{\frac{\gamma}{2}}P_NF(X_{k+1}^N)\|dt\\
		\le&\ \tau \mathbb{E}\int_{t_k}^{t_{k+1}}\!\!\! C(1+(T-t)^{-\frac{\gamma}{2}}+\|\widetilde{X}(t)\|^6_{\gamma})\cdot C\tau^{\gamma-1}\cdot\|F(X_{k+1}^N)\|_\gamma dt\\
		\le&\ C\tau^{\gamma}\int_{t_k}^{t_{k+1}}1+(T-t)^{-\frac{\gamma}{2}}dt.
	\end{aligned}
\end{equation*}

{\bf (a.1.3) Estimation of $a_k^{1,3}$}. By Malliavin integration by parts formula in Lemma \ref{l3.2.4}, Theorems \ref{t3.4.4}, \ref{t3.4.1}, \eqref{eq3.3.4}, Lemma \ref{l3.2.1}, Lemma \ref{l3.4.2} and \eqref{eq3.3.14'r}, we obtain
\begin{align*}
	a_k^{1,3} =&\ |\mathbb{E}\int_{t_k}^{t_{k+1}}\!\!\!\!\int_{t}^{t_{k+1}}\!\!\!(D^2U^N(T-t,\widetilde{X}(t))\mathcal{D}_s\widetilde{X}(t),A_NS_{N,\tau}P_NG(X_{k}^N))_{\mathcal{L}_2^0}dsdt|\\
	=&\ |\mathbb{E}\int_{t_k}^{t_{k+1}}\!\!\!\!\int_{t}^{t_{k+1}}\!\!\!(A_N^{\frac{\gamma}{2}}D^2U^N(T-t,\widetilde{X}(t))A_N^{-\frac{1}{4}-\delta}A_N^{\frac{1}{4}+\delta}\mathcal{D}_s\widetilde{X}(t), A_N^{1-\gamma}S_{N,\tau}A_N^{\frac{\gamma}{2}}P_NG(X_{k}^N))_{\mathcal{L}_2^0}dsdt|\\
	\le&\ \mathbb{E}\int_{t_k}^{t_{k+1}}\!\!\!\!\int_{t}^{t_{k+1}}\!\!\!\|A_N^{\frac{\gamma}{2}}D^2U^N(T-t,\widetilde{X}(t))A_N^{-\frac{1}{4}-\delta}\|_{\mathcal{L}(H)} \|A_N^{\frac{1}{4}+\delta}\mathcal{D}_s\widetilde{X}(t)\|_{\mathcal{L}_2^0}\cdot\\
	&\quad\quad\quad\quad\quad\|A_N^{\gamma-1}S_{N,\tau}\|_{\mathcal{L}(H)}\|A_N^{\frac{\gamma}{2}}P_NG(X_{k}^N)\|_{\mathcal{L}_2^0}dsdt\\
	\le&\ \mathbb{E}\int_{t_k}^{t_{k+1}}\!\!\!\!\int_{t}^{t_{k+1}}\!\!\!C(1+(T-t)^{-\frac{\gamma}{2}}+\|\widetilde{X}(t)\|^{18}_{\gamma}) \cdot C\|A^{\frac{\gamma}{2}}\mathcal{D}_s\widetilde{X}(t)\|_{\mathcal{L}_2^0}\cdot C\tau^{\gamma-1}\cdot \|G(X_{k}^N)\|_{\mathcal{L}_2^\gamma}dsdt\\
	\le&\ C\tau^{\gamma}\mathbb{E}\int_{t_k}^{t_{k+1}}(1+(T-t)^{-\frac{\gamma}{2}}+\|\widetilde{X}(t)\|^{18}_{\gamma})\cdot\|A^{\frac{\gamma}{2}}\mathcal{D}_s\widetilde{X}(t)\|_{\mathcal{L}_2^0}\cdot C(1+\|X_{k}^N\|_{\gamma})dt\\
	\le&\ C\tau^{\gamma}\mathbb{E}\int_{t_k}^{t_{k+1}}1+(T-t)^{-\frac{\gamma}{2}}dt.
\end{align*}
This together with the estimates in {\bf (a.1.1)} and {\bf (a.1.2)} implies
\begin{equation*}
	\sum_{k=0}^{K-1}a_k^1 \le \sum_{k=0}^{K-1}(a_k^{1,1}+a_k^{1,2}+a_k^{1,3}) \le C\tau^{\gamma}.
\end{equation*}

{\bf (a.2) Estimation of $a_k^{2}$}. From \eqref{eq3.3.10}, it follows that
\begin{align*}
	a_k^2
	\le&\ |\mathbb{E}\int_{t_k}^{t_{k+1}}\big(DU^N(T-t,\widetilde{X}(t)), \tau S_{N,\tau}A_N^2S_{N,\tau}^kP_N\xi\big)dt|\\
	&\ + |\mathbb{E}\int_{t_k}^{t_{k+1}}\big(DU^N(T-t,\widetilde{X}(t)), \tau S_{N,\tau}A_N^2\tau\sum_{i=1}^{k-1}S_{N,\tau}^{k-i}P_NF(X_{i+1}^N)\big)dt|\\
	&\ + |\mathbb{E}\int_{t_k}^{t_{k+1}}\big(DU^N(T-t,\widetilde{X}(t)), \tau S_{N,\tau}A_N^2\sum_{i=1}^{k-1}\int_{t_i}^{t_{i+1}}S_{N,\tau}^{k-i}P_NG(X_i^N)dW(s)\big)dt|.
\end{align*}
In a similar way as in {\bf (a.1.1)}, we can prove $\displaystyle\sum_{k=0}^{K-1}a_k^2 \le C\tau^{\gamma}$. This together with \eqref{app.1} and \eqref{app.2} omplies
\begin{equation}\label{eq3.6.a}
	\sum_{k=0}^{K-1}a_k \le C\tau^{\gamma},
\end{equation}
where $C>0$ is a constant independent of both $\tau$ and $N$.

\vskip1mm
\section{Estimation for $b_k$ in \eqref{eq3.6.2}}
\vskip1mm

From \eqref{eq3.3.3.1}, it follows that
\begin{equation*}
	\begin{aligned}
		&\ P_NF(\widetilde{X}(t))-S_{N,\tau}P_NF(X_{k+1}^N)\\ 
		=&\ P_NF(\widetilde{X}(t))-S_{N,\tau}P_NF(\widetilde{X}(t)) + S_{N,\tau}P_NF(\widetilde{X}(t))-S_{N,\tau}P_NF(X_{k+1}^N)\\
		=&\ (I-S_{N,\tau})P_NF(\widetilde{X}(t)) + S_{N,\tau}P_N(F(\widetilde{X}(t))-F(X_{k+1}^N))\\
		=&\ \tau A_NS_{N,\tau}P_NF(\widetilde{X}(t)) + S_{N,\tau}P_N(F(\widetilde{X}(t))-F(X_{k+1}^N)).
	\end{aligned}
\end{equation*}
Thus
\begin{equation}\label{app.3}
	\begin{aligned}
		b_k \le&\ \big|\mathbb{E}\int_{t_k}^{t_{k+1}}\big(DU^N(T-t,\widetilde{X}(t)), \tau A_NS_{N,\tau}P_NF(\widetilde{X}(t))\big)dt\big|\\
		&\ + \big|\mathbb{E}\int_{t_k}^{t_{k+1}}\big(DU^N(T-t,\widetilde{X}(t)), S_{N,\tau}P_N(F(\widetilde{X}(t))-F(X_{k+1}^N))\big)dt\big|\\
		=:&\ b_k^1 + b_k^2.
	\end{aligned}
\end{equation}

{\bf (b.1) Estimation of $b_k^{1}$}. By Theorem \ref{t3.4.3}, \eqref{eq3.3.4}, \eqref{eq3.2.5} and \eqref{eq3.3.14'r}, we have
\begin{align*}
	b_k^1 \le&\ \tau\mathbb{E}\int_{t_k}^{t_{k+1}}|(A_N^{\frac{\gamma}{2}}DU^N(T-t,\widetilde{X}(t)),  A_N^{1-\gamma}S_{N,\tau}A_N^{\frac{\gamma}{2}}P_NF(\widetilde{X}(t)))|dt\\
	\le&\ \tau\mathbb{E}\int_{t_k}^{t_{k+1}}C(1+(T-t)^{-\frac{\gamma}{2}}+\|\widetilde{X}(t)\|_{\gamma}^6)\cdot \|A_N^{1-\gamma}S_{N,\tau}\|_{\mathcal{L}(H)}\|A_N^{\frac{\gamma}{2}}P_NF(\widetilde{X}(t))\|dt\\
	\le&\ \tau\mathbb{E}\int_{t_k}^{t_{k+1}}C(1+(T-t)^{-\frac{\gamma}{2}}+\|\widetilde{X}(t)\|_{\gamma}^6)\cdot C\tau^{\gamma-1}\cdot \|F(\widetilde{X}(t))\|_\gamma dt\\
	\le&\ C\tau^{\gamma}\int_{t_k}^{t_{k+1}}1+(T-t)^{-\frac{\gamma}{2}}dt,
\end{align*}
which implies
\begin{equation}\label{app.4}
	\sum_{k=0}^{K-1}b_k^1 \le C\tau^{\gamma}\int_{0}^{T}1+(T-t)^{-\frac{\gamma}{2}}dt \le C\tau^{\gamma}.
\end{equation}

{\bf (b.2) Estimation of $b_k^{2}$}. By virtue of the mean value theorem and \eqref{eq3.3.9.1}, we obtain
\begin{align*}
	b_k^2 =&\ \mathbb{E}\big|\!\int_{t_k}^{t_{k+1}}\!\!\!\!\big(DU^N(T-t,\widetilde{X}(t)),S_{N,\tau}P_N\int_0^1\!\!DF(\mu(t))(\widetilde{X}(t)-X_{k+1}^N)d\theta\big)dt\big|\\
	\le&\ \mathbb{E}\big|\int_{t_k}^{t_{k+1}}\!\!\!\big(DU^N(T-t,\widetilde{X}(t)), S_{N,\tau}P_N\int_0^1\!\!DF(\mu(t))\!\!\int_{t}^{t_{k+1}}\!\!\!A_NS_{N,\tau}X_{k+1}^Ndrd\theta\big)dt\big|\\
	+&\ \mathbb{E}\big|\int_{t_k}^{t_{k+1}}\!\!\!\!\big(DU^N(T-t,\widetilde{X}(t)), S_{N,\tau}P_N\int_0^1\!\!DF(\mu(t))\!\!\int_{t}^{t_{k+1}}\!\!\!\!S_{N,\tau}P^NF(X_{k+1}^N) drd\theta\big)dt\big|\\
	+&\ \mathbb{E}\big|\int_{t_k}^{t_{k+1}}\!\!\!\!\!\!\big(DU^N(T\!-\!t,\widetilde{X}(t)),S_{N,\tau}P_N\!\!\int_0^1\!\!DF(\mu(t))\!\!\int_{t}^{t_{k+1}}\!\!\!\!\!S_{N,\tau}P^NG(X_k^N)dW(r)d\theta \big)dt\big|,
\end{align*}
where $\mu(t) := \theta\widetilde{X}(t)+(1-\theta)X_{k+1}^N$.
Similarly to the analysis for $a_k^1$, we obtain $\displaystyle\sum_{k=0}^{K-1}b_k^2 \le C\tau^{\gamma}$.
This together with \eqref{app.3} and \eqref{app.4} leads to
\begin{equation}\label{eq3.6.b}
	\sum_{k=0}^{K-1}b_k \le C\tau^{\gamma}.
\end{equation}

\vskip1mm
\section{Estimation for $c_k$ in \eqref{eq3.6.2}}
\vskip1mm

Noticing that
\begin{align*}
	&\ S_{N,\tau}P_NG(X_k^N)Q(S_{N,\tau}P_NG(X_k^N))^* - P_NG(\widetilde{X}(t))Q(P_NG(\widetilde{X}(t)))^*\\
	=&\ [S_{N,\tau}P_NG(X_k^N) - P_NG(\widetilde{X}(t))]Q(P_NG(\widetilde{X}(t)))^*\\
	&\ + S_{N,\tau}P_NG(X_k^N)Q[S_{N,\tau}P_NG(X_k^N) - P_NG(\widetilde{X}(t))]^*\\
	=&\ S_{N,\tau}P_N[G(X_k^N)-G(\widetilde{X}(t))]Q^{\frac{1}{2}}(P_NG(\widetilde{X}(t))Q^{\frac{1}{2}})^*\\
	&\ + [S_{N,\tau}-I]P_NG(\widetilde{X}(t))Q^{\frac{1}{2}}(P_NG(\widetilde{X}(t))Q^{\frac{1}{2}})^*\\
	&\ + S_{N,\tau}P_NG(X_k^N)Q^{\frac{1}{2}}\big(S_{N,\tau}P_N[G(X_k^N)-G(\widetilde{X}(t))]Q^{\frac{1}{2}}\big)^*\\ 
	&\ + S_{N,\tau}P_NG(X_k^N)Q^{\frac{1}{2}}\big([S_{N,\tau}-I]P_NG(\widetilde{X}(t))Q^{\frac{1}{2}}\big)^*\\
	=:&\ \tilde{c}_k^1+\tilde{c}_k^2+\tilde{c}_k^3+\tilde{c}_k^4.
\end{align*}
Thus
\begin{equation}\label{app.5}
	c_k \le \sum_{\ell=1}^4\frac{1}{2}|\mathbb{E}\int_{t_k}^{t_{k+1}}\mathrm{Tr}\{D^2U^N(T-t,\widetilde{X}(t))\tilde{c}_k^\ell\}dt|
	:=\sum_{\ell=1}^4c_k^\ell.
\end{equation}

{\bf (c.1) Estimation of $c_k^{1}$ and $c_k^{3}$}. We mainly focus on the analysis for $c_k^{1}$, and then the estimate for $c_k^{3}$ follows from a similar way. 

Noticing \eqref{eq3.3.9} and applying It\^o's formula for Skorohod integral in \eqref{eq6.2.14} to $G(\widetilde{X}(t))$ for $t \in[t_k,t_{k+1}]$, we have for $h_1,h_2\in H$,
\begin{align*}
	&\ \big(G(\widetilde{X}(t))-G(X_k^N)h_1,h_2\big)\\
	=&\ -\int_{t_k}^tD\big(G(\widetilde{X}(s))h_1,h_2\big)A_NS_{N,\tau}X_{k}^Nds+\int_{t_k}^tD\big(G(\widetilde{X}(s))h_1,h_2\big)S_{N,\tau}P_NF(X_{k+1}^N)ds\\
	&\ + \int_{t_k}^tD\big(G(\widetilde{X}(s))h_1,h_2\big)S_{N,\tau}P_NG(X_{k}^N)dW(s)\\
	&\ + \frac{1}{2}\int_{t_k}^t\mathrm{Tr}\{D^2\big(G(\widetilde{X}(s))h_1,h_2\big)S_{N,\tau}P_NG(X_{k}^N)Q(S_{N,\tau}P_NG(X_{k}^N))^*\}ds\\
	&\ + \int_{t_k}^t\mathrm{Tr}\{D^2\big(G(\widetilde{X}(s))h_1,h_2\big)((\mathcal{D}^-\widetilde{X})_sQ^{\frac{1}{2}})(S_{N,\tau}P_NG(X_{k}^N)Q^{\frac{1}{2}})^*\}ds\\
	=:&\ (\mathcal{A}h_1,h_2)+(\mathcal{B}h_1,h_2)+(\mathcal{C}h_1,h_2)+(\mathcal{E}h_1,h_2)+(\mathcal{G}h_1,h_2).
\end{align*}
Thus, $c_k^1$ is further decomposed as the sum of $c^{1,\Pi}_k$ for $\Pi\in\{\mathcal{A},\mathcal{B},\mathcal{C},\mathcal{E},\mathcal{G}\}$, where
\begin{equation*}
	c^{1,\Pi}_k=: \frac{1}{2}|\mathbb{E}\int_{t_k}^{t_{k+1}}\mathrm{Tr}\{D^2U^N(T-t,\widetilde{X}(t))S_{N,\tau}P_N\Pi Q^{\frac{1}{2}}(P_NG(\widetilde{X}(t))Q^{\frac{1}{2}})^*\}dt|.
\end{equation*}

{\bf (c.1.1) Estimation of $c_k^{1,\mathcal{A}}$}. By \eqref{eq3.3.10}, we have
\begin{align*}
	\mathcal{A} =&\  -\int_{t_k}^tDG(\widetilde{X}(s))A_NS_{N,\tau}S^{k-1}_{N,\tau}P_N\xi ds\\ 
	&\ - \int_{t_k}^tDG(\widetilde{X}(s))A_NS_{N,\tau}\tau\sum_{i=1}^{k-1}S^{k-i}_{N,\tau}P_NF(X_{i+1}^N)ds\\
	&\ -\int_{t_k}^tDG(\widetilde{X}(s))A_NS_{N,\tau}\sum_{i=1}^{k-1}S^{k-i}_{N,\tau}P_NG(X_{i}^N)\Delta W_ids\\
	=:&\ \mathcal{A}_1+\mathcal{A}_2+\mathcal{A}_3.
\end{align*}
Thus $c^{1,\mathcal{A}}_k$ is further decomposed to
\begin{align*}
	c^{1,\mathcal{A}}_k \le&\ \sum_{\ell=1}^3\frac{1}{2}|\mathbb{E}\int_{t_k}^{t_{k+1}}\mathrm{Tr}\{D^2U^N(T-t,\widetilde{X}(t))S_{N,\tau}P_N\mathcal{A}_{\ell} Q^{\frac{1}{2}}(P_NG(\widetilde{X}(t))Q^{\frac{1}{2}})^*\}dt|\\
	=:&\ \sum_{\ell=1}^3\frac{1}{2}c^{1,\mathcal{A}_\ell}_k.
\end{align*}
We will analyze these three terms, respectively.

{\bf (c.1.1.1) Estimation of $c_k^{1,\mathcal{A}_1}$}. By \eqref{eq7.2.3}, \eqref{eq3.3.4}, Lemma \ref{l3.2.1}, Theorem \ref{t3.4.4}, \eqref{eq3.3.15} and \eqref{eq3.3.14'r}, we obtain
\begin{align*}
	c^{1,\mathcal{A}_1}_k
	\le&\ \frac{1}{2}\mathbb{E}\int_{t_k}^{t_{k+1}}\!\!\!\int_{t_k}^t|\mathrm{Tr}\{D^2U^N(T-t,\widetilde{X}(t))S_{N,\tau}P_NDG(\widetilde{X}(s))A_NS^{k}_{N,\tau}P_N\xi Q^{\frac{1}{2}}\\
	&\quad\quad\quad\quad(P_NG(\widetilde{X}(t))Q^{\frac{1}{2}})^*\}|dsdt\\
	\le&\ \frac{1}{2}\mathbb{E}\int_{t_k}^{t_{k+1}}\!\!\!\int_{t_k}^t\|D^2U^N(T-t,\widetilde{X}(t))\|_{\mathcal{L}(H)}\|S_{N,\tau}P_N\|_{\mathcal{L}(H)}\\
	&\quad\cdot\|DG(\widetilde{X}(s))A_N^{\frac{2-\gamma}{2}}S^{k}_{N,\tau}A_N^{\frac{\gamma}{2}}P_N\xi Q^{\frac{1}{2}}\|_{\mathcal{L}_2(H)}\|P_NG(\widetilde{X}(t))Q^{\frac{1}{2}}\|_{\mathcal{L}_2(H)}dsdt\\
	\le&\ \frac{1}{2}\mathbb{E}\int_{t_k}^{t_{k+1}}\!\!\!\int_{t_k}^tC(1+\|\widetilde{X}(t)\|^{12}_{\gamma})\cdot C\cdot C\|A_N^{\frac{2-\gamma}{2}}S^{k}_{N,\tau}A_N^{\frac{\gamma}{2}}P_N\xi\|\cdot\|P_NG(\widetilde{X}(t))\|_{\mathcal{L}_2^0}dsdt\\
	\le&\ \frac{1}{2}C\int_{t_k}^{t_{k+1}}\!\!\!\mathbb{E}\int_{t_k}^t(1+\|\widetilde{X}(t)\|^{12}_{\gamma})\cdot C\cdot Ct_{k}^{-\frac{2-\gamma}{2}}\|\xi\|_{\gamma}\cdot C(1+\|\widetilde{X}(t)\|)dsdt\\
	\le&\ Ct_{k+1}^{-\frac{2-\gamma}{2}}\tau^2,
\end{align*}
this together with \eqref{eq3.3.15} implies $\displaystyle \sum_{k=0}^{K-1}c^{1,\mathcal{A}_1}_k \le C\tau$.

{\bf (c.1.1.2) Estimation of $c_k^{1,\mathcal{A}_2}$}. Similarly to {\bf (c.1.1.1)} and by \eqref{eq3.2.5}, we have
\begin{align*}
	&\ c^{1,\mathcal{A}_2}_k\\
	\le&\ \frac{1}{2}\mathbb{E}\int_{t_k}^{t_{k+1}}\int_{t_k}^t|\mathrm{Tr}\{D^2U^N(T-t,\widetilde{X}(t))\\
	&\quad S_{N,\tau}P_NDG(\widetilde{X}(s))A_N\tau\sum_{i=1}^{k-1}S^{k+1-i}_{N,\tau}P_NF(X_{i+1}^N) Q^{\frac{1}{2}}(P_NG(\widetilde{X}(t))Q^{\frac{1}{2}})^*\}|dsdt\\
	\le&\ \frac{1}{2}\mathbb{E}\int_{t_k}^{t_{k+1}}\int_{t_k}^t\|D^2U^N(T-t,\widetilde{X}(t))\|_{\mathcal{L}(H)}\cdot\|S_{N,\tau}P_N\|_{\mathcal{L}(H)}\\
	&\quad\cdot\|DG(\widetilde{X}(s))A_N\tau\sum_{i=1}^{k-1}S^{k+1-i}_{N,\tau}P_NF(X_{i+1}^N)\|_{\mathcal{L}_2^0}\cdot\|P_NG(\widetilde{X}(t))\|_{\mathcal{L}_2^0}dsdt\\
	\le&\ \frac{1}{2}\mathbb{E}\int_{t_k}^{t_{k+1}}\int_{t_k}^tC(1+\|\widetilde{X}(t)\|^{12}_{\gamma})\cdot C\cdot C\|\tau\sum_{i=1}^{k-1}A_N^{1-\varepsilon}S^{k+1-i}_{N,\tau}A_N^{\varepsilon}P_NF(X_{i+1}^N)\| \cdot C(1+\|\widetilde{X}(t)\|)dsdt\\
	\le&\ C\mathbb{E}\int_{t_k}^{t_{k+1}}\int_{t_k}^t(1+\|\widetilde{X}(t)\|^{12}_{\gamma})\cdot\tau\sum_{i=1}^{k-1}t_{k+1-i}^{-1+\varepsilon}\cdot\|A_N^{\varepsilon}P_NF(X_{i+1}^N)\|\cdot(1+\|\widetilde{X}(t)\|)dsdt\\
	\le&\ C\int_{t_k}^{t_{k+1}}\int_{t_k}^t\tau\sum_{i=1}^{k-1}t_{k+1-i}^{-1+\varepsilon}\mathbb{E}\Big[(1+\|\widetilde{X}(t)\|^{12}_{\gamma})\cdot C\|F(X_{i+1}^N)\|_\gamma\cdot(1+\|\widetilde{X}(t)\|)\Big]dsdt\\
	\le&\ C\int_{t_k}^{t_{k+1}}\int_{t_k}^t\tau\sum_{i=1}^{k-1}t_{k+1-i}^{-1+\varepsilon}\cdot Cdsdt \le C\tau^2,
\end{align*}
which implies $\displaystyle \sum_{k=0}^{K-1}c^{1,\mathcal{A}_2}_k \le C\tau$.

{\bf (c.1.1.3) Estimation of $c_k^{1,\mathcal{A}_3}$}. The estimate for $c^{1,\mathcal{A}_3}_k$ is complicated since there includes an infinite dimensional stochastic integral inside the trace term. To address this, we introduce some basic tools. Let $\Gamma\in\mathcal{L}(H)$, denote by $\Gamma^{l,j}=(\Gamma\varphi_l,\varphi_j)$ and then we have
\begin{equation}\label{eq3.6.7}
	\mathrm{Tr}\{ABC\}=\sum_{l=1}^{\infty}(ABC)^{l,l} = \sum_{l,j,m=1}^{\infty}A^{l,j}B^{j,m}C^{m,l}.
\end{equation}
Meanwhile, the  Wiener process $W(t)$ has a spectral expansion in $L^2(\Omega;H)$
\begin{equation}\label{eq3.6.8}
	W(t)=\sum_{n=1}^{\infty}\beta_n(t)Q^{\frac{1}{2}}\varphi_n, \quad a.s.
\end{equation}
where $\beta_n(t)$ is a family of i.i.d.\ one dimensional Wiener processes. With the help of these terminologies, we have 
\begin{align*}
	c^{1,\mathcal{A}_3}_k
	=&\ \frac{1}{2}|\mathbb{E}\int_{t_k}^{t_{k+1}}\!\!\!\!\mathrm{Tr}\{D^2U^N(T-t,\widetilde{X}(t))S_{N,\tau}P_N\mathcal{A}_3(P_NG(\widetilde{X}(t))Q)^*\}dt|\\
	=&\ \frac{1}{2}|\mathbb{E}\int_{t_k}^{t_{k+1}}\!\!\!\!\mathrm{Tr}\{QG^*(\widetilde{X}(t))P_ND^2U^N(T-t,\widetilde{X}(t))S_{N,\tau}P_N\mathcal{A}_3\}dt|\\
	=&\ \frac{1}{2}|\mathbb{E}\int_{t_k}^{t_{k+1}}\!\!\!\!\sum_{l,j,m=1}^{\infty}\Big(QG^*(\widetilde{X}(t))\Big)^{l,j} \Big(P_ND^2U^N(T-t,\widetilde{X}(t))S_{N,\tau}P_N\Big)^{j,m}\mathcal{A}_3^{m,l}dt|\\
	=&\ \frac{1}{2}|\mathbb{E}\int_{t_k}^{t_{k+1}}\int_{t_k}^t\sum_{l,j,m=1}^{\infty}\Big(QG^*(\widetilde{X}(t))\Big)^{l,j} \Big(P_ND^2U^N(T-t,\widetilde{X}(t))S_{N,\tau}P_N\Big)^{j,m}\\
	&\quad\cdot\Big(DG(\widetilde{X}(s))A_NS_{N,\tau}\sum_{i=1}^{k-1}\int_{t_i}^{t_{i+1}}S^{k-i}_{N,\tau}P_NG(X_{i}^N)dW(r)\Big)^{m,l}dsdt|\\
	=&\ \frac{1}{2}|\mathbb{E}\int_{t_k}^{t_{k+1}}\int_{t_k}^t\sum_{l,j,m,n,p=1}^{\infty} \big(QG^*(\widetilde{X}(t))\big)^{l,j}\Big(P_ND^2U^N(T-t,\widetilde{X}(t))S_{N,\tau}P_N\Big)^{j,m}\\
	&\quad\cdot\partial_pG^{m,l}(\widetilde{X}(s))\Big(A_NS_{N,\tau}\sum_{i=1}^{k-1}\int_{t_i}^{t_{i+1}}S^{k-i}_{N,\tau}P_NG(X_{i}^N)Q^{\frac{1}{2}}\varphi_n,\varphi_p\Big)d\beta_n(r)dsdt|.
\end{align*}
Denote by $\mathcal{D}_s^n\widetilde{X}(t) = \mathcal{D}_s\widetilde{X}(t)Q^{\frac{1}{2}}\varphi_n$ $(n\ge 1)$ and apply Malliavin integration by parts formula in Lemma \ref{l3.2.4} of its expansion form, then we obtain
\begin{align*}
	c^{1,\mathcal{A}_3}_k
	=&\ \frac{1}{2}|\mathbb{E}\int_{t_k}^{t_{k+1}}\int_{t_k}^t\sum_{l,j,m,n,p=1}^{\infty} \big(QDG^*(\widetilde{X}(t))\mathcal{D}_r^n\widetilde{X}(t)\big)^{l,j}\Big(P_ND^2U^N(T-t,\widetilde{X}(t))S_{N,\tau}P_N\Big)^{j,m}\\
	&\quad\quad\quad\quad\cdot\partial_pG^{m,l}(\widetilde{X}(s))\big(A_NS_{N,\tau}\sum_{i=1}^{k-1}\int_{t_i}^{t_{i+1}}S^{k-i}_{N,\tau}P_NG(X_{i}^N)Q^{\frac{1}{2}}\varphi_n,\varphi_p\big)drdsdt|\\
	&\ + \frac{1}{2}|\mathbb{E}\int_{t_k}^{t_{k+1}}\int_{t_k}^t\sum_{l,j,m,n,p=1}^{\infty} \big(QG^*(\widetilde{X}(t))\big)^{l,j}\Big(P_ND^3U^N(T-t,\widetilde{X}(t))\mathcal{D}_r^n\widetilde{X}(t)S_{N,\tau}P_N\Big)^{j,m}\\
	&\quad\quad\quad\quad\cdot\partial_pG^{m,l}(\widetilde{X}(s))\big(A_NS_{N,\tau}\sum_{i=1}^{k-1}\int_{t_i}^{t_{i+1}}S^{k-i}_{N,\tau}P_NG(X_{i}^N)Q^{\frac{1}{2}}\varphi_n,\varphi_p\big)drdsdt|\\
	&\ + \frac{1}{2}|\mathbb{E}\int_{t_k}^{t_{k+1}}\int_{t_k}^t\sum_{l,j,m,n,p,q=1}^{\infty} \big(QG^*(\widetilde{X}(t))\big)^{l,j}\Big(P_ND^2U^N(T-t,\widetilde{X}(t))S_{N,\tau}P_N\Big)^{j,m}\\
	& \cdot\partial_{p,q}G^{m,l}(\widetilde{X}(s))(\mathcal{D}_r^n\widetilde{X}(t),\varphi_q)\big(A_NS_{N,\tau}\sum_{i=1}^{k-1}\int_{t_i}^{t_{i+1}}S^{k-i}_{N,\tau}P_NG(X_{i}^N)Q^{\frac{1}{2}}\varphi_n,\varphi_p\big)drdsdt|\\
	:=&\ c^{1,\mathcal{A}_3,1}_k + c^{1,\mathcal{A}_3,2}_k + c^{1,\mathcal{A}_3,3}_k.
\end{align*}
Again, by the definition of trace, we have
\begin{align*}
	&\ c^{1,\mathcal{A}_3,1}_k\\
	=& \frac{1}{2}|\mathbb{E}\int_{t_k}^{t_{k+1}}\int_{t_k}^t\int_{t_i}^{t_{i+1}}\sum_{l,j,m,n=1}^{\infty} \big(QDG^*(\widetilde{X}(t))\mathcal{D}_r^n\widetilde{X}(t)\big)^{l,j}\Big(P_ND^2U^N(T-t,\widetilde{X}(t))S_{N,\tau}P_N\Big)^{j,m}\\
	&\quad\quad\quad\quad \cdot DG^{m,l}(\widetilde{X}(s))A_NS_{N,\tau}\sum_{i=1}^{k-1}S^{k-i}_{N,\tau}P_NG(X_{i}^N)Q^{\frac{1}{2}}\varphi_ndrdsdt|\\
	=& \frac{1}{2}|\mathbb{E}\int_{t_k}^{t_{k+1}}\int_{t_k}^t\int_{t_i}^{t_{i+1}}\sum_{n=1}^{\infty}\mathrm{Tr}\Big\{QDG^*(\widetilde{X}(t))\mathcal{D}_r^n\widetilde{X}(t) P_ND^2U^N(T-t,\widetilde{X}(t))S_{N,\tau}P_N\\
	&\quad\quad\quad\quad \cdot DG(\widetilde{X}(s))A_NS_{N,\tau}\sum_{i=1}^{k-1}S^{k-i}_{N,\tau}P_NG(X_{i}^N)Q^{\frac{1}{2}}\varphi_n\Big\}drdsdt|.
\end{align*}
Then by \eqref{eq7.2.3}, \eqref{eq3.3.4}, Lemma \ref{l3.2.1}, Theorem \ref{t3.4.4}, H\"older inequality, Lemma \ref{l3.4.2}, \eqref{eq3.3.15} and \eqref{eq3.3.14'r}, we have
\begin{align*}
	&\ c^{1,\mathcal{A}_3,1}_k\\ 
	\le&\ \frac{1}{2}\mathbb{E}\int_{t_k}^{t_{k+1}}\int_{t_k}^t\int_{t_i}^{t_{i+1}}\sum_{n=1}^{\infty}\big|\mathrm{Tr}\Big\{ P_ND^2U^N(T-t,\widetilde{X}(t))S_{N,\tau}P_N\\
	&\quad \cdot DG(\widetilde{X}(s))A_NS_{N,\tau}\sum_{i=1}^{k-1}S^{k-i}_{N,\tau}P_NG(X_{i}^N)Q^{\frac{1}{2}}\varphi_nQ^{\frac{1}{2}}\cdot Q^{\frac{1}{2}}DG^*(\widetilde{X}(t))\mathcal{D}_r^n\widetilde{X}(t)\Big\}\big|drdsdt\\
	\le&\ \frac{1}{2}\mathbb{E}\int_{t_k}^{t_{k+1}}\int_{t_k}^t\int_{t_i}^{t_{i+1}}\sum_{n=1}^{\infty}\|P_ND^2U^N(T-t,\widetilde{X}(t))S_{N,\tau}P_N\|_{\mathcal{L}(H)}\\
	&\quad \cdot\|DG(\widetilde{X}(s))A_NS_{N,\tau}\sum_{i=1}^{k-1}S^{k-i}_{N,\tau}P_NG(X_{i}^N)Q^{\frac{1}{2}}\varphi_n\|_{\mathcal{L}_2^0}\cdot\|DG(\widetilde{X}(t))\mathcal{D}_r^n\widetilde{X}(t)\|_{\mathcal{L}_2^0}drdsdt\\
	\le&\ \frac{1}{2}\mathbb{E}\int_{t_k}^{t_{k+1}}\int_{t_k}^t\int_{t_i}^{t_{i+1}}\sum_{n=1}^{\infty} C(1+\|\widetilde{X}(t)\|_{\gamma}^{12})\\
	&\quad\cdot C\|\sum_{i=1}^{k-1}A_N^{1-\varepsilon}S^{k+1-i}_{N,\tau}P_NA_N^{\varepsilon}G(X_{i}^N)Q^{\frac{1}{2}}\varphi_n\|\cdot C\|\mathcal{D}_r^n\widetilde{X}(t)\|drdsdt\\
	\le&\ C\mathbb{E}\int_{t_k}^{t_{k+1}}\int_{t_i}^{t_{i+1}}\sum_{n=1}^{\infty} (1+\|\widetilde{X}(t)\|_{\gamma}^{12})\cdot C\tau\sum_{i=1}^{k-1}t_{k+1-i}^{-1+\varepsilon}\|A_N^{\varepsilon}G(X_{i}^N)Q^{\frac{1}{2}}\varphi_n\|\cdot \|\mathcal{D}_r\widetilde{X}(t)Q^{\frac{1}{2}}\varphi_n\|drdt\\
	\le&\ C\int_{t_k}^{t_{k+1}}\int_{t_k}^t\tau\sum_{i=1}^{k-1}t_{k+1-i}^{-1+\varepsilon} \mathbb{E}\Big[(1+\|\widetilde{X}(t)\|_{\gamma}^{12})\cdot \sum_{n=1}^{\infty}C\|A^{\frac{\gamma}{2}}G(X_{i}^N)Q^{\frac{1}{2}}\varphi_n\|\cdot \|\mathcal{D}_r\widetilde{X}(t)Q^{\frac{1}{2}}\varphi_n\|\Big]drdt\\
	\le&\ C\int_{t_k}^{t_{k+1}}\int_{t_k}^t\tau\sum_{i=1}^{k-1}t_{k+1-i}^{-1+\varepsilon}\cdot \mathbb{E}\Big[(1+\|\widetilde{X}(t)\|_{\gamma}^{12}) \big(\sum_{n=1}^{\infty}\|A^{\frac{\gamma}{2}}G(X_{i}^N)Q^{\frac{1}{2}}\varphi_n\|^2\big)^{\frac{1}{2}}\cdot\\
	&\quad\cdot \big(\sum_{n=1}^{\infty}\|\mathcal{D}_r\widetilde{X}(t)Q^{\frac{1}{2}}\varphi_n\|^2\big)^{\frac{1}{2}}\Big]drdt\\
	=&\ C\int_{t_k}^{t_{k+1}}\int_{t_k}^t\tau\sum_{i=1}^{k-1}t_{k+1-i}^{-1+\varepsilon}  \mathbb{E}\Big[(1+\|\widetilde{X}(t)\|_{\gamma}^{12})\cdot \|G(X_{i}^N)\|_{\mathcal{L}_2^\gamma}\cdot \|\mathcal{D}_r\widetilde{X}(t)\|_{\mathcal{L}_2^0}\Big]drdt\\
	\le&\ C\tau^2.
\end{align*}

Recalling the definition of $\mathrm{Tr}(T) = \sum_{k=1}^{\infty}(Te_k,e_k)_U$, by \eqref{eq7.2.3}, \eqref{eq3.3.4}, Lemma \ref{l3.2.1}, Proposition \ref{p3.4.5}, H\"older inequality, Lemma \ref{l3.4.2}, \eqref{eq3.3.15} and \eqref{eq3.3.14'r}, we get
\begin{align*}
	&\ c^{1,\mathcal{A}_3,2}\\
	=&\ \frac{1}{2}|\mathbb{E}\int_{t_k}^{t_{k+1}}\int_{t_k}^t\sum_{l,j,m,n=1}^{\infty} \big(QG^*(\widetilde{X}(t))\big)^{l,j}\Big(P_ND^3U^N(T-t,\widetilde{X}(t))\mathcal{D}_r^n\widetilde{X}(t)S_{N,\tau}P_N\Big)^{j,m}\\
	&\quad\quad\quad\quad\cdot DG^{m,l}(\widetilde{X}(s))A_NS_{N,\tau}\sum_{i=1}^{k-1}\int_{t_i}^{t_{i+1}}S^{k-i}_{N,\tau}P_NG(X_{i}^N)Q^{\frac{1}{2}}\varphi_ndrdsdt|\\
	=&\ \frac{1}{2}|\mathbb{E}\int_{t_k}^{t_{k+1}}\int_{t_k}^t\sum_{n=1}^{\infty}\mathrm{Tr} \Big\{P_ND^3U^N(T-t,\widetilde{X}(t))A^{-\frac{1}{4}-\delta}A^{\frac{1}{4}+\delta}\mathcal{D}_r^n\widetilde{X}(t)S_{N,\tau}P_N\\
	&\quad\quad\quad\quad\cdot DG(\widetilde{X}(s))A_NS_{N,\tau}\sum_{i=1}^{k-1}\int_{t_i}^{t_{i+1}}S^{k-i}_{N,\tau}P_NG(X_{i}^N)Q^{\frac{1}{2}}\varphi_nQ^{\frac{1}{2}}Q^{\frac{1}{2}}G^*(\widetilde{X}(t))\Big\}drdsdt|\\
	\le&\ \frac{1}{2}\mathbb{E}\int_{t_k}^{t_{k+1}}\int_{t_k}^t\sum_{n=1}^{\infty} C(1+\|\widetilde{X}(t)\|_{\gamma}^{24})\cdot\|A^{\frac{1}{4}+\delta}\mathcal{D}_r\widetilde{X}(t)Q^{\frac{1}{2}}\varphi_n\|\cdot\|S_{N,\tau}P_N\|_{\mathcal{L}(H)}\\
	&\quad\quad\quad\quad\cdot \|DG(\widetilde{X}(s))A_NS_{N,\tau}\sum_{i=1}^{k-1}\int_{t_i}^{t_{i+1}}S^{k-i}_{N,\tau}P_NG(X_{i}^N)Q^{\frac{1}{2}}\varphi_n\|_{\mathcal{L}_2^0}\cdot\|G(\widetilde{X}(t))\|_{\mathcal{L}_2^0}drdsdt\\
	\le&\ C\int_{t_k}^{t_{k+1}}\int_{t_k}^t\tau\sum_{i=1}^{k-1}t_{k+1-i}^{-1+\varepsilon}  \mathbb{E}\Big[(1+\|\widetilde{X}(t)\|_{\gamma}^{24})\cdot\|A^{\frac{\gamma}{2}}\mathcal{D}_r\widetilde{X}(t)\|_{\mathcal{L}_2^0}\\
	&\quad\cdot\|A^{\frac{\gamma}{2}}G(X_{i}^N)\|_{\mathcal{L}_2^0}\cdot\|G(\widetilde{X}(t))\|_{\mathcal{L}_2^0}\Big]drdsdt\\
	\le&\ C\tau^2.
\end{align*}
In a similar way, we obtain
\begin{align*}
	&\ c^{1,\mathcal{A}_3,3}\\
	=&\ \frac{1}{2}|\mathbb{E}\int_{t_k}^{t_{k+1}}\int_{t_k}^t\sum_{l,j,m,n=1}^{\infty} \big(QG^*(\widetilde{X}(t))\big)^{l,j}\Big(P_ND^2U^N(T-t,\widetilde{X}(t))S_{N,\tau}P_N\Big)^{j,m}\\
	& \quad\cdot D^2G^{m,l}(\widetilde{X}(s))\mathcal{D}_r^n\widetilde{X}(t)A_NS_{N,\tau}\sum_{i=1}^{k-1}\int_{t_i}^{t_{i+1}}S^{k-i}_{N,\tau}P_NG(X_{i}^N)Q^{\frac{1}{2}}\varphi_ndrdsdt|\\
	=&\ \frac{1}{2}|\mathbb{E}\int_{t_k}^{t_{k+1}}\int_{t_k}^t\sum_{n=1}^{\infty}\mathrm{Tr} \Big\{P_ND^2U^N(T-t,\widetilde{X}(t))S_{N,\tau}P_N\\
	& \quad\cdot D^2G(\widetilde{X}(s))\mathcal{D}_r^n\widetilde{X}(t)A_NS_{N,\tau}\sum_{i=1}^{k-1}\int_{t_i}^{t_{i+1}}S^{k-i}_{N,\tau}P_NG(X_{i}^N)Q^{\frac{1}{2}}\varphi_nQ^{\frac{1}{2}}Q^{\frac{1}{2}}G^*(\widetilde{X}(t))\Big\}drdsdt|\\
	\le&\ \frac{1}{2}\mathbb{E}\int_{t_k}^{t_{k+1}}\int_{t_k}^t\sum_{n=1}^{\infty}C(1+\|\widetilde{X}(t)\|_{\gamma}^{12})\cdot\|S_{N,\tau}P_N\|_{\mathcal{L}(H)}\\
	& \quad\cdot \|D^2G(\widetilde{X}(s))\mathcal{D}_r^n\widetilde{X}(t)A_NS_{N,\tau}\sum_{i=1}^{k-1}\int_{t_i}^{t_{i+1}}S^{k-i}_{N,\tau}P_NG(X_{i}^N)Q^{\frac{1}{2}}\varphi_n\|_{\mathcal{L}_2^0}\|G(\widetilde{X}(t))\|_{\mathcal{L}_2^0}drdsdt\\
	\le&\ C\tau^2.
\end{align*}
Thus we have proved $c^{1,\mathcal{A}_3}_k\le C\tau^2$.
According to the estimates in {\bf (c.1.1.1)}-{\bf (c.1.1.3)}, we obtain
\begin{equation*}
	\sum_{k=0}^{K-1}c^{1,\mathcal{A}}_k\le C\tau.
\end{equation*} 

{\bf (c.1.2) Estimation of $c_k^{1,\mathcal{B}}$}. By \eqref{eq3.2.5}, Lemma \ref{l3.2.1} and Theorem \ref{t3.4.4}, we obtain
\begin{align*}
	c_k^{1,\mathcal{B}} \le&\ \frac{1}{2}\mathbb{E}\int_{t_k}^{t_{k+1}}\int_{t_k}^t|\mathrm{Tr}\{D^2U^N(T-t,\widetilde{X}(t))S_{N,\tau}P_N\\
	&\quad\quad\cdot DG(\widetilde{X}(s))S_{N,\tau}P_NF(X_{k+1}^N)Q^{\frac{1}{2}}(P_NG(\widetilde{X}(t))Q^{\frac{1}{2}})^*\}|dsdt\\
	\le&\ \frac{1}{2}\mathbb{E}\int_{t_k}^{t_{k+1}}\int_{t_k}^t\|D^2U^N(T-t,\widetilde{X}(t))\|_{\mathcal{L}(H)}\|S_{N,\tau}P_N\|_{\mathcal{L}(H)}\\
	&\quad\quad\cdot \|DG(\widetilde{X}(s))S_{N,\tau}P_NF(X_{k+1}^N)\|_{\mathcal{L}_2^0}\cdot\|P_NG(\widetilde{X}(t))\|_{\mathcal{L}_2^0}dsdt\\
	\le&\ C\tau^2,
\end{align*}
which implies $\displaystyle \sum_{k=0}^{K-1}c_k^{1,\mathcal{B}} \le C\tau$.

{\bf (c.1.3) Estimation of $c_k^{1,\mathcal{C}}$}. 
By \eqref{eq3.6.7}, \eqref{eq3.6.8} and Malliavin integration by parts formula in Lemma \ref{l3.2.4} of its expansion form, we get
\begin{align*}
	&\ c^{1,\mathcal{C}}_k\\
	=&\ \frac{1}{2}|\mathbb{E}\int_{t_k}^{t_{k+1}}\mathrm{Tr}\{QG^*(\widetilde{X}(t))P_ND^2U^N(T-t,\widetilde{X}(t))S_{N,\tau}P_N\mathcal{C} \}dt|\\
	=&\ \frac{1}{2}|\mathbb{E}\int_{t_k}^{t_{k+1}}\sum_{i,j,m=1}^{\infty} \Big(QG^*(\widetilde{X}(t))\Big)^{i,j}\Big(P_ND^2U^N(T-t,\widetilde{X}(t))S_{N,\tau}P_N\Big)^{j,m}\\
	&\quad\quad\quad\quad\cdot\Big(\int_{t_k}^tDG(\widetilde{X}(s))S_{N,\tau}P_NG(X_{k}^N)dW(s)\Big)^{m,i}dt|\\
	=&\ \frac{1}{2}|\int_{t_k}^{t_{k+1}}\sum_{i,j,m,n,p=1}^{\infty} \mathbb{E}\Big[\big(QG^*(\widetilde{X}(t))\big)^{i,j}\cdot\Big(P_ND^2U^N(T-t,\widetilde{X}(t))S_{N,\tau}P_N\Big)^{j,m}\\
	&\quad\quad\quad\quad\cdot\int_{t_k}^t\partial_p G^{m,i}(\widetilde{X}(s))\big(S_{N,\tau}P_NG(X_{k}^N)Q^{\frac{1}{2}}\varphi_n,\varphi_p\big)d\beta_n(s)\Big]dt|\\
	\le&\ \frac{1}{2}|\int_{t_k}^{t_{k+1}}\!\!\!\!\int_{t_k}^t\sum_{i,j,m,n,p=1}^{\infty}\!\!\!\! \mathbb{E}\Big[\big(QDG^*(\widetilde{X}(t))\mathcal{D}_s^n\widetilde{X}(t)\big)^{i,j}\Big(P_ND^2U^N(T-t,\widetilde{X}(t))S_{N,\tau}P_N\Big)^{j,m}\\
	&\quad\quad\quad\quad\partial_p G^{m,i}(\widetilde{X}(s))\big(S_{N,\tau}P_NG(X_{k}^N)Q^{\frac{1}{2}}\varphi_n,\varphi_p\big)ds\Big]dt|\\
	&\ + \frac{1}{2}|\int_{t_k}^{t_{k+1}}\!\!\!\!\int_{t_k}^t\sum_{i,j,m,n,p=1}^{\infty}\!\!\!\! \mathbb{E}\Big[\big(QG^*(\widetilde{X}(t))\big)^{i,j}\Big(P_ND^3U^N(T\!-\!t,\widetilde{X}(t))\mathcal{D}_s^n\widetilde{X}(t)S_{N,\tau}P_N\Big)^{j,m}\\
	&\quad\quad\quad\quad\partial_p G^{m,i}(\widetilde{X}(s))\big(S_{N,\tau}P_NG(X_{k}^N)Q^{\frac{1}{2}}\varphi_n,\varphi_p\big)ds\Big]dt|\\
	=&\ \frac{1}{2}|\int_{t_k}^{t_{k+1}}\!\!\int_{t_k}^t\sum_{n=1}^{\infty} \mathbb{E}\ \mathrm{Tr}\Big\{QDG^*(\widetilde{X}(t))\mathcal{D}_s^n\widetilde{X}(t)\cdot P_ND^2U^N(T-t,\widetilde{X}(t))S_{N,\tau}P_N\\
	&\quad\quad\quad\quad \cdot DG(\widetilde{X}(s))S_{N,\tau}P_NG(X_{k}^N)Q^{\frac{1}{2}}\varphi_n\Big\}dsdt|\\
	&\ + \frac{1}{2}|\int_{t_k}^{t_{k+1}}\!\!\int_{t_k}^t\sum_{n=1}^{\infty} \mathbb{E}\ \mathrm{Tr}\Big\{QG^*(\widetilde{X}(t))\cdot P_ND^3U^N(T-t,\widetilde{X}(t))\mathcal{D}_s^n\widetilde{X}(t)S_{N,\tau}P_N\\
	&\quad\quad\quad\quad\cdot DG(\widetilde{X}(s))S_{N,\tau}P_NG(X_{k}^N)Q^{\frac{1}{2}}\varphi_n\Big\}dsdt|\\
	\le&\ C\tau^2.
\end{align*}
Here, we have used the same technique as estimating $c_k^{1,\mathcal{A}_3}$.

{\bf (c.1.4) Estimation of $c_k^{1,\mathcal{E}}$}. By the definition of trace and expression of $\mathcal{E}$, we have for $h_1,h_2\in H$,
\begin{align*}
	&\ (\mathcal{E}h_1,h_2) \\
	=&\ \frac{1}{2}\int_{t_k}^t\Big(\mathrm{Tr}\{(S_{N,\tau}P_NG(X_{k}^N)Q^{\frac{1}{2}})^*D^2G(\widetilde{X}(s))S_{N,\tau}P_NG(X_{k}^N)Q^{\frac{1}{2}}\}h_1,h_2\Big)ds\\
	=&\ \frac{1}{2}\int_{t_k}^t\sum_{i=1}^{\infty}\Big(\big(D^2G(\widetilde{X}(s))S_{N,\tau}P_NG(X_{k}^N)Q^{\frac{1}{2}}\varphi_i,S_{N,\tau}P_NG(X_{k}^N)Q^{\frac{1}{2}}\varphi_i\big)h_1,h_2\Big)ds.
\end{align*}
By Lemma \ref{l3.2.1} and \eqref{eq3.3.4}, we have
\begin{align*}
	\|\mathcal{E}\|_{\mathcal{L}_2^0}
	\le&\ \frac{1}{2}\int_{t_k}^t\sum_{i=1}^{\infty}\big|\big(D^2G(\widetilde{X}(s))S_{N,\tau}P_NG(X_{k}^N)Q^{\frac{1}{2}}\varphi_i,S_{N,\tau}P_NG(X_{k}^N)Q^{\frac{1}{2}}\varphi_i\big)\big|ds\\
	\le&\ \frac{1}{2}\int_{t_k}^t\sum_{i=1}^{\infty}C\|S_{N,\tau}P_NG(X_{k}^N)Q^{\frac{1}{2}}\varphi_i\|^2ds\\
	=&\ C\tau\|S_{N,\tau}P_NG(X_{k}^N)\|_{\mathcal{L}_2^0}^2 \le C\tau(1+\|X_k^N\|^2).
\end{align*}
Then by \eqref{eq7.2.3}, \eqref{eq3.3.4} and Theorem \ref{t3.4.4}, we obtain
\begin{align*}
	c^{1,\mathcal{E}}_k
	\le&\  \frac{1}{2}\mathbb{E}\int_{t_k}^{t_{k+1}}|\mathrm{Tr}\{D^2U^N(T-t,\widetilde{X}(t))S_{N,\tau}P_N\mathcal{E} Q^{\frac{1}{2}}(P_NG(\widetilde{X}(t))Q^{\frac{1}{2}})^*\}|dt\\
	\le&\ \frac{1}{2}\mathbb{E}\int_{t_k}^{t_{k+1}}\|D^2U^N(T-t,\widetilde{X}(t))S_{N,\tau}P_N\|_{\mathcal{L}(H)}\|\mathcal{E}\|_{\mathcal{L}_2^0}\|P_NG(\widetilde{X}(t)\|_{\mathcal{L}_2^0}dt\\
	\le&\ \frac{1}{2}\mathbb{E}\int_{t_k}^{t_{k+1}}C(1+\|\widetilde{X}(t)\|^{12}_{\gamma})\cdot C\cdot C\tau(1+\|X_k^N\|^2)\cdot C \cdot C(1+\|\widetilde{X}(t)\|)dt\\
	\le&\ C\tau^2.
\end{align*}

{\bf (c.1.5) Estimation of $c_k^{1,\mathcal{G}}$}. 
By the definition of trace and expression of $\mathcal{G}$, we have for $h_1,h_2\in H$,
\begin{align*}
	&\ (\mathcal{G}h_1,h_2) \\
	=&\ \int_{t_k}^t\Big(\mathrm{Tr}\{(S_{N,\tau}P_NG(X_{k}^N)Q^{\frac{1}{2}})^*D^2G(\widetilde{X}(s))(\mathcal{D}^-\widetilde{X})_sQ^{\frac{1}{2}}h_1,h_2\Big)ds\\
	=&\ \int_{t_k}^t\sum_{i=1}^{\infty}\Big(\big(D^2G(\widetilde{X}(s))(\mathcal{D}^-\widetilde{X})_sQ^{\frac{1}{2}}\varphi_i,S_{N,\tau}P_NG(X_{k}^N)Q^{\frac{1}{2}}\varphi_i\big)h_1,h_2\Big)ds.
\end{align*}
By Lemma \ref{l3.2.1}, \eqref{eq3.3.4} and H\"older inequality, we have
\begin{align*}
	\|\mathcal{E}\|_{\mathcal{L}_2^0}
	\le&\ \int_{t_k}^t\sum_{i=1}^{\infty}\big|\big(D^2G(\widetilde{X}(s))(\mathcal{D}^-\tilde{X})_sQ^{\frac{1}{2}}\varphi_i,S_{N,\tau}P_NG(X_{k}^N)Q^{\frac{1}{2}}\varphi_i\big)\big|ds\\
	\le&\ \int_{t_k}^t\sum_{i=1}^{\infty}C\|(\mathcal{D}^-\tilde{X})_sQ^{\frac{1}{2}}\varphi_i\|\cdot\|S_{N,\tau}P_NG(X_{k}^N)Q^{\frac{1}{2}}\varphi_i\|ds\\
	\le&\ C\tau\sup_{r\in[0,T]}\big(\sum_{i=1}^{\infty}\|(\mathcal{D}^-\tilde{X})_rQ^{\frac{1}{2}}\varphi_i\|^2\big)^{\frac{1}{2}}\big(\sum_{i=1}^{\infty}\|S_{N,\tau}P_NG(X_{k}^N)Q^{\frac{1}{2}}\varphi_i\|^2\big)^{\frac{1}{2}}\\
	=&\ C\tau\sup_{r\in[0,T]}\|(\mathcal{D}^-\tilde{X})_r\|_{\mathcal{L}_2^0}\cdot \|S_{N,\tau}P_NG(X_{k}^N)\|_{\mathcal{L}_2^0}\\ 
	\le&\ C\tau\sup_{r\in[0,T]}\|(\mathcal{D}^-\tilde{X})_r\|_{\mathcal{L}_2^0}\cdot(1+\|X_k^N\|).
\end{align*}
Then by \eqref{eq7.2.3}, \eqref{eq3.3.4}, Theorem \ref{t3.4.4} and Lemma \ref{t3.4.1}, we obtain
\begin{align*}
	c^{1,\mathcal{G}}_k
	\le&\  \frac{1}{2}\mathbb{E}\int_{t_k}^{t_{k+1}}|\mathrm{Tr}\{D^2U^N(T-t,\widetilde{X}(t))S_{N,\tau}P_N\mathcal{G} Q^{\frac{1}{2}}(P_NG(\widetilde{X}(t))Q^{\frac{1}{2}})^*\}|dt\\
	\le&\ \frac{1}{2}\mathbb{E}\int_{t_k}^{t_{k+1}}\|D^2U^N(T-t,\widetilde{X}(t))S_{N,\tau}P_N\|_{\mathcal{L}(H)}\|\mathcal{G}\|_{\mathcal{L}_2^0}\|P_NG(\widetilde{X}(t))\|_{\mathcal{L}_2^0}dt\\
	\le&\ \frac{1}{2}\mathbb{E}\int_{t_k}^{t_{k+1}}C(1+\|\widetilde{X}(t)\|^{12}_{\gamma})\cdot C\cdot C\tau\sup_{r\in[0,T]}\|(\mathcal{D}^-\tilde{X})_r\|_{\mathcal{L}_2^0}\cdot(1+\|X_k^N\|)\\
	&\quad\quad\quad \cdot C \cdot C(1+\|\widetilde{X}(t)\|)dt\\
	\le&\ C\tau^2.
\end{align*}

From the estimates in {\bf (c.1.1)}-{\bf (c.1.5)}, it follows that $\displaystyle\sum_{k=0}^{K-1}c_k^1 \le C\tau$. In a similar way, we can prove $\displaystyle\sum_{k=0}^{K-1}c_k^3 \le C\tau$.

{\bf (c.2) Estimation of $c_k^{2}$ and $c_k^{4}$}.
We will only investigate $c_k^{2}$, and the estimate for $c_k^{4}$ follows from a similar way. 

By \eqref{eq3.3.3.1}, \eqref{eq7.2.3}, Theorem \ref{t3.4.4}, Lemma \ref{l3.2.1}, Lemma \ref{l3.4.2} and \eqref{eq3.3.14'r},  we have
\begin{align*}
	c_k^2
	\le&\ \frac{1}{2}\mathbb{E}\int_{t_k}^{t_{k+1}}|\mathrm{Tr}\{D^2U^N(T-t,\widetilde{X}(t))[S_{N,\tau}-I]P_NG(X_k^N)Q^{\frac{1}{2}}(P_NG(\widetilde{X}(t))Q^{\frac{1}{2}})^*\}|dt\\
	=&\ \frac{1}{2}\mathbb{E}\int_{t_k}^{t_{k+1}}|\mathrm{Tr}\{D^2U^N(T-t,\widetilde{X}(t))\tau A_NS_{N,\tau}P_NG(X_k^N)Q^{\frac{1}{2}}(P_NG(\widetilde{X}(t))Q^{\frac{1}{2}})^*\}|dt\\
	\le&\ \frac{\tau}{2}\mathbb{E}\int_{t_k}^{t_{k+1}}\|D^2U^N(T\!-\!t,\widetilde{X}(t))\|_{\mathcal{L}(H)}|\mathrm{Tr}\{A_N^{1-\gamma}S_{N,\tau}A_N^{\frac{\gamma}{2}}P_NG(X_k^N)Q^{\frac{1}{2}}(A_N^{\frac{\gamma}{2}}P_NG(\widetilde{X}(t))Q^{\frac{1}{2}})^*\}|dt\\
	\le&\ \frac{\tau}{2}\mathbb{E}\!\int_{t_k}^{t_{k+1}} C(1+\|\widetilde{X}(t)\|^{12}_{\gamma})\cdot\|A_N^{1-\gamma}S_{N,\tau}\|_{\mathcal{L}(H)}\|A_N^{\frac{\gamma}{2}}S_{N,\tau}P_NG(X_k^N)\|_{\mathcal{L}_2^0}\cdot\|A_N^{\frac{\gamma}{2}}P_NG(\widetilde{X}(t))\|_{\mathcal{L}_2^0}dt\\
	\le&\ C\tau\mathbb{E}\int_{t_k}^{t_{k+1}} (1+\|\widetilde{X}(t)\|^{12}_{\gamma})\cdot \tau^{\gamma-1}\cdot C\|G(X_k^N)\|_{\mathcal{L}_2^\gamma}\cdot C\|G(\widetilde{X}(t))\|_{\mathcal{L}_2^\gamma}dt\\
	\le&\ C\tau^{\gamma+1}.
\end{align*}
In a similar way, we can prove $c_k^4\le C\tau^{\gamma+1}$. According to the estimates in {\bf (c.1)} and {\bf (c.2)}, we have
\begin{equation}\label{eq3.6.c}
	\sum_{k=0}^{K-1}c_k\le C\tau^{\gamma}.
\end{equation}

\vskip1mm
\section{Estimation for  $d_k$ in \eqref{eq3.6.2}}
\vskip1mm

By \eqref{eq7.2.3}, Theorem \ref{t3.4.4}, Lemma \ref{l3.2.1}, Lemma \ref{l3.4.2}, \eqref{eq3.3.14'r} and Lemma \ref{t3.4.1}, we have
\begin{align*}
	d_k \le&\ \mathbb{E}\int_{t_k}^{t_{k+1}}|\mathrm{Tr}\big\{D^2U^N(T-t,\widetilde{X}(t))((\mathcal{D}^-\widetilde{X})_tQ^{\frac{1}{2}})(S_{N,\tau}P_NG(X_k^N)Q^{\frac{1}{2}})^*\big\}|dt\\
	\le&\ \mathbb{E}\int_{t_k}^{t_{k+1}}\|D^2U^N(T-t,\widetilde{X}(t))\|_{\mathcal{L}(H)}\|(\mathcal{D}^-\widetilde{X})_t\|_{\mathcal{L}_2^0}\|S_{N,\tau}P_NG(X_k^N)\|_{\mathcal{L}_2^0}dt\\
	\le&\ \mathbb{E}\int_{t_k}^{t_{k+1}} C(1+\|\widetilde{X}(t)\|_{\gamma}^{12})\cdot\|(\mathcal{D}^-\widetilde{X})_t\|_{\mathcal{L}_2^0}\cdot C\|G(X_k^N)\|_{\mathcal{L}_2^0}dt\\
	\le&\ C\tau^2.
\end{align*}
Thus we obtain 
\begin{equation}
	\sum_{k=0}^{K-1}d_k\le C\tau.
\end{equation}

\end{document}